\documentclass[12pt,a4paper]{amsart}
\pdfoutput=1
\usepackage[utf8]{inputenc}
\usepackage[english,greek,french]{babel}
\usepackage[T1]{fontenc}
\usepackage{charter}

\usepackage{listings}

\usepackage[dvipsnames]{xcolor}
\usepackage{color}              
\usepackage{indentfirst}        
\usepackage{latexsym,bm}        
\usepackage{array}
\usepackage{amsxtra}
\usepackage{amstext}
\usepackage[arrow, matrix, curve]{xy}
\xymatrixrowsep{15mm}
\xymatrixcolsep{15mm}
\usepackage{xypic}
\usepackage{dsfont} 
\usepackage{amsmath,amsthm,amssymb,amscd,amsfonts,typearea}
\usepackage{tikz}
\usepackage{mathrsfs}
\usepackage{mathtools}
\usepackage{faktor}
\usepackage{xfrac}
\usepackage{enumitem}
\usepackage{moreenum}	
\usepackage{tensor}

\usepackage[
    backend=biber,
    sorting=nyt,
    style=ieee
]{biblatex}
\usepackage[colorlinks, linkcolor=red,
anchorcolor=blue, citecolor=blue]{hyperref}

\usepackage{pdfsync}

\usepackage[a4paper,right=2.5cm,left=2.5cm,
top=2.5cm,bottom=2cm]{geometry}

\DeclareMathOperator{\inj}{inj}

\newcommand{\delbar}{\overline{\partial}}
\newcommand{\lb}{\left\{}
\newcommand{\rb}{\right\}}

\newcommand{\abs}[1]{\left\lvert #1 \right\rvert}
\newcommand{\norm}[1]{\| #1 \|}

\newcommand{\ord}{\mathrm{ord}}

\newcommand{\Dom}{\mathrm{Dom}}

\newcommand{\restr}[1]{{\mkern 1mu \vrule height 2ex\mkern2mu #1}}
\newcommand{\supp}[1]{\mathrm{supp}\, #1}

\renewcommand{\d}{\mathrm{d}}
\newcommand{\C}{\mathbb{C}}
\newcommand{\R}{\mathbb{R}}
\newcommand{\Z}{\mathbb{Z}}
\newcommand{\N}{\mathbb{N}}
\newcommand{\D}{\mathbb{D}}

\newcommand{\E}{\mathbb{E}} 

\renewcommand{\P}{\mathbb{P}} 

\newcommand{\dz}{\mathrm{d}z}
\newcommand{\dzb}{\mathrm{d}\bar{z}}

\newcommand{\dx}{\mathrm{d}x}
\newcommand{\dy}{\mathrm{d}y}
\newcommand{\db}{\overline{\partial}}

\newcommand{\ui}{\mathsf{i}}

\newcommand{\cG}{\mathcal{G}}
\newcommand{\cLL}{\mathcal{L}^{2}}

\newcommand{\bN}{\mathbb{N}}
\newcommand{\bB}{\mathbb{B}}
\newcommand{\fL}{\mathscr{L}}

\newcommand{\bfH}{\mathbf{H}}

\newcommand{\bb}{\boldsymbol}

\DeclareMathOperator{\Div}{Div}

\DeclareMathOperator{\dist}{dist}
\newcommand{\cali}[1]{\mathscr{#1}}

\newcommand{\cC}{\cali{C}}

\newcommand{\mL}{\mathcal{L}}

\newcommand{\FS}{{{_\mathrm{FS}}}}

\newcommand{\longdashrightarrow}[2][->]{
\tikz[baseline=-\the\dimexpr\fontdimen22\textfont2\relax]{
\node[anchor=south,font=\scriptsize, inner ysep=1.5pt,outer xsep=2.2pt](x){#2};
\draw[shorten <=3.4pt,shorten >=3.4pt,dashed,#1](x.south west)--(x.south east);
}
}

\makeatletter
\let\save@mathaccent\mathaccent
\newcommand*\if@single[3]{%
\setbox0\hbox{${\mathaccent"0362{#1}}^H$}%
\setbox2\hbox{${\mathaccent"0362{\kern0pt#1}}^H$}%
\ifdim\ht0=\ht2 #3\else #2\fi
}
\newcommand*\rel@kern[1]{\kern#1\dimexpr\macc@kerna}
\newcommand*\widebar[1]{\@ifnextchar^{{\wide@bar{#1}{0}}}{\wide@bar{#1}{1}}}
\newcommand*\wide@bar[2]{\if@single{#1}{\wide@bar@{#1}{#2}{1}}{\wide@bar@{#1}{#2}{2}}}
\newcommand*\wide@bar@[3]{%
\begingroup
\def\mathaccent##1##2{%
\let\mathaccent\save@mathaccent
\if#32 \let\macc@nucleus\first@char \fi
\setbox\z@\hbox{$\macc@style{\macc@nucleus}_{}$}%
\setbox\tw@\hbox{$\macc@style{\macc@nucleus}{}_{}$}%
\dimen@\wd\tw@
\advance\dimen@-\wd\z@
\divide\dimen@ 3
\@tempdima\wd\tw@
\advance\@tempdima-\scriptspace
\divide\@tempdima 10
\advance\dimen@-\@tempdima
\ifdim\dimen@>\z@ \dimen@0pt\fi
\rel@kern{0.6}\kern-\dimen@
\if#31
\overline{\rel@kern{-0.6}\kern\dimen@\macc@nucleus\rel@kern{0.4}\kern\dimen@}%
\advance\dimen@0.4\dimexpr\macc@kerna
\let\final@kern#2%
\ifdim\dimen@<\z@ \let\final@kern1\fi
\if\final@kern1 \kern-\dimen@\fi
\else
\overline{\rel@kern{-0.6}\kern\dimen@#1}%
\fi
}%
\macc@depth\@ne
\let\math@bgroup\@empty \let\math@egroup\macc@set@skewchar
\mathsurround\z@ \frozen@everymath{\mathgroup\macc@group\relax}%
\macc@set@skewchar\relax
\let\mathaccentV\macc@nested@a
\if#31
\macc@nested@a\relax111{#1}%
\else
\def\gobble@till@marker##1\endmarker{}%
\futurelet\first@char\gobble@till@marker#1\endmarker
\ifcat\noexpand\first@char A\else
\def\first@char{}%
\fi
\macc@nested@a\relax111{\first@char}%
\fi
\endgroup
}
\makeatother

\renewcommand{\bar}{\widebar}

\begin{document}
\selectlanguage{english}


\title[Semipositive line bundle on Punctured Riemann surface]{Semipositive line bundles on punctured Riemann surfaces: Bergman kernels and random zeros}
\author{Bingxiao Liu, Dominik Zielinski}
\address{Universit{\"a}t zu K{\"o}ln,  Department Mathematik/Informatik,
    Weyertal 86-90,   50931 K{\"o}ln, Germany}
    \email{bingxiao.liu@uni-koeln.de}
\email{dzielinski.math@posteo.de}
\thanks{B.\ L.\ is supported by the 
DFG Priority Program 2265 `Random Geometric Systems' (Project-ID 422743078).}

\begin{abstract}
We give an extensive study on the Bergman kernel expansions and the random zeros associated with the high tensor powers of a semipositive line bundle on a complete punctured Riemann surface. We prove several results for the zeros of Gaussian holomorphic sections in the semi-classical limit, including the equidistribution, large deviation estimates, central limit theorem, and number variances.
\end{abstract}

\keywords{Riemann surface; Bergman kernel; semipositive line bundle; random zeros; equidistribution; semi-classical limit}

\maketitle

\maketitle

\theoremstyle{plain}
\newtheorem{theorem}{Theorem}[subsection]
\newtheorem{lemma}[theorem]{Lemma}
\newtheorem{proposition}[theorem]{Proposition}
\newtheorem{corollary}[theorem]{Corollary}

\theoremstyle{remark}
\newtheorem{definition}[theorem]{Definition}
\newtheorem{remark}[theorem]{Remark}
\newtheorem{example}[theorem]{Example}



\thispagestyle{empty}
\tableofcontents

\setcounter{page}{1}
\numberwithin{equation}{subsection}
\renewcommand\thesection{\arabic{section}}
\renewcommand\thesubsection{\arabic{section}.\arabic{subsection}}

\section{Introduction}

This paper aims to give an extensive study on the Bergman kernel expansions and the random zeros under the semi-classical limit associated to the high tensor powers of a semi-positively curved (semipositive for short) line bundle on a complete punctured Riemann surface.

The first half part of this paper, including the results for the spectral gap and Bergman kernel expansions, was done in the Ph.D. thesis of the second named author \cite{Thesis-Zielinski}. Then, following the recent work of the first named author with Drewitz and Marinescu \cite{Drewitz_2023, DrLM:2023aa, Drewitz:2024aa}, we applied these results to study the zeros of the Gaussian holomorphic sections for the semipositive line bundles, including equidistribution, large deviation estimates, the central limit theorem, and number variances.

An effective approach for Bergman kernel expansions is the method of analytic localization as explained in detail by Ma and Marinescu in their book \cite{MM07}. A key ingredient in their method is the spectral gap of Kodaira Laplacians that holds for the uniformly positive line bundles on complete Hermitian manifolds (the metrics are always taken to be smooth unless we say otherwise). However, for semipositive line bundles (the Chern curvature form is nonnegative), there are examples (see \cite{MR2023951}) of compact Hermitian manifolds with complex dimension $\geqslant 2$ such that the spectral gap does not hold. For the semipositive line bundles on a compact Riemann surface, a certain spectral gap always holds, provided that the Chern curvature admits at least a strictly positive point. Recently, Marinescu and Savale \cite{Marinescu2023, MS23} worked out precisely the spectral gap by subelliptic estimates for this setting under the assumption that Chern curvature vanishes at most to finite order on the compact Riemann surface. Then they obtained the asymptotic expansions of the Bergman kernel functions, that is, the on-diagonal Bergman kernels. Their result shows that the expansion factors at the vanishing points of the Chern curvature are different from the non-vanishing points. Here, we extend further their work to the case of complete punctured Riemann surfaces and give the results for the near-diagonal expansions of Bergman kernels. Note that, for semipositive or big line bundles with singular metrics on complex manifolds of general dimension, there are also other approaches such as $L^2$-estimates for $\bar{\partial}$-operator to study the Bergman kernels; see \cite{MR2016088, CM15, DMM16}.

The complete punctured Riemann surfaces that are the subject of this paper have already been examined by Auvray, Ma, and Marinescu \cite{AMM16, AMM21, AMM22}, where they give the expansions of Bergman kernels for the high tensor powers of a uniformly positive line bundle under the assumption of Poincar\'{e} metric near the punctures. The important examples for this model of Riemann surfaces are arithmetic surfaces, on which the holomorphic sections correspond to cusp forms (see \cite{AMM21} or \cite[Section 4]{Drewitz_2023}). Note that for positive line bundles on punctured Riemann surfaces equipped with non-smooth metrics, Coman, Klevtsov, and Marinescu \cite{MR3951075} obtained the estimates and the leading term of the Bergman kernel functions and then discussed several interesting applications.

In \cite{Drewitz_2023}, the first named author with Drewitz and Marinescu applied the results from \cite{AMM16, AMM21} to study the zeros of random holomorphic sections for a positive line bundle on the complete punctured Riemann surface. In particular, estimates for large deviations and hole probabilities were established, following the seminal work of Shiffman, Zelditch, and Zrebiec \cite{SZZ08}. In this paper, we investigate the above problems under the semipositive condition; see Theorems \ref{thm:1.3.2-24}, \ref{thm:5.3.1}, and Proposition \ref{prop:1.2.3}. Moreover, we go further to work out the smooth statistics such as number variance and central limit theorem for the random zeros; see Theorems \ref{thm:3.5.1ss} and \ref{thm:1.5.2-24ss}. We will see that the existence of vanishing points of the Chern curvature form requires more techniques in the proofs, but eventually, they will not contribute to the principal behaviors of random zeros. It remains interesting to study the subprincipal behaviors of random zeros to identify the contribution of vanishing points.

The random zeros as point processes on Riemann surfaces provide a valuable model for quantum chaotic dynamics as in \cite{Bogomolny_1996, MR1649013}. In \cite{MR2738347, MR3021794}, Zeitouni and Zelditch studied the large deviation principle for zeros for compact Riemann surfaces; we also refer to \cite{MR4692882, Dinh:2024aa, Wu:2024aa} for recent breakthroughs on this topic, in particular, the hole probabilities of random zeros on compact Riemann surfaces (cf. Proposition \ref{prop:1.2.3}).

Shiffman and Zelditch \cite{MR1675133} first established the general framework for the random zeros of holomorphic sections in K\"{a}hler geometry, by using the Bergman kernel expansions. Then in their series of work \cite{MR1675133, SZ08, SZZ08, MR2742043, MR4293941}, the equidistribution, the large deviation, the number variance, and the central limit theorem for random zeros were proven for the positive line bundles on compact K\"{a}hler manifolds. The first named author with Drewitz and Marinescu in their work \cite{Drewitz_2023,DrLM:2023aa,Drewitz:2024aa} extended the aforementioned results to the uniformly positive line bundles on non-compact Hermitian manifolds. In particular, a probabilistic Berezin-Toeplitz quantization was introduced in \cite{DrLM:2023aa, Drewitz:2024aa} by considering square-integrable Gaussian holomorphic sections.

Note that Dinh and Sibony \cite{MR2208805} gave a different approach for the equidistribution of random zeros which also provides estimates for the speed of convergence, see \cite{DMS12, DMM16}. We also refer to the survey \cite{MR3895931} for more references on the topics of random zeros in complex geometry.

Now, we give in detail the geometric setting and the main results of this paper.
\subsection{Semipositive line bundles over punctured Riemann surfaces}\label{ss:1.1amiens}

Let $\bar{\Sigma}$ be a compact Riemann surface, and let
$D=\{a_{1},\cdots, a_{N}\} \subset \bar{\Sigma}$ be a finite set of points. We consider the punctured Riemann surface 
$\Sigma=\bar{\Sigma}\setminus D$ together with a Hermitian form
$\omega_{\Sigma}$ on $\Sigma$. We always fix an imaginary unit
$\ui=\sqrt{-1}$.

Let $T\Sigma$ denote the real tangent
bundle of $\Sigma$, and let $J\in 
\mathrm{End}(T\Sigma)$ denote the complex structure of $\Sigma$. Then we have the bidegree splitting
\begin{equation}
	T\Sigma\otimes_{\R}\C=T^{(1,0)}\Sigma\oplus T^{(0,1)}\Sigma.
	\label{eq:1.1.1thalys}
\end{equation}
Then $\omega_{\Sigma}$ is a real $(1,1)$-form such that 
$\omega_{\Sigma}(\cdot, J\cdot)$ is a Riemannian metric $g^{T\Sigma}$ on $T\Sigma$. 
Moreover, $\omega_{\Sigma}$ is K\"{a}hler. Let $\nabla^{T\Sigma}$ 
denote the Levi-Civita connection associated with $g^{T\Sigma}$, then 
it preserves the splitting \eqref{eq:1.1.1thalys}, we write it as
\begin{equation}
	\nabla^{T\Sigma}=\nabla^{T^{(1,0)}\Sigma}\oplus\nabla^{T^{(0,1)}\Sigma}.
\end{equation}
In particular, $\nabla^{T^{(1,0)}\Sigma}$ is exactly the 
Chern connection on the holomorphic line bundle $T^{(1,0)}\Sigma$ 
associated with the Hermitian metric 
$h^{T^{(1,0)}\Sigma}(\cdot,\cdot)=g^{T\Sigma}(\cdot,\bar{\cdot})$.

Let $L$ be a holomorphic line bundle on $\bar{\Sigma}$, 
and let $h$ be a singular Hermitian metric on $L$ such that:
\begin{enumerate}[label=(\greek*)]
	\item \label{itm:1}$h$ is smooth over $\Sigma$ and for all $j \in 
	\{1, \ldots, N\}$ there exists a trivialization of $L$ in the neighborhood $\bar{V_j}$ of $a_j$ in $\bar{\Sigma}$ with associated coordinate $z_j$ ($a_{j}$ corresponds to $z_{j}=0$) such that
	\begin{align*}
		\abs{1}_{h}^2 (z_j) = \abs{\log \abs{z_j}^2}.
	\end{align*}
	\item \label{itm:2} The Chern curvature $R^L = (\nabla^L)^2$ of $h$ satisfies
	\begin{itemize}
		\item[(i)] On $\Sigma$, we have $\ui R^{L}\geqslant 0$.
		\item[(ii)] For each $j \in \{1,\ldots, N\}$, we have $\ui R^L =\omega_{\Sigma}$ on $V_j := \bar{V}_j \setminus \{a_{j}\}$.
		\item[(iii)] $R^L$ vanishes at most to finite order at any point $x \in \Sigma$, that is,
		\begin{align*}
			\ord_x(R^L) &:= \min \lb \ell \in \N: J^\ell(\Lambda^2 T^\ast \Sigma) 
			\ni j^\ell_x R^L \neq 0 \rb \, <\infty,
		\end{align*}
		where $J^\ell(\Sigma; \Lambda^2T^\ast \Sigma)$ denotes the 
		$\ell$-th jet bundle over $\Sigma$ (see \hyperref[ssec:jet]{Appendix}). 
	\end{itemize}
\end{enumerate}

By assumptions \ref{itm:1} and \ref{itm:2} - (ii), in the local 
coordinate $z_{j}$ on $V_{j}$, we have $\omega_{\Sigma} = 
\omega_{\mathbb{D}^{*}}$ is the Poincar{\'e} metric on punctured unit 
disc given as follows
\begin{equation}
	\omega_{\mathbb{D}^{*}}=\frac{\ui\dz\wedge 
	\dzb}{|z|^{2}\log^{2}(|z|^{2})}.
	\label{eq:1.1.1}
\end{equation}
Then $(\Sigma, \omega_\Sigma)$ is complete, and the volume of 
$\Sigma$ with respect to the Riemannian volume form $\omega_{\Sigma}$ is finite. Let $\mathrm{dist}(\cdot,\cdot)$ denote the Riemannian distance on $\Sigma$.
\\

One typical example of a semipositive line bundle as described above is from branched coverings. If $f: \Sigma \to \Sigma^0$ is a branched covering of a Riemann surface $\Sigma^0$ with
branch points $\{y_1, \ldots , y_M \} \subset  \Sigma$, the Hermitian holomorphic line bundle on $\Sigma$, that is defined as the pullback of a positive one on $\Sigma^0$, becomes semipositive with curvature vanishing at the branch points (see \cite[Example 17]{Marinescu2023}).
\\

For $x\in\Sigma$, we set
\begin{equation}
	\rho_{x}=2+\ord_{x}(R^{L})\in \N_{\geqslant 2}.
	\label{eq:1.1.1bis}
\end{equation}
The function $x \mapsto \rho_x$ is upper 
semi-continuous on $\Sigma$, and the 
assumptions \ref{itm:2} - (ii) and (iii) infer that
\begin{equation}
	\rho_{\Sigma}:=\max_{x \in \Sigma} \rho_x < \infty 
	\label{eq:defofr}
\end{equation}
The semi-positivity in assumption \ref{itm:2} - (i) implies that $\rho_{x}$ is even for all $x \in \Sigma$, and so is $\rho_{\Sigma}$. 
Moreover, we have a decomposition $\Sigma = \bigcup_{j=2}^{\rho_{\Sigma}} 
\Sigma_j$, with $\Sigma_j := \{x \in \Sigma: \rho_x = j\}$; each $\Sigma_{\leqslant j} = \bigcup_{j^\prime=2}^j \Sigma_{j^\prime}$ is open. In particular, $\Sigma_2$ is an open dense subset of $\Sigma$. Note that $\ui R^{L}$ is strictly positive on
$\Sigma_2$, consequently, we have
\begin{equation}
	\label{eq:1.1.4}
	\mathrm{deg}(L) = \int_\Sigma \frac{\ui}{2\pi} R^L > 0,
\end{equation}
so that $L$ is ample, hence positive, over $\overline{\Sigma}$ (see also \cite{MR615130}). 
\\

From now on, we also fix a holomorphic line bundle $E$ over $\Sigma$ 
with a smooth Hermitian metric $h^{E}$, and we assume that $(E,h^{E})$ is identical to the trivial complex line bundle with the trivial Hermitian metric on each $V_{j}$ (in assumption \ref{itm:2}).

For $p\geqslant 1$, we denote by $h_{p}:=h^{\otimes p}\otimes h^{E}$ the metric induced by 
$h$ on $L^{p}\otimes E:=L^{\otimes p}\otimes E$ on $\Sigma$. Let 
$H^{0}(\Sigma, L^{p}\otimes E)$ 
be the space 
of holomorphic sections of $L^{p}\otimes E$ on $\Sigma$ and let
$\mathcal{L}^{2}(\Sigma, L^{p}\otimes E)$ be the space of 
$\mathcal{L}^{2}$-sections of $L^{p}\otimes E$ on $\Sigma$ with respect to $h_p$ and $\omega_\Sigma$. Set
\begin{equation}
	H^{0}_{(2)}(\Sigma, L^{p}\otimes E)
	=H^{0}(\Sigma, L^{p}\otimes E)\cap\mathcal{L}^{2}(\Sigma, L^{p}\otimes E),
	\label{eq:1.1.2}
\end{equation}
which is equipped with the associated $\mathcal{L}^{2}$-metric. 
Then by the integrability near the punctures, the sections in 
$H^{0}_{(2)}(\Sigma, L^{p}\otimes E)$ extend to holomorphic sections of 
$L^{p}$ over $\overline{\Sigma}$:
\begin{equation}
	H^{0}_{(2)}(\Sigma, L^{p}\otimes E)\subset H^{0}(\overline{\Sigma}, L^{p}\otimes E).
	\label{eq:1.1.8ab}
\end{equation}
Moreover, for $p\geqslant 2$, elements in $H^{0}_{(2)}(\Sigma, L^{p}\otimes E)$ are 
exactly the sections in $H^{0}(\overline{\Sigma}, L^{p}\otimes E)$ that vanish on the 
puncture divisor $D$ (cf. \cite[Remark 3.2]{AMM21} \cite[Section 4]{AMM22}).
Let $g$ denote the genus of $\overline{\Sigma}$. Then by the Riemann-Roch formula for $p\gg 1$, we have
    \begin{equation}
        d_p:=\dim H^0_{(2)}(\Sigma, L^p\otimes E) =p\deg(L) + \deg (E) +1-g-N
    \end{equation}

Let 
\begin{equation}
B_p: \mathcal{L}^{2}(\Sigma, L^{p}\otimes E)\longrightarrow  H^{0}_{(2)}(\Sigma, L^{p}\otimes E)
\label{eq:1.1.9amiens}
\end{equation} 
denote the orthogonal projection, which is known as Bergman projection. We will denote its Schwartz kernel, the Bergman kernel, 
by $B_{p}(x,y)$ for $x, y\in \Sigma$. If $S^{p}_{j}$, 
$j \in \{1,\ldots, d_{p}\}$ is an orthonormal basis of $H^{0}_{(2)}(\Sigma, 
L^{p}\otimes E)$ with respect to the $\mathcal{L}^{2}$-inner product, then 
\begin{equation}
	B_{p}(x,y)=\sum_{j=1}^{d_{p}}S^{p}_{j}(x)\otimes 
	S^{p,*}_{j}(y)\in (L^{p}\otimes E)_{x}\otimes (L^{p}\otimes 
	E)^{\ast}_{y},
	\quad\text{for $x,y\in \Sigma$,}
\end{equation}
where the duality is defined by $h_{p}$. In particular, $B_{p}(x):=B_{p}(x,x)$ is a nonnegative smooth function in $x\in\Sigma$, which is called the Bergman kernel function.

\subsection{Spectral gap and Bergman kernel expansion}

With the geometric setting described in the previous section, one of the main objects of investigation in this paper is the asymptotic expansion of the Bergman kernels $B_p(x,y)$ as $p\to +\infty$. There are two ingredients in our approach: the first one extends the result of Marinescu and Savale \cite{Marinescu2023, MS23} for a semipositive line bundle on a compact Riemann surface to our punctured Riemann surface, from which we prove a spectral gap for the Kodaira Laplacians; the second is the technique of analytic localization developed by Dai--Liu--Ma \cite{DLM06} and Ma--Marinescu 
\cite{MM07}, which is inspired by the work of Bismut--Lebeau 
\cite{BL91} in local index theory. In order to deal with the Bergman kernel near the punctures, we will follow the seminal work of Auvray, Ma, and Marinescu \cite{AMM16, AMM21}.

\begin{theorem}[Spectral gaps]\label{thm:1.2.1}
	Let $\Sigma$ be a punctured Riemann surface, and let $L$ 
	be a holomorphic line bundle as above such that $L$ carries a singular 
	Hermitian metric $h$ satisfying conditions \ref{itm:1} and 
	\ref{itm:2}. Let $E$ be a holomorphic line bundle on 
	$\Sigma$ equipped with a smooth Hermitian metric $h^{E}$ such 
	that $(E,h^{E})$ on each chart $V_{j}$ is exactly trivial 
	Hermitian line bundle.
	Consider the Dirac operator $D_p$ and Kodaira Laplacian $\square_p$ as in Subsection \ref{ss2.1P}. 
	Then there exist constants $C_1, C_2 \in \R_{>0}$ independent of $p$, 
	such that for all $s \in 
	\Omega^{0,1}_{\mathrm{c}}(\Sigma,L^{p}\otimes E)$,
	\begin{enumerate}[label=(\roman*)]
		\item \label{dirac}the Dirac operators are bounded from below, \begin{equation}
		\norm{D_p s}^2_{\cLL} \geqslant 2(C_1 p^{\sfrac{2}{\rho_{\Sigma}}} - 
		C_2) \norm{s}^2_{\cLL} \, ,
		\label{eq:dirac}
	\end{equation}
	\item \label{spec} for $p\in\N$, we have
	\begin{equation}
		\mathrm{Spec}(\square_p) \subset \{0\} \cup \left[C_1 
		p^{\sfrac{2}{\rho_{\Sigma}}} - C_2, +\infty\right[.
		\label{eq:spec}
	\end{equation}
\end{enumerate}
In particular, we have the first $\cLL$-Dolbeault cohomology group (see Subsection \ref{ss2.1P})
$$H^{1}_{(2)}(\Sigma,L^{p}\otimes E) = 0$$
for $p \gg 0$.
\end{theorem}

The proof of the spectral gap will be given in Subsection \ref{subsection:2.2}. As a consequence, we have the following pointwise expansions for the Bergman kernel functions, which extend the result of Marinescu and Savale \cite[Theorem 3]{Marinescu2023} to our non-compact setting.

\begin{theorem}[Asymptotic expansion of Bergman kernel functions]\label{thm:bkexpansion}
	We assume the same conditions on $\Sigma$, $L$ and $E$ as in Theorem \ref{thm:1.2.1}. Fix $\rho_{0}\in\{2, 4, \ldots,\rho_{\Sigma}\}$, and let $W: [0,1]\ni s\mapsto W(s)\in \Sigma$ be a smooth path such 
	that $W(s)\in \Sigma_{\rho_{0}}$ for all $s\in [0,1]$. Then for every $r \in \N$, there exists a smooth function $b_{r}(x)$ in $x\in W([0,1])$ such that for any $k\in\N$, we have the following asymptotic expansion of Bergman kernel functions uniformly on $W([0,1])$,
	\begin{equation}
		\label{eq:bkexpansion}
		B_p(x) = p^{\sfrac{2}{\rho_{0}}} \left[ \sum_{r=0}^k b_r(x) 
		p^{-\sfrac{2r}{\rho_0}} \right] + 
		\mathcal{O}(p^{-\sfrac{2k}{\rho_0}}) \, ,
	\end{equation}
	where the expansion holds in any $\mathscr{C}^{\ell}$-norms on $W([0,1])$ 
	with $\ell\in\N$. Moreover, for $x\in W([0,1])$,
	\begin{equation}
		b_{0}(x)=B^{j^{\rho_{0}-2}_{x}R^{L}}(0,0)>0,
  \label{eq:1.2.4june}
	\end{equation}
	where $j^{\rho_{0}-2}_{x}R^{L}\in 
	\ui S^{\rho_{0}-2}\R^{2}\otimes\Lambda^{2}(\R^{2})^{\ast}$ is 
	defined 
	as the $(\rho_{0}-2)$-degree homogeneous part of the Taylor expansion 
	of $R^{L}$ in the geodesic normal coordinate centered at $x$, and 
	$B^{j^{\rho_{0}-2}_{x}R^{L}}$ is the model Bergman projection that will be defined in Subsecton \ref{ss3.2P}.
	
	For $t\in\;]0,1[\;$, $\gamma\in \;]0,\frac{1}{2}[\;$, $\ell, m\in\N$, and $V_{j}$ described in assumption 
\ref{itm:1} with coordinate $z_{j}$ (it is clear that $\rho_{z_{j}}=2$), the following expansions hold uniformly in 
$\mathscr{C}^{\ell}$-norm for points $z_{j}\in \D^{\ast}(a_{j},\frac{1}{6})\setminus 
\D^{\ast}(a_{j}, te^{-p^{\gamma}})$,
\begin{equation}
	B_{p}(z_{j})=\frac{p-1}{2\pi}+\mathcal{O}(p^{-m}).
 \label{eq:24-1.2.5}
\end{equation}
	\end{theorem}

Define the nonnegative bounded smooth function $\boldsymbol{c}$ on $\Sigma$ as follows,
\begin{equation}
 \boldsymbol{c}(x)=\frac{\ui R^L_{x}}{\omega_{\Sigma}(x)}\geqslant 0. 
 \label{eq:1.2.6june}
\end{equation}
Then for the points $x\in \Sigma_2$ (that is $ \boldsymbol{c}(x) >0$), the function given in \eqref{eq:1.2.4june} is
\begin{equation}
    b_{0}(x)=\frac{\boldsymbol{c}(x)}{2\pi}.
\end{equation}
In particular, as in \eqref{eq:24-1.2.5}, $b_{0}(x)=\frac{1}{2\pi}$ (or, equivalently, $\boldsymbol{c}(x)=1$) near the punctures.

For $t\in\;]0,1[\;$, $\gamma\in \;]0,\frac{1}{2}[\;$, set 
\begin{equation}
    \Sigma_{p,t,\gamma}=\Sigma \setminus \bigcup_{j=1}^{N} 
\mathbb{D}^{\ast}(a_{j},te^{-p^\gamma}),
\label{eq:1.2.8-intro-24}
\end{equation}
where 
$\mathbb{D}^{\ast}(a_{j},te^{-p^\gamma})$ denote the punctured (open) disc 
of radius $te^{-p^\gamma}$ centered at $a_{j}$ 
in the coordinate $z_{j}\in V_{j}$ described in assumption 
\ref{itm:1}. Then we have the convergence of subsets
$$\lim_{p\to +\infty}\Sigma_{p,t,\gamma}=\Sigma.$$
As a consequence of Theorem \ref{thm:bkexpansion}, we have the following uniform upper bound on $B_p(x)$ when $x$ stays in $\Sigma_{p,t,\gamma}$.

\begin{corollary} \label{lem:bkupper}
	Set
 \begin{equation}
		C_0:=\sup_{x \in \Sigma} 
		\frac{\boldsymbol{c}(x)}{2\pi} \geqslant \frac{1}{2\pi}.
  \label{eq:1.2.9-24}
	\end{equation}
Then for any fixed $t\in\;]0,1[\;$, $\gamma\in \;]0,\frac{1}{2}[\;$, we have for $p\gg 1$,
	\begin{equation}
			\sup_{x\in \Sigma_{p,t,\gamma}} B_p(x)
		\leqslant C_0
		\left(1+o(1)\right)p,
	\end{equation}
	where the small o-term $o(1)$ is uniform in $x\in \Sigma_{p,t,\gamma}$ as $p\to +\infty$.
\end{corollary}

In the pointwise expansion of $B_p(x)$, the leading term grows as $p^{2/\rho_x}$ ($\rho_x\geqslant 2$). Corollary \ref{lem:bkupper} describes this upper bound for the point $x\in \Sigma_{p,t,\gamma}$, which still keeps at least an exponentially small distance from the punctures.
However, our assumptions about punctures implies that a global supremum of $B_p(x)$ on $\Sigma$ will behave like $p^{3/2}$, as $p \to +\infty$, following the work of Auvray--Ma--Marinescu \cite{AMM21} for the Poincar\'{e} punctured disc.
\begin{proposition}\label{cor:1.2.3}
	We assume the same conditions on $\Sigma$, $L$ and $E$ as in Theorem \ref{thm:1.2.1} with the number of punctures $N\geqslant 1$. We have
	\begin{equation}
	\sup_{x\in\Sigma} 
	B_{p}(x)=\left(\frac{p}{2\pi}\right)^{\sfrac{3}{\,2}}+\mathcal{O}(p).
 \label{eq:1.2.11june}
\end{equation}
\end{proposition}

The proofs of Theorem \ref{thm:bkexpansion}, Corollary \ref{lem:bkupper}, and Proposition \ref{cor:1.2.3} will be presented in Subsection \ref{ss:4.4june}. In Theorem \ref{thm:4.4.1june} we also obtain the pointwise expansions of the derivatives of $B_p(x)$. Moreover, considering the Kodaira maps defined with $H^0_{(2)}(X,L^p\otimes E)$, a version of Tian's approximation theorem \cite{Tia90} will be given in Subsection \ref{ss5.2Tian}.

In \cite[Section 3.1]{MS23}, on a compact Riemann surface equipped with a semipositive line bundle, the uniform estimates of the upper and lower bounds for the Bergman kernel functions were discussed (in this case, Proposition \ref{cor:1.2.3} does not apply), and the analogous results can be smoothly extended to our setting. Here, we will not discuss such uniform estimates, but we will focus on the near-diagonal expansions of $B_p$, Theorems \ref{thm:24-4.3.1} and \ref{thm:4.3.3feb24}, and their consequences for the study of random zeros in $\Sigma$. More precisely, we will be concerned with the semi-classical limit of the zeros of the Gaussian holomorphic sections for the higher tensor powers of $L$ but associated to a semipositive Hermitian metric on $L$. The following three subsections are dedicated to explain our results for random zeros, which lie in the framework of the smooth statistics of random point processes in $\Sigma$.

Now, as an extension of \cite[Proposition 5.3]{AMM21}, we give off-diagonal estimates for the Bergman kernels; see Subsection \ref{ss3.3P} for a proof. Fix $0<r<e^{-1}$, and fix a smooth function $\eta:\Sigma\to [ 1,\,\infty[$ such that $\eta(z)=|\log{|z|^{2}}|$ for $z\in \mathbb{D}^{\ast}_{r}$ near each puncture. 
\begin{proposition}[Off-diagonal estimates on Bergman kernels]\label{Prop: off-diag-Berg}
Fix a sufficiently small $\varepsilon>0$. Given $m,\ell\in \bN$, $\gamma>1/2\;$, there exists $C_{\ell,m,\gamma}>0$ such that for $z,z'\in\Sigma$, $\dist(z,z')\geqslant \varepsilon$, we have
	\begin{equation}
		\left|\eta(z)^{-\gamma}\eta(z')^{-\gamma}B_{p}(z,z')\right|_{\mathscr{C}^{m}(h_{p})}\leqslant C_{\ell,m,\gamma}p^{-\ell},
		\label{eq:3.3.8P}
	\end{equation}    
 where $|\cdot|_{\mathscr{C}^{m}(h_{p})}$ is the $\mathscr{C}^{m}$-norm induced by $g^{T\Sigma}$, $h_p$ and the corresponding connections.
\end{proposition}

\subsection{Equidistribution of zeros of Gaussian holomorphic sections}

Recall that, with the assumptions described in Subsection \ref{ss:1.1amiens}, $H^0_{(2)}(\Sigma, L^p\otimes E)$ equipped with the $\cLL$-inner product is a Hermitian vector space of dimension $d_p<\infty$.

For a non-trivial holomorphic section $s_p\in H^0_{(2)}(\Sigma, L^p\otimes E)$, the zeros of $s_p$ consist of isolated points in $\Sigma$. We consider the divisor
\begin{equation}
    \Div(s_p):=\sum_{x\in\Sigma,\, s_p(x)=0} m_x\cdot x,
\end{equation}
where $m_x$ denotes the multiplicity of $x$ as a zero of $s_p$ (or vanishing order). Then we define the following measure on $\Sigma$, 
\begin{equation}
[\Div(s_p)]:=\sum_{s_p(x)=0} m_x \delta_x,
\end{equation}
where $\delta_x$ denotes the Dirac mass at $x$.

Then the Poincar\'{e}-Lelong formula states an identity for the distributions on $\Sigma$,
\begin{equation}
    [\Div(s_p)]=\frac{\ui}{2\pi} \partial\overline{\partial} \log |s_p(x)|^2_{h_p}+p c_1(L,h)+c_1(E,h^E).
    \label{eq:1.3.3-24june}
\end{equation}
At the same time, we introduce the following norm for the distributions on $\Sigma$: let $T$ be a distribution on $\Sigma$, for any open susbet $U\subset \Sigma$, define
\begin{equation}
    \|T\|_{U,-2}:=\sup_{\varphi}|\langle T,\varphi\rangle|,
    \label{eq:1.3.4june-24}
\end{equation}
where the supremum is taken over all the smooth test functions $\varphi$ with support in $U$ and such that their $\mathscr{C}^2$-norm satisfies $\|\varphi\|_{\mathscr{C}^2}\leqslant 1$.

In the sequel, our main object is to study the asymptotic behaviours of $[\Div(s_p)]$ for random sequences of $s_p$'s as $p\to +\infty$, which can be viewed as a random point process on $\Sigma$. Let us start with the Gaussian holomorphic sections.

\begin{definition}[Standard Gaussian holomorphic sections]\label{defn:Gaussian}
    On $H^0_{(2)}(\Sigma, L^p\otimes E)$, we define the standard Gaussian probability measure $\P_p$ associated to the $\cLL$-inner product. Let $\bb{S}_p$ be the random variable valued in $H^0_{(2)}(\Sigma, L^p\otimes E)$ with the law $\P_p$, which is called the standard Gaussian holomorphic sections of $(L^p\otimes E, h_p)$ over $\Sigma$. We also set the product probability space 
    $$ (H_\infty, \P_\infty):=\prod_p \left(H^0_{(2)}(\Sigma, L^p\otimes E),\P_p\right)$$
    whose elements are the sequences $\{s_p\}_p$ of holomorphic sections.

    We have an equivalent definition. Let $\{S^p_j\}_{j=1}^{d_p}$ be an orthonormal basis of $H^0_{(2)}(\Sigma, L^p\otimes E)$ and let $\{\eta^p_j\}_{j=1}^{d_p}$ be a vector of independent and identically distributed (i.i.d.) standard complex Gaussian variables (that is $\mathcal{N}_\C(0,1)$), then we can also write
    \begin{equation}
        \bb{S}_p=\sum_{j=1}^{d_p} \eta^p_j S^p_j.
    \end{equation}
    Note that these random variables are taken independently for different $p$'s. We will always use equally the above two models to state our results.
\end{definition}

Now we can give the equidistribution results for the random zeros $[\Div(\bb{S}_p)]$, which states that the measures defined from random zeros will asymptotically converge to the semipositive smooth measure $c_1(L,h)$ on $\Sigma$. The proof will be given in Subsection \ref{ss:equidis}, and we refer to Definition \ref{defn:speed} for the notion of convergence speed.
\begin{theorem}[Equidistribution of {$[\Div(\bb{S}_p)]$}]\label{thm:1.3.2-24}
   	We assume the same conditions on $\Sigma$, $L$ and $E$ as in Theorem \ref{thm:1.2.1}. 
    \begin{enumerate}[label=(\roman*)]
        \item\label{1.3.2-1} The expectation $\E[[\Div(\bb{S}_p)]]$, as a measure on $\Sigma$\,, exists, and as $p\to +\infty$, we have the weak convergence of measures 
       \begin{equation}
           \frac{1}{p}\E[[\Div(\bb{S}_p)]] \longrightarrow c_1(L,h),
           \label{eq:1.3.6-24}
       \end{equation}
  and for any relatively compact open subset $U$ in $\Sigma$, the above convergence has the convergence speed $\mathcal{O}(\log p /p)$ on $U$, that is, there exists a constant $C_U>0$ such that
       $$\left\|\frac{1}{p}\E \left[[\Div(\bb{S}_p)] \right] - c_1(L,h)\right\|_{U,-2}\leqslant C_U\frac{\log p}{p}.$$
       \item\label{1.3.2-2}  For $\P_\infty$-almost every sequence $\{s_p\}_p$\,, we have the weak convergence of measures on $\Sigma$,
              \begin{equation}
           \frac{1}{p}[\Div(s_p)] \longrightarrow c_1(L,h).
       \end{equation}
       Moreover, given any relatively compact open subset $U\subset \Sigma$\, , for $\P_\infty$-almost every sequence $\{s_p\}_p$\,, the above convergence on $U$ has convergence speed $\mathcal{O}(\log p /p)$.
    \end{enumerate}
\end{theorem}

 In order to obtain the convergence speed in Theorem \ref{thm:1.3.2-24} - \ref{1.3.2-2}, we need to use a result - Theorem \ref{thm:5.3.1-24} - of Dinh, Marinescu, and Schmidt \cite{DMS12} (see also \cite[Theorems 1.1 and 3.2]{DMM16}), motivated by the ideas of Dinh and Sibony \cite{MR2208805}.

\subsection{Normalized Bergman kernel and large deviations of random zeros}

Now we consider the normalized Bergman kernel, which will play the role of correlation functions of $\bb{S}_p$ (in Definition \ref{defn:Gaussian}), viewed as the holomorphic Gaussian fields on $\Sigma$. The normalized Bergman kernel is defined as
\begin{equation}
	N_{p}(x,y)=\frac{|B_{p}(x,y)|_{h_{p,x}\otimes 
	h_{p,y}^{\ast}}}{\sqrt{B_{p}(x,x)}\sqrt{B_{p}(y,y)}}, \quad  x,y\in \Sigma.
	\label{eq:1.3.1}
\end{equation}
Due to the positive of $L$ on $\overline{\Sigma}$, for any compact subset $K$ of $\Sigma$ and all sufficiently large $p\gg 1$, the function $N_{p}(x,y)$ is smooth on $K\times K$ with values in $[0,1]$.

Let $\inj^U$ denote the injectivity radius for a subset $U\subset \Sigma$ (see \eqref{eq:4.2.1june}). Then we have the following near-diagonal expansions of $N_{p}(x,y)$ only for the points $x,y\in \Sigma_2$. At a vanishing point $x$ of $R^L$, due to the lack of the explicit formula for the model Bergman kernel $B^{R^L_0}_x$, such near-diagonal expansions remain unclear.

\begin{theorem}\label{thm:1.3.1}
	Let $U$ be a relatively compact open subset of $\Sigma_2\subset \Sigma$ (hence $\ui R^L$ is strictly positive on $\overline{U}$), and set
 $$\varepsilon_0:=\inf_{x\in U} \bb{c}(x)>0,$$
 where $\boldsymbol{c}(x)=\ui R^L_x/\omega_{\Sigma}(x)$ is a strictly positive function on $\Sigma_2$.
 Then there exists $\delta_U \in \;]0, \inj^{U}/4[\;$ such that we have 
	the following uniform estimate on the normalized Bergman kernel: fix $k\geqslant 1$ and $b\geqslant \sqrt{12k/\varepsilon_{0}}\,$, then we have
 \begin{enumerate}[label=(\roman*)]
     \item\label{1.4.1-1}  There exists $C>0$ such that for all $p$ with $b\sqrt{\log p/{p}}\leqslant \delta_U\,$,
and all $x,y\in U$ 
with $\operatorname{dist}(x,y) \geqslant b\sqrt{\log p/{p}}$ we have
$N_{p}(x,y)\leqslant C p^{-k}$.

\item\label{1.4.1-2}   There exist functions 
$$R_p:\left\{(x,y)\in U\times U:\operatorname{dist}(x,y) \leqslant 
b\sqrt{\tfrac{\log p}{p}}
\,\right\}\to\R$$
such that  $\sup |R_p|\to 0$ as $p\to\infty\;$, and such that for all sufficiently large $p$,
	\begin{equation}
		N_{p}(x,y)= 
		 (1+ R_{p}(x,y))\exp \left\{-\frac{\boldsymbol{c}(x)p}{4}\operatorname{dist}(x,y)^{2} \right\}.
		\label{eq:1.5.2July}
	\end{equation}

\item\label{1.4.1-3}   Moreover, for any $\varepsilon \in \;]0,1/2]\,$, there exists $C=C(U,b,k,\varepsilon)>0$ such that for all sufficiently large $p\;$,
 \begin{equation}
     \sup |R_p| \leqslant Cp^{-1/2 +\varepsilon}.
 \end{equation}
 \end{enumerate}
\end{theorem}

In the case of compact K\"{a}hler 
manifolds with positive line bundles, such results were established in \cite[Propositions 2.6 and 2.7]{SZ08} and in \cite[Proposition 2.1]{SZZ08}. In the non-compact complete Hermitian manifolds with uniformly positive line bundles, by applying the Bergman kernel expansion obtained by Ma and 
Marinescu \cite[Theorems 4.2.1 and 6.1.1]{MM07}, such results are proven in \cite[Theorems 1.8 and 5.1]{Drewitz_2023} (see also \cite[Theorem 3.13]{DrLM:2023aa}). Note that, comparing with  \cite[Theorems 1.8]{Drewitz_2023}, we have improved some estimates in our Theorem \ref{thm:1.3.1}. For normalized Berezin-Toeplitz kernels, the analogous result was given in \cite[Theorem 1.20 and Corollary 1.21]{Drewitz:2024aa}.

Recall that the Gaussian holomorphic section $\bb{S}_p$ is constructed in Definition \ref{defn:Gaussian}. For any open subset $U\subset \Sigma$, set
\begin{equation}
    \mathcal{N}^U_p(\bb{S}_p):=\int_U [\Div(\bb{S}_p)]=\sum_{x\in U, \bb{S}_p(x)=0} m_x.
\end{equation}
Then $\mathcal{N}^U_p(\bb{S}_p)$ is a random variable valued in $\bN$.

Note that $c_1(L,h)$ defines a nonnegative smooth measure on $\Sigma$, for any open subset $U$, we set
\begin{equation}
    \mathrm{Area}^{L}(U):=\int_U c_1(L,h).
\end{equation}

As a consequence of Theorem \ref{thm:1.3.1}, we obtain the following results for random zeros, which generalize \cite[Corollary 1.2 and Thoerem 1.4]{SZZ08} and \cite[Theorem 1.5, Corollary 1.6]{Drewitz_2023}. Their proof will be given in Subsection \ref{ss:LDE}.
\begin{theorem}[Large deviation estimates or concentration inequalities]\label{thm:5.3.1}
   	We assume the same conditions on $\Sigma$, $L$ and $E$ as in Theorem \ref{thm:1.2.1}.
    \begin{enumerate}[label=(\roman*)]
        \item\label{1.4.2-1} If $U$ is a relatively compact open subset in $\Sigma$, then for any $\delta>0$, there exists a constant $C_{\delta,U}>0$ such that for $p\gg 0$
the following holds:
\begin{equation}\label{eq:1.4.4-24}
    \P_{p} \left( \left\| \frac{1}{p}[\Div(\bb{S}_p)] - c_1(L,h) \right\|_{U,-2} > \delta \right) \leqslant e^{-C_{\delta,U}p^{2}}.
\end{equation}
        \item\label{1.4.2-2} If $U$ is an open set of $\Sigma$ with $\partial U$ having zero
measure with respect to some given smooth volume measure on 
$\overline{\Sigma}$\, ($U$ might not be relatively compact in $\Sigma$), then for any $\delta>0$, there exists a 
constant $C'_{\delta,U}>0$ such that for $p\gg 0$
the following holds:
\begin{equation}\label{eq:1.4.5DLM}
    \P_{p}\left( \left| \frac{1}{p}\mathcal{N}^{U}_{p}(\bb{S}_{p}) - \mathrm{Area}^{L}(U) \right| > \delta \right) \leqslant e^{-C'_{\delta,U}p^{2}}.
\end{equation}
As a consequence, for $\P_\infty$-almost every sequence $\{s_p\}_p\in H_\infty\;$, we have
\begin{equation}
    \frac{1}{p}\mathcal{N}^{U}_{p}(s_{p})\longrightarrow
\mathrm{Area}^{L}(U).
\label{eq:1.4.8-24}
\end{equation}
    \end{enumerate}
\end{theorem}

\begin{proposition}[Hole probabilities]\label{prop:1.2.3}	
If $U$ is a nonempty open set of $\Sigma$ with $\partial U$ having zero measure in $\Sigma$, then there exists a constant $C_{U}>0$ such that for $p\gg 0,$
	\begin{equation}
		\P_{p} \left(\mathcal{N}^{U}_{p}(\bb{S}_{p})=0 \right)
		\leqslant e^{-C_{U}p^{2}}.
  \label{eq:1.4.9-24june}
	\end{equation}

	 If $U$ is a relatively compact open subset of $\Sigma$ such that 
	$\partial U$ has zero measure in $\Sigma$, and if there exists a 
	section $\tau\in H^{0}_{(2)}(\Sigma, L)$ such that it does not vanish in $\overline{U}\subset \Sigma$, then there 
	exists $C'_{U,\tau}>0$ such that for $p\gg 0$,
	\begin{equation}
		\P_{p} \left(\mathcal{N}^{U}_{p}(\bb{S}_{p})=0 \right) \geqslant e^{-C'_{U,\tau}p^{2}}.
		\label{eq:1.2.4}
	\end{equation}
\end{proposition}

\subsection{Number variance and central limit theorem}
Under the geometric assumptions in Subsection \ref{ss:1.1amiens}, set
\begin{equation}
    \Sigma_\ast:=\bigcup_{j\geqslant 4} \Sigma_j=\{z\in\Sigma\;:\; R^L_z=0\}
\end{equation}
for the set of points in $\Sigma$ where the curvature vanishes.
Then it is known that the compact set $\Sigma_\ast$ has a measure zero with respect to $\omega_\Sigma$ (see also Lemma \ref{lm:5.5.1-24june}).

\begin{definition}\label{defn:1.5.1-24}
     Let $\varphi$ be a real $\cC^3$-function on $\Sigma$, we define a $\cC^1$-function $\fL(\varphi)$ on $\Sigma_2$ (we have to exclude the vanishing points of $c_1(L,h)$) by the following identity
     \begin{equation}
         \ui\partial\overline{\partial}\varphi = \fL(\varphi) c_1(L,h).
     \end{equation}
     In fact, up to a constant factor, $\fL(\varphi)$ is exactly the action of the Laplacian operator on $\varphi$ where the Laplacian operator is associated with the Hermitian metric $c_1(L,h)$ on $\Sigma_2$.

     To shorten our statements, we introduce the following class of test functions on $\Sigma$:
     \begin{equation}
     \mathcal{T}^3(L,h):= \left\{\varphi\in \cC^3_c(\Sigma, \R) : \partial\overline{\partial}\varphi \equiv 0 \text{ in a tubular} \text{ neighbourhood of $\Sigma_\ast$} \right\}. 
     \end{equation}
     Then for $\varphi\in \mathcal{T}^3(L,h)$, the real function $\fL(\varphi)$ is well-defined globally on $\Sigma$ that is identically zero near $\Sigma_\ast$.
\end{definition}

Recall that the definition of convergence in distribution is given as the pointwise convergence of the distribution functions towards the distribution function of the limiting random variable in all points of continuity.
The following result shows the asymptotic normality of the random zeros in $\Sigma$ under semi-classical limit, whose proof will be given in Subsection \ref{ss:CLT-24}.

\begin{theorem}[Central limit theorem]\label{thm:3.5.1ss}
We assume the same conditions on $\Sigma$, $L$ and $E$ as in Theorem \ref{thm:1.2.1}. Let $\varphi\in  \mathcal{T}^3(L,h) $ be such that $\partial\overline{\partial}\varphi\not\equiv 0$, set
	\begin{equation}
Y_{p}(\varphi):= \left\langle[\mathrm{Div}(\bb{S}_{p})],\varphi \right\rangle,
\label{eq:1.5.1-24}
	\end{equation}
then as $p\to \infty$, the distribution of the random 
variables
	\begin{equation}
\frac{Y_{p}(\varphi)-\E[Y_{p}(\varphi)]}{\sqrt{\mathrm{Var}[Y_{p}(\varphi)]}}
	\end{equation}
converges weakly to $\mathcal{N}_{\R}(0,1)$, standard real normal distribution.
\end{theorem}

Such kind of results as above were obtained by Sodin--Tsirelson \cite[Main Theorem]{STr} for Gaussian holomorphic functions and by Shiffman--Zelditch \cite[Theorem 1.2]{MR2742043} for positive line bundles on compact K\"{a}hler manifolds. Moreover, as pointed out in \cite[Remark 3.17]{DrLM:2023aa}, this result also holds for the standard Gaussian holomorphic sections $\{\boldsymbol{S}_p\}_p$ on noncompact Hermitian manifolds. Then in \cite[Theorem 1.17]{Drewitz:2024aa}, the first named author with Drewitz and Marinescu obtained a central limit theorem for the zeros of square-integrable Gaussian holomorphic sections via Berezin-Toeplitz quantization on complete Hermitian manifolds. All proofs of these results are based on the seminal result of Sodin and Tsirelson in \cite[Theorem 2.2]{STr} for the non-linear functionals of the Gaussian process (see Theorem \ref{thm:STr2004}).

 Note that in Theorem \ref{thm:3.5.1ss}, we need to take the test function $\varphi\in  \mathcal{T}^3(L,h) $. Since $\varphi$ does not necessarily vanish near $\Sigma_\ast$, such a kind of test function still allows variables $Y_p(\varphi)$ to contain the contributions of points in $\Sigma_\ast$. 
 
 Shiffman and Zelditch \cite{SZ08, MR2742043} established the framework to compute the asymptotics of $\mathrm{Var}[Y_{p}(\varphi)]$ on a compact K\"{a}hler manifold, in particular, they obtained a pluri-bipotential for it. Their method can be easily adapted to our setting, so that in Subsection \ref{ss:nv}, we will prove the following theorem.

\begin{theorem}[Number variance]\label{thm:1.5.2-24ss}
We assume the same conditions on $\Sigma$, $L$ and $E$ as in Theorem \ref{thm:1.2.1}. Fix any $\varepsilon \in \;]0,1/2]\,$. Let $\varphi\in  \mathcal{T}^3(L,h) $ be such that $\partial\overline{\partial}\varphi\not\equiv 0$, and let $Y_p(\varphi)$ be given as in \eqref{eq:1.5.1-24}, then we have the formula for $p\gg 0$,
\begin{equation}\label{eq:6.12feb24}
\mathrm{Var}[Y_{p}(\varphi)]=\frac{\zeta(3)}{4\pi^2 p }\int_\Sigma \abs{\fL(\varphi)(z)}^2 c_1(L,h)(z)+\mathcal{O}(p^{-3/2+\varepsilon}),
\end{equation}
where
$$\zeta(3)=\sum_{k=1}^\infty \frac{1}{k^{3}}\cong 1.202056903159594 \ldots$$ 
is the Ap\'{e}ry's constant.
\end{theorem}

With the same assumptions in Theorem \ref{thm:3.5.1ss}, by \eqref{eq:1.3.6-24}, we have 
$$p^{-1}\E[Y_{p}(\varphi)]\longrightarrow \langle c_1(L,h_L),\varphi\rangle=\int_\Sigma \varphi c_1(L,h)$$ as $p \to +\infty$. Therefore, as a consequence of Theorem \ref{thm:3.5.1ss} and \eqref{eq:6.12feb24} (also with Khintchine’s theorem \cite[Theorem 1.2.3]{MR691492}), we get the following result. 
\begin{corollary}\label{C:univ}
Under the same geometric assumptions of Theorem \ref{thm:3.5.1ss}, and take $\varphi\in  \mathcal{T}^3(L,h) $ with $\partial\overline{\partial}\varphi\not\equiv 0$,
the distributions of the real random variables 
\begin{equation}\label{eq:Div_loc}
    \sqrt{p} \, \left\langle[\mathrm{Div}(\bb{S}_{p})]-p c_1(L,h_L),\varphi \right\rangle , \; p\in\N,\,
\end{equation}
converge weakly to $\mathcal{N}_{\R}(0,\sigma(U,h,\varphi))$
as $p \to +\infty$, where
\begin{equation}
\sigma(U,h,\varphi):= \frac{\zeta(3)}{4\pi^2}
\int_\Sigma |\fL(\varphi)(z)|^2 c_1(L,h)(z)>0.
\end{equation}
\end{corollary}

\subsection*{Acknowledgments}
 The second author would like to express his gratitude to his Ph.D. advisor Prof. George Marinescu. The authors thank Dr. Nikhil Savale for many useful discussions.

\section{Semipositive line bundles and Spectral gap of Kodaira Laplacian}
In this section, we introduce the Dirac operators and Kodaira Laplacians on $\Sigma$. Following the work of Ma--Marinescu \cite{MM07}, of Auvray--Ma--Marinescu \cite{AMM21}, and of Marinescu--Savale \cite{Marinescu2023}, we prove the spectral gaps stated in Theorem \ref{thm:1.2.1}. Finally, we combine this spectral gap with a result of Hsiao and Marinescu \cite{MR3194375} to obtain the leading term of the Bergman kernel functions $B_p(x)$ on $\Sigma$. 
\subsection{{$\cLL$}-Dolbeault cohomology and Kodaira Laplacian}\label{ss2.1P}

Let
$\Omega^{0,\bullet}_{\mathrm{c}}(\Sigma, L^{p}\otimes E)$ denote the 
set of the smooth 
sections of $\Lambda^{\bullet}(T^{\ast(0,1)}\Sigma)\otimes 
L^{p}\otimes E$ on $\Sigma$ with compact support, and for $s\in 
\Omega^{0,\bullet}_{\mathrm{c}}(\Sigma, L^p\otimes E)$, the 
$\cLL$-norm of $s$ is given by
\begin{equation}
	\left\| s \right\|^{2}_{\cLL}:=\int_{\Sigma} \left| s \right|_{h_{p}}^{2} \, \omega_{\Sigma}.
	\label{eq:1.2.1}
\end{equation}
Let $\Omega^{0,\bullet}_{(2)}(\Sigma,L^p\otimes E)$ be the 
Hilbert space defined as the completion of 
$(\Omega^{0,\bullet}_{\mathrm{c}}(\Sigma,L^p\otimes E), \|\cdot\|_{\cLL})$, 
in particular, $\mathcal{L}^{2}(\Sigma, 
L^{p}\otimes E)=\Omega^{0,0}_{(2)}(\Sigma, L^p\otimes E)$.
As in \eqref{eq:1.1.2}, let $H^{0}_{(2)}(X,L^{p}\otimes E)$ denote the space of $\cLL$-holomorphic sections of $L^{p}\otimes E$ on $\Sigma$, 
which, by \eqref{eq:1.1.8ab}, is a finite-dimensional vector space equipped with the $\cLL$-inner 
product.

We consider the 
$\cLL$-Dolbeault complex, 
\begin{equation}
	0\rightarrow\Omega^{0,0}_{(2)}(\Sigma,L^p\otimes E) 
	\xrightarrow[\hspace{10mm}]{\delbar_p} 
	\Omega^{0,1}_{(2)}(\Sigma, L^p\otimes E)\rightarrow 0,
\end{equation}
where $\delbar_{p}$ is taken to be the maximal extension, that is, with 
the domain
\begin{equation}
	\Dom(\delbar_{p}):=\{s\in 
	\Omega^{0,0}_{(2)}(\Sigma,L^p\otimes E)\;:\; \delbar_{p}s\in 
	\Omega^{0,1}_{(2)}(\Sigma,L^p\otimes E) \}.
	\label{eq:2.1.2ab}
\end{equation}
Let $\delbar_{p}^{\ast}$ denote the maximal extension of the formal adjoint of $\delbar_{p}$ with 
respect to the $\cLL$-metrics, then since $(\Sigma, 
\omega_{\Sigma})$ is complete,  
$\delbar_{p}^{\ast}$ coincides with the Hilbert adjoint of 
$\delbar_{p}$. Let $H^{q}_{(2)}(\Sigma,
L^p)$, $q=0,1$, denote the $\cLL$-Dolbeault cohomology groups.

The Dirac operator $D_p$ and the Kodaira Laplacian operator $\square_p$ are 
given by
\begin{equation}
	\label{eq:1.2.5}
	\begin{split}
		D_p &:= \sqrt{2}(\delbar_p + \delbar_p^\ast), \\
		\square_p &:= \frac{1}{2} (D_p)^2  = \delbar_p \delbar_p^\ast + 
		\delbar_p^\ast \delbar_p \, .
	\end{split}
\end{equation}
Note that $\square_p: 
\Omega^{0,\bullet}_{\mathrm{c}}(\Sigma,L^p\otimes E) 
\longrightarrow \Omega^{0,\bullet}_{\mathrm{c}}(\Sigma, L^p\otimes E)$ is essentially self-adjoint, so it has a 
unique self-adjoint extension which we still denote by $\square_p$, the 
domain of this extension is $\Dom(\square_p) = \lb s \in 
\Omega^{0,\bullet}_{(2)}(\Sigma, 
L^p\otimes E)\;:\; \square_p(s) \in 
\Omega^{0,\bullet}_{(2)}(\Sigma, 
L^p\otimes E)\rb$.

Note that $D_p$ interchanges and $\square_p$ preserves the 
$\Z$-grading of $\Omega^{0,\bullet}_{\mathrm{c}}(\Sigma,L^p\otimes E)$. Then
\begin{equation}
	\label{eq:2.1.7}
	\begin{split}
		&\square^0_p:= {\square_p}_\restr{\Omega^{0,0}(\Sigma,
		L^p\otimes E)}=\delbar_p^{\ast} \delbar_p\, , \\
		&\square^1_p:= {\square_p}_\restr{\Omega^{0,1}(\Sigma,
		L^p\otimes E)}=\delbar_p \delbar_p^{\ast} \, .
	\end{split}
\end{equation}
Moreover, the completeness of $(\Sigma, g^{T\Sigma})$ infers that, for $q=0,1$,
\begin{equation}
	\label{eq:2.1.8}
	\ker \square^{q}_{p}\cong H^{q}_{(2)}(\Sigma, L^p\otimes E).
\end{equation}

For $x\in \Sigma$, $v\in T_{x}\Sigma$, by splitting
\eqref{eq:1.1.1thalys}, we write $ v 
= v^{(1,0)} + v^{(0,1)} \in T_{x}^{(1,0)}\Sigma \oplus T_{x}^{(0,1)}\Sigma$; 
we denote by $\bar{v}^{(1,0)\ast} \in T_{x}^{(0,1)\ast}\Sigma$ the metric 
dual of $v^{(1,0)}$. The Clifford multiplication endomorphism $c : 
T_{x}\Sigma \to 
\mathrm{End}(\Lambda^{\bullet}(T_{x}^{\ast(0,1)}\Sigma))$ is then defined as 
\begin{equation}
	v \mapsto c(v) := \sqrt{2}(\bar{v}^{(1,0)\ast} \wedge 
	-\iota_{v^{(0,1)}}),
	\label{eq:2.1.9thalys}
\end{equation}
where $\iota$ is the contraction operator.

If $\{e_{1},e_{2}\}$ is a local orthonormal frame of 
$(T\Sigma,g^{T\Sigma})$, then the Dirac operators in \eqref{eq:1.2.5} can then be written as follows:
\begin{equation}
	\label{eq:dirac-clifford}
	D_p = \sum_{j=1}^2 c(e_{j})\nabla^{\Lambda^{0,\bullet} \otimes 
	L^p\otimes E}_{e_{j}},
\end{equation}
where $\nabla^{\Lambda^{0,\bullet} \otimes 
L^p\otimes E}$ denote the Hermitian metric induced by $\nabla^{T\Sigma}$ 
and the Chern connections $\nabla^{L}$, $\nabla^{E}$.

Set $\omega=\frac{1}{\sqrt{2}}(e_{1}-\ui e_{2})$, it forms an orthonormal 
frame of $T^{(1,0)}\Sigma$. Let $\bar{\omega}^{\ast}$ denote the 
metric dual of $\omega$. By \cite[Theorem 1.4.7]{MM07}, let $\Delta^{\Lambda^{0,\bullet} \otimes 
L^p\otimes E}$ denote the Bochner Laplacian associated with $\nabla^{\Lambda^{0,\bullet} \otimes 
L^p\otimes E}$, we have the following formula for 
$\square_{p}$,

\begin{equation}
	\begin{split}
		\square_{p}=&\frac{1}{2}\Delta^{\Lambda^{0,\bullet} \otimes 
		L^p \otimes E} + \frac{r^{\Sigma}}{4} \, \bar{\omega}^{\ast} \wedge \iota_{\bar{\omega}}\\
		&+ p \left( R^{L}(\omega,\bar{\omega}) \, \bar{\omega}^{\ast} \wedge \iota_{\bar{\omega}}-\frac{1}{2}R^{L}(\omega,\bar{\omega}) \right) +  \left( R^{E}(\omega,\bar{\omega}) \, \bar{\omega}^{\ast} \wedge \iota_{\bar{\omega}}-\frac{1}{2}R^{E}(\omega,\bar{\omega}) \right),
	\end{split}
	\label{eq:2.1.12cologne}
\end{equation}
where $r^{\Sigma}=2R^{T^{(1,0)}\Sigma}(\omega,\bar{\omega})$ is the 
scalar curvature of $(\Sigma,g^{T\Sigma})$. Note that $r^{\Sigma}$ is 
a bounded 
function on $\Sigma$ which is constant near punctures. In particular, near the 
punctures, 
\begin{equation}
	R^{E}(\omega,\bar{\omega}) \, \bar{\omega}^{\ast} \wedge 
	\iota_{\bar{\omega}}-\frac{1}{2}R^{E}(\omega,\bar{\omega})=0,
\end{equation}
and we have more explicit formula for $\square_{p}$ as given 
in \cite[(4.15)]{AMM21}.

\subsection{Spectral gap: proof of Theorem {\ref{thm:1.2.1}}}\label{subsection:2.2}
Now we consider the action of $\square_{p}$ on 
$\Omega^{0,1}_{\mathrm{c}}(\Sigma, L^{p}\otimes E)$. Then since we 
assume that $\ui R^{L}$ is nonnegative, i.e., 
$R^{L}(\omega,\bar{\omega})\geqslant 0$, then, on $(0,1)$-forms,
\begin{equation}
	p(R^{L}(\omega,\bar{\omega}) \, \bar{\omega}^{\ast} \wedge \iota_{\bar{\omega}}-\frac{1}{2}R^{L}(\omega,\bar{\omega}))\geqslant \frac{1}{2}p R^{L}(\omega,\bar{\omega})\geqslant 0.
	\label{eq:2.2.1aug}
\end{equation}
For the points such that $R^{L}$ does not vanish, the above 
term clearly admits a local lower bound growing linearly in $p$.

Under the assumption that $R^{L}$ is semipositive and vanishes up to a finite order, the arguments from \cite[sub-elliptic estimates (2.12) and Proof of Theorem 1]{Marinescu2023} prove that for a compact subset 
$K\subset \Sigma$, there exist constants $C_{1}>0$, $C_{2}>0$ such 
that for $p\gg 1$ and for $s\in 
\Omega^{0,1}_{\mathrm{c}}(\Sigma, L^{p}\otimes E)$ with 
$\supp(s) \subset K$, 
\begin{equation}
	(C_{1}p^{2/\rho_{\Sigma}}-C_{2})\|s\|_{\cLL}\leqslant 
	\Big\|\frac{1}{2}\Delta^{\Lambda^{0,\bullet} \otimes 
	L^p\otimes E}s\Big\|_{\cLL}.
	\label{eq:2.2.2aug}
\end{equation}

We will combine the above considerations to prove Theorem \ref{thm:1.2.1}.

\begin{proof}[Proof of Theorem \ref{thm:1.2.1}] 
	For $s \in \Omega^{0,1}_{\mathrm{c}}(\Sigma, L^{p}\otimes E)$ and a domain $A \subset \Sigma$, set
	\begin{align*}
		\norm{s}_A^2 := \int_A \abs{s}^2_{h_{p}} \omega_\Sigma \, ;
	\end{align*}
	observe that $A \subset B$ implies $\norm{\cdot}_A \leqslant 
	\norm{\cdot}_B$. We fix a compact subset $K$ of $\Sigma$ such that 
	outside of $K$ we have $\ui R^L > c_{K}\omega_{\Sigma}$ with some 
	constant $c_{K}>0$. Then $R^{L}$ can only vanish at the points in $K$. Let $U \subset \Sigma$ be an open relatively compact neighbourhood of 
	$K$. Take smooth functions $\phi_1, \phi_2 : \Sigma \to [0,1]$ such 
	that
	\begin{align*}
		\supp(\phi_1) &\subset U \\
		\supp(\phi_2) &\subset \Sigma \setminus K \, ,
	\end{align*}
	with $\phi_1 \equiv 1$ on $K$ and $\phi_1^2 + \phi_2^2 \equiv 1$ on 
	$\Sigma$. Note that near the punctures, $\phi_{2}$ takes the constant 
	value $1$, then $\norm{\delbar 
	\phi_2}^2_{\mathscr{C}^0}<\infty$, where $\mathscr{C}^{0}$-norm 
	is taken with respect to $g^{T^{*(0,1)}\Sigma}$ for a 
	$(0,1)$-form on $\Sigma$.
	
	The assumption on $(E,h^{E})$ that it is the trivial line bundle 
	near punctures implies that there exists a constant $c_{0}>0$ 
	such that for $x\in\Sigma$, we have
	\begin{equation}
		R^{E}(\omega,\bar{\omega}) \, \bar{\omega}^{\ast} \wedge 
		\iota_{\bar{\omega}}-\frac{1}{2}R^{E}(\omega,\bar{\omega})\geqslant 
		-c_{0}\mathrm{Id}_{T^{\ast(0,1)}\Sigma\otimes L^{p}\otimes E}.
		\label{eq:2.2.3P}
	\end{equation}

	First, we apply \eqref{eq:2.2.2aug} to the sections with support 
	contained in $U$. Then by \eqref{eq:2.1.12cologne}, 
	\eqref{eq:2.2.1aug}, \eqref{eq:2.2.3P} and using the same arguments as in 
	\cite[Proposition 14]{Marinescu2023}, we get that there exist 
	constant $c_1, c_2 \in \R_{>0}$ such that for $s\in 
	\Omega^{0,1}_{\mathrm{c}}(\Sigma, L^{p}\otimes E)$,
	\begin{equation}
		\label{eq:1.2.7}
		(c_1p^{\sfrac{2}{\rho_{\Sigma}}}-c_2) \norm{\phi_1 s}_{U}^2 \leqslant  \norm{\delbar_p^\ast (\phi_1 s)}_{U}^2\, .
	\end{equation}
	
	On the other hand, since $\ui R^{L}(\omega,\bar{\omega})>c_{K}\omega_{\Sigma}$ on 
	the support of $\phi_{2}$, then by \eqref{eq:2.2.3P} and \cite[Theorem 6.1.1, 
	(6.1.7)]{MM07}, there exists
	a constant $c_3 >0$, such that for sufficiently large $p\in 
	\N$
	\begin{equation}
		\label{eq:1.2.8}
		c_3 p \norm{\phi_2s}_{\Sigma \setminus K}^2 \leqslant \norm{\delbar_p^\ast (\phi_2 s)}_{\Sigma \setminus K}^2 \, .
	\end{equation}

	Let $\nabla^{\Lambda^{0,\bullet}\otimes L^{p}\otimes E}$ be the 
	connection on $\Lambda^{\bullet}(T^{\ast(0,1)}\Sigma) \otimes 
	L^p\otimes E$ 
	that is induced by the holomorphic Hermitian connection
	$\nabla^{T^{(1,0)}\Sigma}$ and $\nabla^{L^p\otimes E}$, and let $0 \neq w \in 
	T^{(1,0)}\Sigma$ be a local unit frame, defined on some open 
	set $V$. Because $\Sigma$ is K{\"a}hler, by \cite[Lemma 1.4.4]{MM07}, 
	we have locally $\delbar_p^\ast = 
	-\iota_{\bar{w}}\nabla^{\Lambda^{0,\bullet}\otimes L^p\otimes E}_{\bar{w}}$ for $p \in \N$. As a consequence,
	\begin{equation}
		\label{eq:1.2.10}
  \begin{split}
		&\norm{\delbar_p^\ast (\phi_1 s)}^2_U \leqslant \norm{\delbar 
		\phi_j}^2_{\mathscr{C}^0}\cdot \norm{s}^2_{\cLL} + \norm{\phi_1 
		\delbar_p^\ast s}^2_{\cLL} \, ,\\
  &\norm{\delbar_p^\ast (\phi_2 s)}^2_{\Sigma\setminus K} \leqslant \norm{\delbar 
		\phi_j}^2_{\mathscr{C}^0}\cdot \norm{s}^2_{\cLL} + \norm{\phi_2 
		\delbar_p^\ast s}^2_{\cLL} \, .
    \end{split}
	\end{equation} 
	
	Combining \eqref{eq:1.2.7} - \eqref{eq:1.2.10}, for sufficiently 
	large $p\in\N$,
	\begin{equation}
		\left(\min\big\{c_{1}p^{\sfrac{2}{\rho_{\Sigma}}}-c_{2},\; c_{3}p\big\}-\norm{\delbar 
		\phi_1}^2_{\mathscr{C}^0}-\norm{\delbar 
		\phi_2}^2_{\mathscr{C}^0}\right)\|s\|^{2}_{\cLL}\leqslant 
		\|D_{p}s\|^{2}_{\cLL}.
	\end{equation}
	Since $\rho_{\Sigma}\geqslant 2$, the above inequality infers that there exist constants $C_{1}>0$, $C_{2}>0$ such 
	that for $p\in\N$,
	\begin{equation}
		\begin{split}
			\norm{D_p s}^2_{\cLL} \geqslant 2(C_1 p^{\sfrac{2}{\rho_{\Sigma}}} - C_2) 
			\norm{s}^2_{\cLL} \,.
		\end{split}
	\end{equation}
	This proves \eqref{eq:dirac}.

	Observe that $\mathrm{Spec}(\square_p)= 
	\mathrm{Spec}(\square_p^0)\cup \mathrm{Spec}(\square_p^1) \subset 
	\R_{\geqslant 0}$. For $s\in \Omega^{0,1}_{\mathrm{c}}(\Sigma, 
	L^{p}\otimes E)$,
	\begin{equation}
		\|D_{p}s\|^{2}_{\cLL}=2\langle \square_p s,s\rangle.
	\end{equation}
	Then we get $\mathrm{Spec}(\square_p^1)\subset [C_1 
	p^{\sfrac{2}{\rho_{\Sigma}}} - C_2, +\infty[$, and $H^{1}_{(2)}(\Sigma, 
	L^p\otimes E) = 0$ for $p \gg 0$.
	
	Now take $s \in \Omega^{(0,0)}_{\mathrm{c}}(\Sigma, L^p\otimes E)$, 
	applying \eqref{eq:dirac} to $\delbar_p s$ gives
	\begin{equation}
		\|\square_p^{0} s\|^{2}_{\cLL}\geqslant 
		(C_{1}p^{\sfrac{2}{\rho_{\Sigma}}}-C_{2})\langle\square_p^{0} s,s\rangle.
		\label{eq:1.2.16april}
	\end{equation}
	As a consequence, $\mathrm{Spec}(\square_p^0)\subset \{0\}\cup \left[C_1 
	p^{\sfrac{2}{\rho_{\Sigma}}} - C_2, +\infty\right[$, so that we get \eqref{eq:spec}. This completes the proof of our theorem.
\end{proof}

\subsection{Leading term of Bergman kernel function: a result of Hsiao--Marinescu}

For an arbitrary holomorphic line bundle on a Hermitian manifold, Hsiao and Marinescu \cite{MR3194375} studied the asymptotic expansions of kernel functions of the spectral projections for the low-energy forms. In particular, they refined and generalized the local holomorphic Morse inequalities by Berman \cite{Berman2004}. 

Generally speaking, fix $k_0\geqslant 3$, Hsiao and Marinescu considered the spectral projection $P_{[0,p^{-k_0}]}$ from $\cLL(\Sigma, L^p\otimes E)$ onto the spectral space of the Kodaira Lapacian $\square_p$ associated with the interval $[0,p^{-k_0}]$. Similarly to the Bergman kernel function, let $P_{[0,p^{-k_0}]}(x)$ denote the corresponding spectral kernel function. In \cite[Theorem 1.3 and Corollary 1.4]{MR3194375}, Hsiao and Marinescu obtained a local holomorphic Morse inequality for $P_{[0,p^{-k_0}]}(x)$ as $p \to +\infty$. In particular, the leading term in the expansion was computed.

In the present paper, the spectral gap \eqref{eq:spec} implies that for $p\gg 1$, we have
\begin{equation}
    P_{[0,p^{-k_0}]}=B_p,\;  P_{[0,p^{-k_0}]}(x)=B_p(x), x\in\Sigma.
\end{equation}
Then \cite[Theorem 1.3 and Corollary 1.4]{MR3194375} applies to $B_p(x)$. Note that their results are stated for the sections of $L^p$, but by \cite[Remark 1.11-(II)]{MR3194375}, these conclusions also hold true for $L^p\otimes E$ in our case. 

\begin{theorem}[Hsiao and Marinescu {\cite[Corollary 1.4]{MR3194375}}]\label{thm:HM-thm}
We assume the same conditions on $\Sigma$, $L$ and $E$ as in Theorem \ref{thm:1.2.1}. Recall that the function $\boldsymbol{c}(x)$ on $\Sigma$ is defined in \eqref{eq:1.2.6june}. Then
\begin{enumerate}[label=(\roman*)]
\item Let $\mathbf{1}_{\Sigma_2}$ denote the characteristic function of the open subset $\Sigma_2\subset\Sigma$. For any $x\in \Sigma$, we have
\begin{equation}
    \lim_{ p \to +\infty} \frac{1}{p}B_p(x)=\mathbf{1}_{\Sigma_2}(x) \frac{\boldsymbol{c}(x)}{2\pi}.
    \label{eq:2.3.1june}
\end{equation}
\item Let $K$ be a compact subset of $\Sigma$ and take $\varepsilon >0$, then there exists $p_0\in\N$ such that for any $p\geqslant p_0$, we have for $x\in K$,
\begin{equation}
    B_p(x)\leqslant \left(\varepsilon + \mathbf{1}_{\Sigma_2}(x) \frac{\boldsymbol{c}(x)}{2\pi}\right)p.
    \label{eq:2.3.2june}
\end{equation}
\end{enumerate}
\end{theorem}

It is clear that we can recover the pointwise convergence \eqref{eq:2.3.1june} from our Theorem \ref{thm:bkexpansion}. Moreover, the results stated in Corollary \ref{lem:bkupper} and Proposition \ref{cor:1.2.3} extend the upper bound in \eqref{eq:2.3.2june} for our punctured Riemann surface.

\section{Bergman kernel near the punctures}
In this section, we begin to explain the technique of analytic localization to compute the Bergman kernel $B_p(z,z')$, where the spectral gap in Theorem \ref{thm:1.2.1} plays an essential role. Subsequently, we obtain global off-diagonal estimates for $B_p(z,z')$. Then we will apply the work of Auvray, Ma, and Marinescu \cite{AMM16, AMM21, AMM22} to get the asymptotic expansion of the Bergman kernel function $B_p(z)$ when $z$ is near the punctures. The near-diagonal expansion of $B_p$ and the proof of Theorem \ref{thm:bkexpansion} will be given in the next section.

We introduce the following notation. 
For $m \in \N$ and $s \in \mathscr{C}^\infty(\Sigma, L^p\otimes E)$, 
$z \in \Sigma$, set
\begin{equation}
	|s|_{\mathscr{C}^m(h_{p})}(z) := \left(|s|_{h_{p}} + 
	|\nabla^{p,\Sigma}s|_{h_{p},\omega_\Sigma} + \ldots + 
	|(\nabla^{p,\Sigma})^m s|_{h_{p},\omega_\Sigma}\right)(z),
\end{equation}
where $\nabla^{p,\Sigma}$ is the connection on $(T\Sigma)^{\otimes \ell} 
\otimes L^p\otimes E$, for every $\ell \in \Z_{\geqslant 0}$, induced by the 
Levi-Civita connection associated to $\omega_\Sigma$ and the 
Chern connection that corresponds to the metric $h_{p}$, and 
$|\cdot|_{h_{p},\omega_{\Sigma}}$ denotes the Hermitian metric on $(T\Sigma)^{\otimes \ell} 
\otimes L^p\otimes E$ induced by $g^{T\Sigma}$ and $h_{p}$.
Then for any subset $U\subset \Sigma$, define the norm $\| \cdot 
\|_{\mathscr{C}^m(U,h_{p})}$ on $U$ as follows,
\begin{equation}
	\|s\|_{\mathscr{C}^{m}(U,h_{p})}:=\sup_{z\in 
	U}|s|_{\mathscr{C}^m(h_{p})}(z).
	\label{eq:3.3.2dec}
\end{equation}
If $U=\Sigma$, we write simply
$\|s\|_{\mathscr{C}^{m}(h_{p})}:=\|s\|_{\mathscr{C}^{m}(\Sigma,h_{p})}$. 
Similarly, we also define the analogue norms for the sections on 
$\D^{\ast}$, $\Sigma\times\Sigma$, etc.

For $k\geqslant 1$, let $\bfH^{k}(\Sigma, \omega_{\Sigma}, L^{p}\otimes E,h_{p})$ 
denote the Sobolev space of sections of $(L^{p}\otimes E, h_{p})$ that are $\cLL$-integrable up 
to order $k$. For $s\in \bfH^{k}(\Sigma, \omega_{\Sigma}, 
L^{p}\otimes E,h_{p})$, set
\begin{equation}
	\|s\|^{2}_{\bfH^{k}_{p}}:=\int_{\Sigma}\left(|s|^{2}_{h_{p}}(z)+\left|\nabla^{p,\Sigma}s\right|^{2}_{h_{p},\omega_{\Sigma}}(z)+\cdots+\left|(\nabla^{p,\Sigma})^{k}s\right|^{2}_{h_{p},\omega_{\Sigma}}(z)\right)\omega_{\Sigma}(z)<\infty.
	\label{eq:3.3.1P}
\end{equation}

\subsection{Localization of the problem and off-diagonal estimates}\label{ss3.3P}
In this subsection, we explain how to localize the computations for the 
Bergman kernel $B_{p}$ on $\Sigma$ by the technique of analytic 
localization. For this method, we need two key ingredients: the first one is the 
spectral gap, which is already 
given by Theorem \ref{thm:1.2.1} for our case; the second is the elliptic 
estimates for $\square^{0}_{p}$ as $p$ grows (cf. \cite[Lemma 1.6.2]{MM07}), it is clear by the 
definition of $\square^{0}_{p}$ that they hold 
true on any compact subsets of $\Sigma$. Due to the seminal work 
of Auvray, Ma and Marinescu \cite{AMM16, AMM21}, the necessary elliptic 
estimates for $\square^{0}_{p}$ near the punctures were also
established. Finally, using the finite propagation speed for wave
operators, we can localize the computations of $B_{p}(z,z')$ in our case to the
problems well considered in \cite{AMM16, AMM21} (for computations near punctures) and in \cite{MM07}, \cite{Marinescu2023, MS23} (for computations away from punctures).

Now we give more details. We start with an elliptic estimate proved in \cite[Proposition 4.2]{AMM21}. 
Note that in \cite{AMM21}, they take $(E,h^{E})$ to be a trivial
line bundle on $\Sigma$ and assume that $(L,h)$ is uniformly (strictly) positive on $\Sigma$, but with the same model near 
punctures on $\Sigma$, neither the twist by $E$ nor the positivity of $(L,h)$ away from punctures play any role in the proof of this estimate, so that it extends easily to our 
case.
\begin{proposition}[{\cite[Proposition 4.2]{AMM21}}]
	For any $k\in\N^{\ast}$, there exists $C=C(k,h)$ such that 
	for $p\gg 1$ and all $s\in \bfH^{2k}(\Sigma, \omega_{\Sigma}, 
	L^{p}\otimes E,h_{p})$,
	\begin{equation}
		\|s\|^{2}_{\bfH^{2k}_{p}}\leqslant 
		C\sum_{j=0}^{k}p^{4(k-j)}\|(\square^{0}_{p})^{j}s\|^{2}_{\cLL}
		\label{eq:3.3.2P}
	\end{equation}
\end{proposition}

Fix a small $\varepsilon > 0$. Let $\psi:\R \to [0,1]$ be a smooth even function such that
\begin{equation}
	\psi(v) =\begin{cases} 1 &,\  |v| \leqslant \varepsilon /2 \\ 
	0 &,\  |v| \geqslant \varepsilon \end{cases},
	\label{eq:3.3.5P}
\end{equation}
and define
\begin{align*}
	\varphi(a) = \left( \int_{-\infty}^\infty \psi(v) \d v \right)^{-1} \cdot \int_{-\infty}^\infty e^{iva} \psi(v) \d v \, 
\end{align*}
which is an even function with $\varphi(0)=1$ and lies in the 
Schwartz space $\mathcal{S}(\R)$. 	

For $p> 0$, set 
$\varphi_p(s):=\mathbf{1}_{[\frac{1}{2}\sqrt{C_{1}}p^{\sfrac{1}{\rho_{\Sigma}}},\,\infty[}\,(|s|)\varphi(s)$, 
where $C_1$ is the constant in the spectral gap of Theorem \ref{thm:1.2.1}.

Note that $\varphi$ and $\varphi_{p}$ are even functions. We 
consider the
bounded linear operators $\varphi(D_{p})$, $\varphi_{p}(D_{p})$ 
acting on $\mathcal{L}_{2}^{0,0}(\Sigma, L^{p}\otimes E)$ defined via the 
functional calculus of $\square^{0}_{p}$. In particular, we have
\begin{equation}
	\varphi(D_p) = \frac{1}{2\pi} \int_{\R} 
	\cos\left(\xi\sqrt{\square^{0}_{p}}\right) \hat{\varphi}(\xi) \mathrm{d} \xi \, ,
	\label{eq:3.3.6bis}
\end{equation}
where $\hat{\varphi}$ denotes the Fourier transform of $\varphi$ and is a multiple of the function 
$\psi$ defined in \eqref{eq:3.3.5P}. Then for $p\gg 0$ with $C_{1}p^{\sfrac{2}{\rho_{\Sigma}}}-C_{2}\geqslant 
\frac{C_{1}}{4}p^{\sfrac{2}{\rho_{\Sigma}}}$, we have 
\begin{equation}
	\varphi(D_{p})-B_{p}=\varphi_{p}(D_{p}).
	\label{eq:3.3.6P}
\end{equation}

Let $\varphi_{p}(D_{p})(z,z')$ denote the Schwartz integral kernel of 
$\varphi_{p}(D_{p})$, which is clearly smooth on 
$\Sigma\times\Sigma$. We have the following estimates as an 
extension of \cite[Proposition 5.3]{AMM21}. Fix $0<r<e^{-1}$, recall that the smooth function $\eta:\Sigma \longrightarrow [ 1,\,\infty[$ is such that $\eta(z)=|\log{|z|^{2}}|$ for $z\in \mathbb{D}^{\ast}_{r}$ near each punctures.

\begin{proposition}\label{prop:3.3.2}
	For $\ell,\ m\geqslant 0$, $\gamma>\frac{1}{2}$, there exists 
	$C_{\ell,m,\gamma}>0$ such that for any $p>1$, we have
	\begin{equation}
		\left\|\eta(z)^{-\gamma}\eta(z')^{-\gamma}\varphi_{p}(D_{p})(z,z')\right\|_{\mathscr{C}^{m}(h_{p})}\leqslant C_{\ell,m,\gamma}p^{-\ell}.
		\label{eq:3.3.7P}
	\end{equation}
\end{proposition}
\begin{proof}
	Note that $\varphi(s)$ when is a Schwartz function on $\R$, then for 
	any $k\in\bN$, there exists $M_{k}>0$ such that for $s\in\R$,
	\begin{equation}
		|s^{k}\varphi(s)|\leqslant M_{k}.
	\end{equation}
	Then 
	\begin{equation}
		|\varphi_{p}(s)|\leqslant M_{k}\left(\frac{4}{C_{1}}\right)^{k/2} 
		p^{-\sfrac{k}{\rho_{\Sigma}}}.
		\label{eq:3.3.10P}
	\end{equation}
	Combining \eqref{eq:3.3.10P} with the estimate \eqref{eq:3.3.2P} and the definition of 
	$\varphi_{p}(D_{p})$, we conclude that for any $k, \ell \in \bN$, 
	there exists $C_{k,\ell}>0$ such that for $s\in 
	\mathcal{L}^{0,0}_{2}(\Omega, L^{p}\otimes E)$,
	\begin{equation}
		\|\varphi_{p}(D_{p})s\|_{\bfH^{k}_{p}}\leqslant 
		C_{k,\ell}p^{-\ell}\|s\|_{\cLL}.
	\end{equation}
	Using the above inequality, the proof of \eqref{eq:3.3.7P} 
	follows from the same arguments given in the proof of 
	\cite[Proposition 5.3]{AMM21}, which also need the Sobolev 
	embeddings \cite[Lemma 2.6]{AMM21} for the sections on 
	$\Sigma$ and $\Sigma\times\Sigma$.
\end{proof}

Now Proposition \ref{Prop: off-diag-Berg} is a consequence of Proposition \ref{prop:3.3.2}.
\begin{proof}[Proof of Proposition \ref{Prop: off-diag-Berg}]
We take $\varepsilon$ in \eqref{eq:3.3.5P} the same as fixed one in Proposition \ref{Prop: off-diag-Berg}. By \eqref{eq:2.1.12cologne}, the second order term 
	of $\square^{0}_{p}$ is 
	$\frac{1}{2}\Delta^{\Lambda^{0,\bullet}\otimes L^{p}\otimes E}$. Thus by 
	the finite propagation speed for the wave operators (cf. 
	\cite[Appendix Theorem D.2.1]{MM07}) in \eqref{eq:3.3.6bis} and our assumptions on $\psi$ in \eqref{eq:3.3.5P}, we get that for $z\in \Sigma$, 
	the support of $\varphi(D_{p})(z,\cdot)$ is included in 
	$\mathbb{B}^{\Sigma}(z,\frac{\varepsilon}{\sqrt{2}})$, and
	$\varphi(D_{p})(z,\cdot)$ depends only on the restriction of 
	$\square^{0}_{p}$ on 
	$\mathbb{B}^{\Sigma}(z,\frac{\varepsilon}{\sqrt{2}})$. In 
	particular, if $z,z'\in \Sigma$ are such that $d(z,z')\geqslant 
	\varepsilon$, then
	\begin{equation}
		\varphi(D_{p})(z,z')=0,
	\end{equation}
	so that \eqref{eq:3.3.8P} follows from \eqref{eq:3.3.6P} and 
	\eqref{eq:3.3.7P}. This completes our proof.
\end{proof}

\subsection{Bergman kernel for Poincar\'{e} punctured unit 
disc}\label{ss3.1P}
The Bergman kernel for Poincar\'{e} punctured unit disc is our model 
for the Bergman kernel $B_{p}$ near the punctures of $\Sigma$, which 
is also a central object studied by Auvray--Ma--Marinescu in \cite{AMM16, AMM21}. Now we 
recall the main results proved in \cite[Section 3]{AMM21}.

We consider the Poincar\'{e} punctured unit disc as follows,
\begin{align*}
	(\D^\ast, \omega_{\D^\ast}, \underline{\C}, h_{\D^\ast}) \, ,
\end{align*}
where $h_{\D^\ast} = |\log(|z|^2)|h_0^\C$ with $h_0^\C$ the flat 
Hermitian metric on the trivial line bundle $\underline{\C} \to 
\D^\ast$. Let $z\in\D^{\ast}$ denote the natural coordinate.

For $p\in\bN^{\ast}$, consider the Hermitian metric 
$h_{p,\D^{\ast}}:=|\log(|z|^2)|^{p}h_0^{\C}$ on 
$\underline{\C}$. Define
\begin{equation}
	H^{p}_{(2)}(\D^{\ast}):=H^{0}_{(2)}(\D^\ast, 
	\omega_{\D^\ast}, \underline{\C}, h_{p,\D^\ast}),
\end{equation}
to be the space of $\cLL$-integrable holomorphic functions on 
$\D^{\ast}$ (with respect to the Hermitian metric $h_{p,\D^{\ast}}$).
We denote by $B^{\D^\ast}_p$ the corresponding Bergman 
kernel.

By \cite[Lemma 3.1]{AMM21}, for $p\geqslant 2$, a canonical orthonormal basis of 
$H^{p}_{(2)}(\D^{\ast})$ is given as follows
\begin{equation}
	\left\{\left(\frac{\ell^{p-1}}{2\pi 
	(p-2)!}\right)^{1/2}z^{\ell}\,:\, \ell\in\bN^{\ast}\right\}.
\end{equation}
Then for $p\geqslant 2$, $z,z'\in\D^{\ast}$, we have
\begin{equation}
	B^{\D^{\ast}}_{p}(z,z')=\frac{\left|\log(|z'|^{2})\right|^{p}}{2\pi(p-2)!}\sum_{\ell=1}^{\infty} \ell^{p-1}z^{\ell}\bar{z'}^{\ell}.
\end{equation}
Then the Bergman kernel function has the formula as follows
\begin{equation}
	B^{\D^{\ast}}_{p}(z)=\frac{\left|\log(|z|^{2})\right|^{p}}{2\pi(p-2)!}\sum_{\ell=1}^{\infty} \ell^{p-1}|z|^{2\ell}.
	\label{eq:3.1.4P}
\end{equation}

More explicit evaluations are worked out in \cite[Section 3]{AMM21} for the right-hand side of \eqref{eq:3.1.4P}. In \cite[Proposition 
3.3]{AMM21}, they proved that for any $0<a<1$ and any $m\geqslant 0$, 
there exists $c=c(a)>0$ such that 
\begin{equation}
\label{eq:3.1.5}
	\left\|B_{p}^{\D^{\ast}}(z)-\frac{p-1}{2\pi}\right\|_{\mathscr{C}^{m}(\{a\leqslant |z|<1\},\omega_{\D^{\ast}})}=\mathcal{O}(e^{-cp}),\,\text{as}\, p \to +\infty.
\end{equation}
More generally, for $0<a<1$ and $0<\gamma<\frac{1}{2}$, there exists 
$c=c(a,\gamma)>0$ such that 
\begin{equation}
	\left\|B_{p}^{\D^{\ast}}(z)-\frac{p-1}{2\pi}\right\|_{\mathscr{C}^{m}(\{ae^{-p^{\gamma}}\leqslant |z|<1\},\omega_{\D^{\ast}})}=\mathcal{O}(e^{-cp^{1-2\gamma}}),\,\text{as}\, p \to \infty.
 \label{eq:24-3.1.6}
\end{equation}

Another seminal result proved by Auvray, Ma and Marinescu is the supremum value of $B_{p}^{\D^{\ast}}(z)$. In \cite[Corollary 
3.6]{AMM21}, they proved that 
\begin{equation}
\label{eq:3.1.7}
	\sup_{z\in\D^{\ast}} 
	B_{p}^{\D^{\ast}}(z)=\left(\frac{p}{2\pi}\right)^{\sfrac{3}{\, 2}}+\mathcal{O}(p).
\end{equation}
Their calculations also showed that the points $z$ where 
$B_{p}^{\D^{\ast}}(z)$ approaches its supremum have exponentially small norm $|z|$ as $p \to \infty$.

\subsection{Bergman kernel expansions near a puncture}
Now we consider the chart $V_{j}$ described in our assumption \ref{itm:2}. Fix $0<r<e^{-1}$; we view $\D^{\ast}_{r}$ as a subset of $V_{j}$ with the local complex coordinate $z_{j}$ on $V_{j}$. Then we have the identification of 
geometric data
\begin{equation}
	(V_{j}, \omega_{\Sigma}, L^{p}\otimes E, 
	h_{p})|_{\D^{\ast}_{r}}\cong (\D^\ast, 
	\omega_{\D^\ast}, \underline{\C}, 
	h_{p,\D^\ast})|_{\D^{\ast}_{r}},
	\label{eq:3.3.3dec}
\end{equation}
where the right-hand side is the Poincar\'{e} punctured unit disc 
described in Subsection \ref{ss3.1P}. Let $\square^{0}_{\D^{\ast},p}$ denote the Kodaira Laplacian operator for the Poincar\'{e} punctured unit disc acting on $\mathcal{L}^{0,0}_{2}(\D^{\ast},\omega_{\D^{\ast}},\underline{\C}, 
h_{p,\D^{\ast}})$. Then restricting to $\D^{\ast}_{r}$, $\square^{0}_{\D^{\ast},p}$ coincides with operator $\square^{0}_{p}$.

Note that by \cite[Corollary 5.2]{AMM21}, 
$\square^{0}_{\D^{\ast},p}$ has a spectral gap, i.e.\ ,  there exists 
$C'>0$ such that for $p\gg 0$,
\begin{equation}
	\mathrm{Spec}(\square^{0}_{\D^{\ast},p})\subset\{0\}\cap 
	[C'p,+\infty[.
\end{equation}
Then for $\square^{0}_{\D^{\ast},p}$, we can proceed as in Subsection \ref{ss3.3P}. More precisely, fix $0<\varepsilon<\frac{r}{2}$ to 
define $\psi$ in \eqref{eq:3.3.5P} and the corresponding function 
$\varphi$. Then for $p\geqslant 1$,
\begin{equation}
	\varphi(D_{\D^{\ast},p})-B^{\D^{\ast}}_{p}=\varphi_{p}(D_{\D^{\ast},p}).
	\label{eq:4.1.3S}
\end{equation}

By the finite propagation speed, as explained in the proof of Proposition \ref{prop:3.3.2}, for $z,z'\in \D^{\ast}_{r/2}$, we have
\begin{equation}
	\varphi(D_{\D^{\ast},p})(z,z')=\varphi(D_{p})(z,z').
\end{equation}
Therefore, on $\D^{\ast}_{r/2}\times\D^{\ast}_{r/2}$, we have
\begin{equation}
	B_{p}(z,z')-B^{\D^{\ast}}_{p}(z,z')=\varphi_{p}(D_{\D^{\ast},p})(z,z')-\varphi_{p}(D_{p})(z,z').
	\label{eq:4.1.5X}
\end{equation}
Note that, in fact, both terms in the right-hand side of \eqref{eq:4.1.5X} 
satisfy the estimate \eqref{eq:3.3.7P} on 
$\D^{\ast}_{r/2}\times\D^{\ast}_{r/2}$. Then we can proceed as in 
\cite[Section 6]{AMM21} since the computations are local, we see 
that the results of \cite[Theorems 1.1 \& 1.2]{AMM21} still holds in 
our setting. More precisely, we have the following 
results.
\begin{theorem}[{\cite[Theorems 1.1 \& 
	1.2]{AMM21}}]\label{thm:3.1}
	Fix any $\ell, m \in \N_{\geqslant 0}$.
	For any $\alpha>0$, there exists a constant $C = 
	C(\ell,m,\alpha) > 0$ such that on $\D^{\ast}_{r/2}\times\D^{\ast}_{r/2}$
	\begin{equation}
		\left|B_p(z,z') - B^{\D^\ast}_p(z,z')\right|_{\mathscr{C}^m} \leqslant 
		Cp^{-\ell}\left|\log(|z|^2)\right|^{-\alpha}\left|\log(|z'|^2)\right|^{-\alpha}.
	\end{equation}
	Moreover, for every $\delta > 0$, there exists a constant $C' = 
	C'(\ell,m,\delta) > 0$, such that for all $p \in \Z_{>0}$ and $z_{j} 
	\in \D^{\ast}_{r/2}\,$,
	\begin{equation}
		\label{eq:bk-on-puncture}
		\left|B_p - B^{\D^\ast}_p\right|_{\mathscr{C}^m}(z_j) \leqslant 
		C'p^{-\ell}\left|\log(|z_j|^2)\right|^{-\delta}.
	\end{equation}
\end{theorem}

The behavior of $B^{\D^{\ast}}_{p}$ has been described in Subsection 
\ref{ss3.1P}, combining with the above theorem, we get the asymptotic 
expansion of $B_{p}$ on $\D^{\ast}_{r/2}$ as $p \to +\infty$.

\section{Bergman kernel expansion on $\Sigma$ for semipositive line bundles}\label{ss:bkestimates}
In addition to the off-diagonal estimates in Proposition 
\ref{Prop: off-diag-Berg}, we continue to study the near-diagonal expansion 
of $B_{p}$ via the local models that will 
be described explicitly in Subsection \ref{ss3.2P}. Then we can proceed 
as in \cite[Sections 4.1 \& 4.2]{MM07} to conclude the desired expansions. Finally, we will give the proofs of Theorem \ref{thm:bkexpansion}, Corollary \ref{lem:bkupper}, and Proposition \ref{cor:1.2.3}.

\subsection{Model Dirac and Kodaira Laplacian operators on 
$\C$}\label{ss3.2P}
Alongside the Kodaira Laplacians of our interest, we need to 
introduce certain model operators which play an important role in our 
calculations. We always equip $\mathbb{R}^2$ with the standard Euclidean metric and the standard complex 
structure such that $\R^{2}\cong \C$. Let $z=x+\ui y\in\C$ denote the usual 
complex coordinate, and let $\{e_{1}:=\frac{\partial}{\partial 
x},\ e_{2}=\frac{\partial}{\partial y}\}$ be the standard Euclidean 
basis of $\R^{2}$. Now fix an even integer $\rho'\geqslant 2$.

Let $R$ be a non-trivial $(1,1)$-form on $\mathbb{R}^2$ whose 
coefficient with respect to the frame $\dz\wedge \dzb$ is given by a real nonnegative homogeneous polynomial of degree ${\rho'-2}$.

We define a smooth $1$-form $a^{R} \in \Omega^1(\mathbb{R}^2)$ by
\begin{equation}
	a^{R}_{v_1}(v_2) := \int_0^1 R_{t v_1}(v_2, t v_1) 
	\mathrm{d} t \, , 
\end{equation}
where $v_{1}\in\R^{2}$ and $v_{2}\in T_{v_{1}}\R^{2}\cong \R^{2}$.
Set 
\begin{equation}
	\nabla^{R} = \mathrm{d} - a^{R}\, ,
\end{equation}
it is a unitary connection on the trivial Hermitian line bundle 
$\underline{\C}$ over $\mathbb{R}^2$. In particular, the curvature 
form of 
$\nabla^{R}$ is exactly given by $R$.
Let $\Delta_{R}$ denote the corresponding Bochner Laplacian. 

Take $\delbar$ to be the standard $\delbar$-operator on 
$\mathbb{R}^2\cong \C$; then the $(0,1)$ part of the connection 
$\nabla^{R}$ is $\delbar_\C := \delbar - \left( a^{R} \right)^{0,1}$. 
Let $\delbar_{\C}^{\ast}$ denote the formal adjoint of $\delbar_{\C}$ 
with respect to the standard inner product on $\mathbb{R}^2$.

The following operators are called the model Dirac operator and model Kodaira Laplacian on $\mathbb{R}^2$, corresponding to $R$:
\begin{equation}
	\label{eq:modeloperators}
		{D}_{R} := \sqrt{2} \left( \delbar_\C + \delbar_\C^\ast \right),\; 
		{\square}_{R} := \frac{1}{2} \left( D_{R} \right)^2 \, .
\end{equation}
This model Kodaira Laplacian ${\square}_{R}$ is related to the model Bochner Laplacian by the Lichnerowicz formula
\begin{equation}
	\square_{R} = \frac{1}{2}\Delta_{R} + \frac{1}{2}c\left(R\right)
\end{equation}
with	$c\left(R\right) = R(e_1,e_2) c(e_1)c(e_2).$
We always identify $\Delta_{R}$ and $\square_{R}$ with their unique self-adjoint extensions that act on the $\mathcal{L}^{2}$-sections over 
$\R^{2}$.

Recall that $\square^{0}_{R}$ denotes the restriction of 
$\square_{R}$ on $(0,0)$-sections. In \cite[Proposition 18 in Appendix]{Marinescu2023}, it was proved that there 
exists a constant $c_{R}>0$ such that
\begin{equation}
	\mathrm{Spec}(\square^{0}_{R})\subset\{0\}\cup \left[c_{R},+\infty\right[\;.
	\label{eq:3.2.8}
\end{equation}

Consider the following first-order differential operators
\begin{equation}
	b=-2\frac{\partial}{\partial z}+\frac{1}{\rho'}\ui R(e_{1},e_{2})\bar{z},\; 
	b^{+}=2\frac{\partial}{\partial \bar{z}}+\frac{1}{\rho'}\ui R(e_{1},e_{2})z.
\end{equation}
Then we have
\begin{equation}
	\square^{0}_{R}=\frac{1}{2}bb^{+}.
\end{equation}
Moreover, for $s\in\mathcal{L}^{0,0}_{2}(\R^{2},\underline{\C})$, 
$s\in\ker \square^{0}_{R}$ if and only if $b^{+}s\equiv 0$.

Consider the $\mathcal{L}^{2}$-orthogonal projection
\begin{equation}
	B^{R} : \mathcal{L}^{0,0}_2(\mathbb{R}^2, \underline{\C}) 
	\longrightarrow \ker \square^{0}_{R} \, .
	\label{eq:3.2.9S}
\end{equation}
Let $B^{R}(z,z')$, $z,z'\in \R^{2}$ denote the Schwartz integral kernel 
of the above projection, which is a smooth function on 
$\R^{2}\times\R^{2}$. We also set
\begin{equation}
    B^R(z)=B^R(z,z).
\end{equation}

The following lemma was already known in \cite[the text above Proposition 19]{Marinescu2023}, which can also be viewed as a consequence of the lower bound for the Bergman kernel proved by Catlin \cite{Cat89} by considering the local models. Here we also give a direct proof to shed light on the space $\ker \square^{0}_{R}$.

\begin{lemma}\label{prop:3.2.1}
	For a nontrivial semipositive $R$ as above, $B_{R}$ is an even function, i.e.\,, for 
	$z,z'\in\R^{2}$
	we have $B^{R}(z,z')=B^{R}(-z,-z')$. Moreover,
	\begin{equation}
		B^{R}(0)>0,
		\label{eq:3.2.8P}
	\end{equation}
	and the quantity $B^{R}(0)$ depends on $R$ 
	smoothly (with $R$ having the coefficients as above of a given degree $\rho'-2$).
\end{lemma}
\begin{proof}
	Set $\omega=\frac{1}{\sqrt{2}}(e_{1}-\ui e_{2})$. Note that
	\begin{equation}
		\psi(x,y):=R(\omega,\overline{\omega})=\ui R(e_{1},e_{2})
	\end{equation}
	is, by our assumption, a real homogeneous nonnegative polynomial in $x,y$ of degree $\rho'-2$. In particular, it is an even function 
	in $(x,y)\in\R^{2}$. So that we get the even parity for $B^{R}$ by our 
	construction.
	
	Let $\Psi(x,y)$ be a homogeneous polynomial in $x,y$ of degree 
	$\rho'$ such that
	\begin{equation}
		\frac{\partial 
		\Psi}{\partial\bar{z}}(x,y)=\frac{1}{\rho'}\psi(x,y)z.
	\end{equation}
	Note that for any fixed $\lambda\in\C$, $\Psi+\lambda z^{\rho'}$ also 
	satisfies the above equation. Moreover, we have
	\begin{equation}
		-\frac{1}{2}\Delta^{\R^{2}}\Re(\Psi)=\psi(x,y)\geqslant 0,
	\end{equation}
	where $\Delta^{\R^{2}}=-(\frac{\partial^{2}}{\partial 
	x^{2}}+\frac{\partial^{2}}{\partial y^{2}})$. 
	The real part $\varphi:=\Re(\Psi)$ is a subharmonic, non-harmonic real homogeneous polynomial in $x,y$ of degree $\rho'$.
	
	A straightforward observation is as follows: if $g$ is an entire 
	function on $\C$ such that $|g|^{2}e^{-\varphi}$ is integrable 
	on $\C$ (with respect to the standard Lebesgue measure), then
	\begin{equation}
		ge^{-\frac{1}{2}\Psi}\in \ker \square^{0}_{R}.
	\end{equation}
	
	This way, we change our problem to study the weighted Bergman kernel on $\C$ associated to the real subharmonic 
	function $\frac{1}{2}\varphi$ as in \cite{Chr91}. By \cite[Proposition 
	1.10]{Chr91}, $\ker \square^{0}_{R}$ is an infinite dimensional 
	subspace of $\mathcal{L}^{0,0}_2(\mathbb{R}^2, \underline{\C})$. 
	In particular, there exists a nontrivial entire function $g$ on 
	$\C$ such that $ge^{-\frac{1}{2}\Psi}\in \ker \square^{0}_{R}$. 
	If $g(0)\neq 0$, then $ge^{-\frac{1}{2}\Psi}$ does not vanish at 
	$z=0$. If $g(0)=0$, we write $g(z)=z^{k}f(z)$, where 
	$k\in\N^{\ast}$, $f$ is also an entire function with $f(0)\neq 0$. Then 
	the integrability of $|g|^{2}e^{-\varphi}$ implies that of 
	$|f|^{2}e^{-\varphi}$, so that
	$fe^{-\frac{1}{2}\Psi}\in \ker \square^{0}_{R}$ and it does not 
	vanish at point $z=0$. As a consequence, we have
	\begin{equation}
		B^R(0)=B^{R}(0,0)>0
	\end{equation}
	by the variational characterization of the Bergman kernel.
		
	Analogously to \cite[(4.2.22)]{MM07}, by the spectral gap 
	\eqref{eq:3.2.8}, for $t>0$, we have 
	\begin{equation}
		\exp(-t\square^{0}_{R})-B^{R}=\int_{t}^{\infty}\square^{0}_{R}\exp(-s\square^{0}_{R})\mathrm{d}s.
		\label{eq:4.1.16P}
	\end{equation}
	Then
	\begin{equation}
		B^{R}(0,0)=\exp(-t\square^{0}_{R})(0,0)-\int_{t}^{\infty}\{\square^{0}_{R}\exp(-s\square^{0}_{R})\}(0,0)\mathrm{d}s.
		\label{eq:3.2.17paris}
	\end{equation}
	Now we replace $R$ by a smooth family of non-trivial $(1,1)$-forms 
	on $\R^{2}$ whose coefficients with respect to $\dz\wedge\dzb$ 
	are given by nonnegative real homogeneous polynomials in $x,y$ of degree 
	${\rho'-2}$. 
	Then locally in the parametrization space for this family $R$, the 
	spectral gaps $c_{R}$ in \eqref{eq:3.2.8}, as $R$ varies, admit a 
	uniform lower bound $c>0$ (see \cite[Appendix: Proposition 18]{Marinescu2023}). Combining with the smooth dependence of the heat 
	kernels of $\square^{0}_{R}$ on $R$ (see Duhamel's formula 
	\cite[Theorem 2.48]{BGV04}), 
	$\int_{t}^{\infty}\{\square^{0}_{R}\exp(-s\square^{0}_{R})\}(0,0)\mathrm{d}s$ depends 
	continuously on $R$ for any given $t>0$. As a consequence of \eqref{eq:3.2.17paris}, 
	we conclude that $B^{R}(0,0)$ depends smoothly on $R$. This way, 
	we complete our proof of the lemma.
\end{proof}

\begin{example}
	We consider a simple but nontrivial example $R(x,y)=y^{2}\dz\wedge\dzb$, $\rho'=4$, then we can rewrite it as
	\begin{equation}
		R(x,y)=-2\ui y^{2}\dx\wedge \dy.
	\end{equation}
	Then
	\begin{equation}
		a^{R}_{z}:=\int_{0}^{1}t^{3}(2\ui y^{2}x\dy - 
		2\ui y^{3}\dx)dt=\frac{\ui}{2}y^{2}(x\dy-y\dx),
	\end{equation}
	and
	\begin{equation}
		(a^{R})^{0,1}_{z}=-\frac{1}{4}y^{2}z \dzb.
	\end{equation}
	
	An explicit computation shows that 
	$\delbar^{\ast}_{\C}=-2\iota_{\frac{\partial}{\partial \bar{z}}}\frac{\partial}{\partial 
	z}+\frac{1}{2}y^{2}\bar{z}\iota_{\frac{\partial}{\partial \bar{z}}}$, 
	and that
	\begin{equation}
		\begin{split}
			\square_{R}=&\frac{1}{2}\Delta^{\R^{2}}
			-\frac{1}{2}y^{2}(z\frac{\partial}{\partial 
			z}-\bar{z}\frac{\partial}{\partial 
			\bar{z}})+\frac{\ui}{2}xy\\
			&+\frac{1}{8}y^{4}|z|^{2}-y^{2}+2y^{2}\dzb\wedge\iota_{\frac{\partial}{\partial\bar{z}}}.
		\end{split}
	\end{equation}
	Note that the 
	differential operator
	\begin{equation}
		-\frac{1}{2}y^{2}(z\frac{\partial}{\partial 
		z}-\bar{z}\frac{\partial}{\partial 
		\bar{z}})+\frac{\ui}{2}xy=\frac{\ui}{2}y^{2}(y\frac{\partial}{\partial x}- x\frac{\partial}{\partial y})+\frac{\ui}{2}xy
		\label{eq:3.2.13}
	\end{equation}
	is formally self-adjoint with respect to the standard 
	$\mathcal{L}^{2}$-metric on the functions over $\R^{2}$.
	
	In this example, we have
	\begin{equation}
		b=-2\frac{\partial}{\partial z}+\frac{1}{2}y^{2}\bar{z},\; 
		b^{+}=2\frac{\partial}{\partial \bar{z}}+\frac{1}{2}y^{2}z.
	\end{equation}
	Then
	\begin{equation}
		\square^{0}_{R}=\frac{1}{2}bb^{+}.
	\end{equation}
	
	Note that
	\begin{equation}
		\Re\{|z|^{4}-|z|^{2}z^{2}-\frac{1}{3}|z|^{2}\bar{z}^{2}+\frac{1}{2}z^{4}\}\geqslant \frac{1}{24}x^{4}+\frac{1}{6}y^{4}.
		\label{eq:3.2.16}
	\end{equation}
	Consider the following $\mathcal{L}^{2}$-function on $\C$
	\begin{equation}
		f(z)=\exp\left\{-\frac{1}{16}\left(|z|^{4}-|z|^{2}z^{2}-\frac{1}{3}|z|^{2}\bar{z}^{2}+\frac{1}{2}z^{4}\right)\right\}.
	\end{equation}
	We have $f(0)=1$, and $f\in \ker \square^{0}_{R}$. Moreover, we have
	\begin{equation}
		B^{R}(0)\geqslant\frac{1}{\|f\|_{\cLL}}.
	\end{equation}
\end{example}


\subsection{Construction of local models}
This subsection is a continuation of Subsection \ref{ss3.3P} on the technique of analytical localization, and we will use the same notation as introduced in Subsection \ref{ss3.3P}. In order to compute the asymptotic
expansion of $B_{p}(z)$ as $p \to +\infty$, we need to construct a model Kodaira Laplacian associated with the local geometry near $z$. The machinery of the construction was explained in detail in \cite[Sections 1.6 \& 4.1]{MM07}, and for a compact Riemann surface equipped with a semipositive line bundle, Marinescu and Savale already used this construction in \cite{Marinescu2023, MS23}. In the sequel, we will give more details in order to work out more explicitly the near-diagonal expansions of $B_p$.

Note that $(\Sigma, g^{T\Sigma})$ is complete and hence by the Hopf-Rinow theorem geodesically complete. Thus the exponential map
\begin{align*}
	T_{z}\Sigma \ni Z \mapsto \mathrm{exp}^\Sigma_{z}(Z) \in \Sigma \,
\end{align*}
is well-defined for all $z \in \Sigma$. For an open subset $U \subset 
\Sigma$, set
\begin{equation}
\begin{split}
    	\inj^{U} := \inf_{z \in U} \sup\{\varepsilon >0\;:\;\exp^{U}_z 
	\text{ is a diffeomorphism of } \qquad &\\
 \mathbb{B}^{T_z 
	\Sigma}(0,\varepsilon)
 \text{ onto its image in } U\},&
\end{split}
 \label{eq:4.2.1june}
\end{equation}
which is called the injectivity radius of $U$. If $U$ contains any punctures, we 
always have $\inj^{U}=0$ since the injective radius of a point $z\in 
U$ goes to $0$ as $z$ approaches any puncture in $U$. If $U$ is 
relatively compact in $\Sigma$, then $\inj^{U}>0$.

Fix a point $z_0 \in \Sigma$ and fix an 
open neighborhood $U_0 \subset \Sigma$ of $z_0$ that is
relatively compact in $\Sigma$. Hence $\inj^{U_{0}} > 0$. Let 
$\{e_1,e_2\}$, $\{\mathfrak{e}\}$, and $\{\mathfrak{f}\}$ be 
orthonormal bases for $T_{z_0}\Sigma$, $E_{z_0}$ and $L_{z_0}$ 
respectively, and let $\{w =\frac{1}{\sqrt{2}}(e_1-\ui e_2)\}$ be 
an orthonormal basis for $T^{(1,0)}_{z_0}\Sigma$. Fix some 
$\varepsilon < \inj^{U_{0}} / 4$ such that the vanishing order of $R^{L}$ 
on $\bB^{\Sigma}(z_{0},4\varepsilon)$ is at most $\rho_{z_{0}}-2$. 
Since $\varepsilon$ does not exceed the injectivity radius of $U_{0}$, the exponential map
\begin{equation}
	\label{eq:exponantialmap}
	T_{z_0}\Sigma \supset \mathbb{B}^{T_{z_0}\Sigma}(0,4\varepsilon) 
	\ni Z \mapsto \mathrm{exp}^\Sigma_{z_0}(Z) \in 
	\mathbb{B}^\Sigma({z_0}, 4\varepsilon) \subset \Sigma
\end{equation}
is a diffeomorphism of open balls; it yields a local chart via
\begin{equation}
	\label{eq:identification}
	\R^2 \ni (Z_1,Z_2) \longmapsto Z_1 e_1 + Z_2 e_2 \in T_{z_0}\Sigma \, ,
\end{equation}
called the normal coordinate system (centered at $z_0$).

We always identify $\bB^{T_{z_0}\Sigma} (0,4\varepsilon)$ with 
$\bB^\Sigma(z_0,4\varepsilon)$ via \eqref{eq:exponantialmap}. For $Z \in 
\bB^{T_{z_0}\Sigma}(0,4\varepsilon)$ we identify $L_Z, E_Z$ and 
$\Lambda^{\bullet}(T^{\ast(0,1)}_Z \Sigma)$ to $L_{z_0}, E_{z_0}$ and 
$\Lambda^{\bullet}(T^{\ast(0,1)}_{z_0} \Sigma)$, respectively, by parallel 
transport with respect to $\nabla^L, \nabla^E$ and 
$\nabla^{\Lambda^{\bullet}(T^{\ast(0,1)} \Sigma)}$ along $\gamma_Z: [0,1] \ni u \mapsto 
\mathrm{exp}_{z_0}^\Sigma(uZ)$. This way, we trivilize the 
bundles $L$, $E$, $\Lambda^{\bullet}(T^{\ast(0,1)} \Sigma)$ near 
$z_{0}$. In particular, we will still denote by $\{e_{1}, e_{2}\}$, 
$\{\mathfrak{e}\}$, and $\{\mathfrak{f}\}$ the respective 
orthonormal smooth frames of the vector bundles on point $Z$, defined as 
the parallel transports as above of the vectors $\{e_{1}, e_{2}\}$, 
$\{\mathfrak{e}\}$, and $\{\mathfrak{f}\}$ from $z_{0}$.

With the above local trivializations, we write the connection 
$\nabla^{\Lambda^{0,\bullet} \otimes L^{p} \otimes E}$ 
as follows
\begin{equation}
	\nabla^{\Lambda^{0,\bullet} \otimes L^{p} \otimes E} = \mathrm{d} - \left( a^{\Lambda^{0,\bullet}} + p a^L + a^E \right)
\end{equation}
where $\mathrm{d}$ denotes the ordinary differential operator, 
and $a^{\Lambda^{0,\bullet}}, a^E, a^L$ are respectively the 
local connection $1$-forms of $\nabla^{\Lambda^{0,\bullet}}, 
\nabla^E, \nabla^L$ in this trivialization. Note that these 
connection $1$-forms are purely imaginary.

In coordinate $(Z_{1}, Z_{2})$, we write
\begin{equation}
	a^{L}=\sum_{i=1}^{2} a_{i}^{L}\mathrm{d}Z_{i}.
	\label{eq:4.2.4P}
\end{equation}

Let $R_{ij}^L$ 
denote the coefficients of the curvature form $R^{L}$ with respect to the frame $\mathrm{d}Z_{i}\wedge \, \mathrm{d}Z_{j}$, $i,j=1,2$. We have
\begin{equation}
	R^{L}_{11}=R^{L}_{22}\equiv 0,\; R^{L}_{12}=-R^{L}_{21}.
\end{equation}
Then we can write
\begin{equation}
	R^{L}_{Z}=R^{L}_{12,Z} \, \mathrm{d}Z_1 \wedge \mathrm{d}Z_2.
\end{equation}
Similarly, we define $R_{ij,Z}^{\Lambda^{0,\bullet}}$ and $R_{ij, 
Z}^{E}$. Moreover, we have the 
following relations for $Z \in B^{T_{z_0}\Sigma}(0,\varepsilon)$
\begin{equation}
	a_{i,Z}^{L} = \sum_{j=1}^{2}\int_0^1 t Z^j R_{ij, 
			tZ}^{L} \, \mathrm{d}t \,.
\end{equation}
The analogous identities also hold for $a^{\Lambda^{0,\bullet}}, a^E$.

On the other hand, in these normal coordinates, we find that the curvature $R^L$ of $\nabla^L$ has
the following Taylor expansion at the origin
\begin{equation}
	R^L_{Z} = \sum_{|\alpha| = \rho_{z_0}-2} R_{12;\alpha}^L Z^\alpha 
	\mathrm{d}Z_1 \wedge \mathrm{d}Z_2 + 
	\mathcal{O}(|Z|^{\rho_{z_0}-1}) =: R_{0,Z}^L + 
	\mathcal{O}(|Z|^{\rho_{z_0}-1}) \, ,
	\label{eq:4.2.8P}
\end{equation}
where the $(\mathrm{d}Z_1 \wedge \mathrm{d}Z_2)$-coefficient of $R_0^L$ 
is the product of $-\ui$ and a positive homogeneous even polynomial 
of order $\rho_{z_0}-2$ in $Z$.

Now we construct the local model for $B_{p}$ at $z_{0}$.
Set $\Sigma_0 := T_{z_0}\Sigma \cong \R^2$, and let $Z=(Z_{1},Z_{2})$ 
denote the natural coordinate on $\Sigma_{0}$. Let $(L_0,h_{0})$, 
$(E_0,h^{E_{0}})$ denote 
the trivial line bundles on $\Sigma_{0}$ given by $(L_{z_0}, 
h_{z_{0}})$, $(E_{z_0}, h^{E}_{z_{0}})$ respectively. 
We equip $\Sigma_0$ with $J_{0}$ the almost complex structure on 
$\Sigma_{0}$ that coincides with the pullback of the complex structure $J$ on $\Sigma$ 
by the map \eqref{eq:exponantialmap} in $\bB^\Sigma(z_0,2\varepsilon)$, and is equal to $J_{z_0}$ outside  
$\bB^\Sigma(z_0,4\varepsilon)$. Meanwhile, let $g^{T\Sigma_{0}}$ be the 
Riemannian metric on $\Sigma_{0}$ that is compatible with $J_{0}$ 
and that coincides with the Riemannian metric $g^{T\Sigma}$ on 
$\bB^\Sigma(z_0,2\varepsilon)$, 
and equals to $g^{T\Sigma}_{z_{0}}$ outside 
$\bB^\Sigma(z_0,4\varepsilon)$. In fact, $J_{0}$ is integrable, and
the triplet $(\Sigma_{0},J_{0},g^{T\Sigma_{0}})$ becomes a Riemann surface 
equipped with a complete K\"{a}hler metric $\omega_{\Sigma_{0}}$ induced by $g^{T\Sigma_{0}}$.

Let $T^{\ast(0,1)}\Sigma_{0}$ denote the 
anti-holomorphic cotangent bundle of $(\Sigma_{0},J_{0})$, and let
$\widetilde{\nabla}^{\Lambda^{0,\bullet}}$ denote the Hermitian 
connection on $\Lambda^{\bullet}(T^{\ast(0,1)}\Sigma_{0})$ associated 
with the Levi-Civita connection of $(T\Sigma_{0}, g^{T\Sigma_{0}})$. 
Note that on 
$\bB^{T_{z_{0}}\Sigma}(0,2\varepsilon)$, the pair  
$(\Lambda^{\bullet}(T^{\ast(0,1)}\Sigma_{0}), 
\widetilde{\nabla}^{\Lambda^{0,\bullet}})$ coincides with 
$(\Lambda^{\bullet}(T^{\ast(0,1)}\Sigma),\nabla^{\Lambda^{\bullet}(T^{\ast(0,1)} \Sigma)})$ via the 
identification \eqref{eq:exponantialmap}, and outside 
$\bB^{T_{z_{0}}\Sigma}(0,4\varepsilon)$, the connection
$\widetilde{\nabla}^{\Lambda^{0,\bullet}}$ is given by the trivial 
connection on the trivial bundle 
$\Lambda^{\bullet}(T^{\ast(0,1)}_{z_{0}}\Sigma)$. We can always 
trivialize $T^{\ast(0,1)}\Sigma_{0}$ by the parallel transport along 
the geodesic rays starting at $0$, so that for $Z\in\Sigma_{0}$, 
$T^{\ast(0,1)}_{Z}\Sigma_{0}\cong T^{\ast(0,1)}_{z_{0}}\Sigma$.

Fix an even smooth function $\chi \in C^\infty (\R,[0,1])$ with $\chi = 1$ on $[-2,2]$ and $\supp \chi \subset [-4,4]$. 
We defined a nonnegative curvature form as follows, for 
$Z\in\Sigma_{0}$,
\begin{equation}
	\widetilde{R}^{L_0}_{Z} := 
	\chi\left(\frac{|Z|}{\varepsilon}\right) R^{L}_{Z} + \left( 1 - 
	\chi\left(\frac{|Z|}{\varepsilon}\right)\right) R_{0,Z}^{L}\, ,
\end{equation}
where $R^{L}_{0}$ is defined in \eqref{eq:4.2.8P}.
On $\Sigma_{0}$, define a $1$-form
\begin{equation}
	\tilde{a}^{L_{0}} = \sum_{i=1}^{2} \tilde{a}_i^{L_{0}} \, \mathrm{d}Z_{i}, \; \tilde{a}_i^{L_{0}}(Z) := \int_0^1 t Z^j \widetilde{R}_{ij, 
	tZ}^{L_0} \, \mathrm{d}t.
\end{equation}

Then we set
\begin{equation}
	\label{eq:modified-connections-locally}
	\begin{split}
		& \widetilde{\nabla}^{E_{0}} = \mathrm{d} - \chi\left(\frac{|Z|}{\varepsilon}\right) a^E \, ,  \\
		&\widetilde{\nabla}^{L_{0}} = \mathrm{d} - \tilde{a}^{L_{0}} 
		\, .
	\end{split}
\end{equation}
They are Hermitian connections on the line bundle $L_{0}$, $E_{0}$ 
respectively. Moreover, the 
curvature form of $\widetilde{\nabla}^{L_{0}}$ is exactly
$\widetilde{R}^{L_{0}}$.

As in \eqref{eq:1.1.1bis}, we define for $Z\in\Sigma_{0}$,
\begin{equation}
	\tilde{\rho}_Z := 2 + \ord_Z(\widetilde{R}^{L_0}) \, .
\end{equation}
Since both the vanishing order of $R^{L}$ on 
$\bB^{\Sigma}(z_{0},4\varepsilon)$ and the vanishing order 
$R^{L}_{0}$ on $\Sigma_{0}$ are at most $\rho_{z_{0}}-2$, we get
\begin{equation}
	\tilde{\rho}_{Z}\leqslant \rho_{z_{0}}.
\end{equation}
In particular, $\tilde{\rho}_0 = \rho_{z_0}$, and if $\widetilde{R}^{L_0}(Z) \neq 0$, we have $\tilde{\rho}_Z = 
2$.

Under the above setting on $\Sigma_{0}$, we can define the 
corresponding Dirac and Kodaira Laplacian operators. Note that we 
can use the formulae in \eqref{eq:modeloperators}, or equivalently we use 
the connections $\widetilde{\nabla}^{\Lambda^{0,\bullet}}$, 
$\widetilde{\nabla}^{L_{0}}$, $\widetilde{\nabla}^{E_{0}}$ to define 
the Dirac operator $\widetilde{D}_p$ by \eqref{eq:dirac-clifford}.
Then we have the operators
\begin{equation}
	\label{eq:modoperators}
	\begin{split}
		\widetilde{D}_p &: \Omega^{0, \bullet}_{\mathrm{c}}(\Sigma_0, 
		L_0^p\otimes E_{0}) \longrightarrow \Omega^{0, \bullet}_{\mathrm{c}}(\Sigma_0, L_0^p\otimes E_{0})\, , \\
		\widetilde{\square}_p := \frac{1}{2} (\widetilde{D}_p)^2 &: 
		\Omega^{0, \bullet}_{\mathrm{c}}(\Sigma_0, L_0^p\otimes 
		E_{0}) \longrightarrow \Omega^{0, \bullet}_{\mathrm{c}}(\Sigma_0, 
		L_0^p\otimes E_{0})\, .
	\end{split}
\end{equation}
They extend uniquely to self-adjoint operators acting on 
$\cLL$-sections over $\Sigma_{0}$. By 
construction, the differential operators $\widetilde{D}_p$ and 
$\widetilde{\square}_p$ coincide with $D_p$ and $\square_p$ 
respectively on 
$\bB^{T_{z_{0}}\Sigma}(0,2\varepsilon)\cong\bB^\Sigma(z_{0},2\varepsilon)$.

Let $\widetilde{\Delta}^{\Lambda^{0,\bullet}\otimes L^{p}_{0}\otimes 
E_{0}}$ be the Bochner Laplacian associated to the connection 
$\widetilde{\nabla}^{\Lambda^{0,\bullet}\otimes L^{p}_{0}\otimes 
E_{0}}$. Analogous to \eqref{eq:2.1.12cologne}, we have
\begin{equation}
	\begin{split}
		\widetilde{\square}_{p} = 
		&\frac{1}{2}\widetilde{\Delta}^{\Lambda^{0,\bullet} \otimes 
		L^{p}_{0} \otimes E_{0}} + \frac{r^{\Sigma_0}}{4} \, \bar{\omega}^{\ast} \wedge \iota_{\bar{\omega}}\\
		&\qquad + p \left( \widetilde{R}^{L_0}(\omega,\bar{\omega}) \, \bar{\omega}^{\ast}  \wedge \iota_{\bar{\omega}}-\frac{1}{2}\widetilde{R}^{L_0}(\omega,\bar{\omega}) \right) + \left( \widetilde{R}^{E_0}(\omega,\bar{\omega}) \, \bar{\omega}^{\ast} \wedge \iota_{\bar{\omega}}-\frac{1}{2}\widetilde{R}^{E_0}(\omega,\bar{\omega}) \right),
	\end{split}
	\label{eq:modified-lichnerowicz}
\end{equation}
where $\omega$ denote a unit frame of $T^{\ast(1,0)}\Sigma_{0}$, the 
function $r^{\Sigma_0}$ is the scalar curvature of $(\Sigma_0, 
g^{T\Sigma_0})$, and $\widetilde{R}^{E_0}$ is the curvature  
form of $\widetilde{\nabla}^{E_0}$. Furthermore, $r^{\Sigma_{0}}$, 
$R^{E_{0}}$ vanishes identically outside $\bB^{T_{z_{0}}\Sigma}(0, 
4\varepsilon)$. 

By \eqref{eq:modified-lichnerowicz}, $\widetilde{\square}_p$ 
preserves the degree of $\Lambda^{\bullet}(T^{\ast(0,1)}\Sigma)$. For 
$j=0,1$, let $\widetilde{\square}^{j}_p$ denote the restriction of 
$\widetilde{\square}_p$ on 
$\Omega^{0,j}_{(2)}(\Sigma_{0},L_{0}^{p}\otimes E_{0})$. By the same sub-elliptic estimate proved in 
\cite[(4.13)]{Marinescu2023} for $\widetilde{\Delta}^{\Lambda^{0,\bullet} \otimes 
		L^{p}_{0} \otimes E_{0}}$ as an analogue of 
		\eqref{eq:2.2.2aug}, we get that there exist constants 
		$C'_1,\, C'_2 > 0$, such that
\begin{equation}
	\label{eq:1.3.9}
	\begin{split}	
	&\mathrm{Spec}(\widetilde{\square}^{0}_p) \subset \{0\} \cup \left[C'_1 
	p^{\sfrac{2}{\rho_{z_{0}}}} - C'_2, +\infty\right[ \, ,\\
	&\mathrm{Spec}(\widetilde{\square}^{1}_p) \subset \left[C'_1 
	p^{\sfrac{2}{\rho_{z_{0}}}} - C'_2, +\infty\right[ \, .
	\end{split}
\end{equation}

Set 
\begin{equation}
	H^0_{(2)}(\Sigma_0, L^{p}_0 \otimes E_0):= 
	\mathrm{ker}(\widetilde{\square}^{0}_p).
\end{equation}
Consider the orthogonal projection
\begin{equation}
	\widetilde{B}_{z_{0},p} : \mathcal{L}_{2}^{0,0}(\Sigma_0, L^{p}_0 \otimes 
	E_0) \longrightarrow H^0_{(2)}(\Sigma_0, L^{p}_0 \otimes E_0).
\end{equation}
Let $\widetilde{B}_{z_{0},p}(Z,Z')$ denote the Schwartz kernel of 
$\widetilde{B}_{z_{0},p}$ with respect to the volume element induced 
by $g^{T\Sigma_{0}}$. It is clearly smooth on $\Sigma_{0}\times 
\Sigma_{0}$.

Then we can proceed as in Subsection \ref{ss3.3P}, in particular, by 
Proposition \ref{prop:3.3.2}, we get that for $\ell,\ m\geqslant 0$, there exists 
	$C_{\ell,m}>0$ such that for any $p>1$, we have
	\begin{equation}
		\left\|B_{p}(z,z')-\widetilde{B}_{z_{0},p}(z,z')\right\|_{\mathscr{C}^{m}(\bB^{\Sigma}(z_{0},\varepsilon)\times \bB^{\Sigma}(z_{0}, \varepsilon), h_{p})}\leqslant C_{\ell,m,\gamma}p^{-\ell}.
		\label{eq:4.2.19P}
	\end{equation}
In a shorter notation, we will write for the above statement that
\begin{equation}
	B_{p}-\widetilde{B}_{z_{0},p}=\mathcal{O}(p^{-\infty}), \, 
	\text{ on }\, \bB^{\Sigma}(z_{0},\varepsilon)\times 
	\bB^{\Sigma}(z_{0}, \varepsilon).
	\label{eq:4.2.20S}
\end{equation}

\subsection{Near-diagonal expansion of Bergman kernel}\label{ss4.2W}

The next step is to compute the asymptotic expansion of 
$\widetilde{B}_{z_{0},p}$ around $z_{0}$ as $p \to +\infty$, 
where we can apply the standard method via the rescaling technique as 
in \cite[Subsections 4.1.3 - 4.1.5]{MM07}. One difference is that the 
curvature form $\widetilde{R}^{L_{0}}$ has vanishing order 
$\rho_{z_{0}}-2$ at $Z=0$, so that the rescaling factor will be
\begin{equation}
	t=p^{-\sfrac{1}{\rho_{z_{0}}}}.
\end{equation}

Fix a unit vector $e_{L,z_{0}}$ of $(L_{z_{0}},h_{z_{0}})$. This 
way, we always trivialize $L_{0}^{p}$ as $\C$. Similarly for the line bundle $E_0$. Now, we consider the 
operator $\widetilde{\square}^{0}_{p}$, $p\in\N^{\ast}$, as a family 
of differential operators acting on 
$\mathscr{C}^{\infty}(\R^{2}, \C)$. Let
$\langle\cdot,\cdot\rangle_{\cLL}$ denote the $\cLL$ - inner product on $\mathscr{C}^{\infty}(\R^{2}, \C)$ associated with the Riemannian metric $g^{T\Sigma_{0}}$ and 
$h^{E}_{0}$,  then $\widetilde{\square}^{0}_{p}$ is self-adjoint with
respect to this $\cLL$-inner product.

Meanwhile, we can equip $\R^{2}\cong T_{z_{0}}\Sigma$ with the flat 
Riemnnian metric $g^{T_{z_{0}}\Sigma}$, let $\mathrm{dV}_{0}$ denote 
the corresponding volume form. Let $\kappa(Z)$ be the smooth positive 
function on $\R^{2}$ defined by the equation
\begin{equation}
	\omega_{\Sigma_{0}}(Z)=\kappa(Z) \, \mathrm{dV}_{0}(Z).
	\label{eq:4.2.22S}
\end{equation}
Then $\kappa(0)=1$ and for $Z$ outside $\bB(0,4\varepsilon)$, 
$\kappa(Z)=1$.
Let $\langle\cdot,\cdot\rangle_{\cLL,0}$ denote the standard $\cLL$-inner product on 
$\mathscr{C}^{\infty}(\R^{2}, \C)$.

For $s\in \mathscr{C}^{\infty}(\R^{2}, \C)$, $Z\in\R^{2}$, for 
$t=p^{-\sfrac{1}{\rho_{z_{0}}}}\,$, set
\begin{equation}
	\begin{split}
		&(S_{t}s)(Z):=s(Z/t);\\
		&\mathfrak{L}_{t}:=S_{t}^{-1}\kappa^{1/2}t^{2}\widetilde{\square}^{0}_{p} \kappa^{-1/2}S_{t};\\
		&\mathfrak{L}_{0}:=\square^{0}_{R^{L}_{0}},
	\end{split}
	\label{eq:4.2.23S}
\end{equation}
where the operator $\square_{R^{L}_{0}}$ is the model Kodaira 
Laplacian defined in \eqref{eq:modeloperators} acting on 
$\mathscr{C}^{\infty}(\R^{2}, \C)$ associated to the 
$(1,1)$-form $R^{L}_{0}$ given in \eqref{eq:4.2.8P} with 
$\rho'=\rho_{z_{0}}$. Recall that $B^{R^{L}_{0}}(Z,Z')$ denotes the Bergman 
kernel associated to $\square_{R^{L}_{0}}$ defined by \eqref{eq:3.2.9S}. Moreover, by \eqref{eq:modeloperators}, \eqref{eq:4.2.8P} and 
\eqref{eq:4.2.22S}, both $\mathfrak{L}_{t}$, 
$\mathfrak{L}_{0}$ are self-adjoint with respect to the $\cLL$-metric 
$\langle\cdot,\cdot\rangle_{\cLL,0}$.

By \eqref{eq:1.3.9} and \eqref{eq:4.2.23S}, we get that there exist 
constants $\mu_{0}>0$ and
$t_{0}\in\;]0,1]$ such that for $t\in \;]0,t_{0}]$,
\begin{equation}
	\mathrm{Spec}(\mathfrak{L}_{t})\subset\{0\}\cup \left[\mu_{0},+\infty\right[.
	\label{eq:4.2.4S}
\end{equation}
As explained in Subsection \ref{ss3.2P}, $\mathfrak{L}_{0}$ also 
admits a spectral gap with a constant $c_{R^{L}_{0}}>0$. 

Define the orthogonal projection $\mathfrak{B}_{0,t,z_{0}}: 
(\mathcal{L}^{0,0}_{2}(\R^{2}, \C), 
\langle\cdot,\cdot\rangle_{\cLL,0}) \longrightarrow \ker 
\mathfrak{L}_{t}$, and let $\mathfrak{B}_{0,t,z_{0}}(Z,Z')$ denote the 
smooth kernel of $\mathfrak{B}_{0,t,z_{0}}$ with respect to 
$\mathrm{dV}_{0}$. By \eqref{eq:4.2.23S} with 
$t=p^{-\sfrac{1}{\rho_{z_{0}}}} \leqslant t_{0}$, we have
\begin{equation}
	\widetilde{B}_{z_{0},p}(Z,Z')=t^{-2}\kappa^{-\frac{1}{2}}(Z)\mathfrak{B}_{0,t,z_{0}}(Z/t,Z'/t)\kappa^{-\frac{1}{2}}(Z').
	\label{eq:4.2.25P}
\end{equation}

The structure of the differential operator $\mathfrak{L}_{t}$ is 
exactly the same as the rescaled operator defined in 
\cite[(4.1.29)]{MM07}, so that the computations in the proof of 
\cite[Theorem 4.1.7]{MM07} still hold (with the vanishing order 
$\rho_{z_{0}}-2$ of $\widetilde{R}^{L_{0}}$ at $Z=0$). We can 
conclude the analogue results in \cite[Theorem 4.1.7]{MM07} for our $\mathfrak{L}_{t}$, as explained 
in \cite[Subsection 4.1]{Marinescu2023}. More precsiely, there exist polynomials $\mathcal{A}_{i,j,r}$, $\mathcal{B}_{i,r}$, 
$\mathcal{C}_{r}$ ($r\in\N, i,j\in\{1,2\}$) in $Z=(Z_{1},Z_{2})$ with 
the following properties:
\begin{itemize}
	\item their coefficients are polynomials in 
	$R^{T\Sigma}$, $R^{L}$, $R^{E}$ and their derivatives at 
	$z_{0}$ up to order $r+\rho_{z_{0}}-2$;
	\item $\mathcal{A}_{i,j,r}$ is a homogeneous polynomial in $Z$ of 
	degree $\mathrm{deg}_{Z}\ \mathcal{A}_{i,j,r}=r$, we also have
	\begin{equation}
		\mathrm{deg}_{Z}\ \mathcal{B}_{i,r}\leqslant r+\rho_{z_{0}}-1,\; 
		\mathrm{deg}_{Z}\ \mathcal{C}_{r}\leqslant r+2\rho_{z_{0}}-2.
		\label{eq:4.2.26S}
	\end{equation}
	Moreover,
	\begin{equation}
		\mathrm{deg}_{Z}\ \mathcal{B}_{i,r} -(r-1)=\mathrm{deg}_{Z}\ 
		\mathcal{C}_{r} -r=0 \;\mathrm{mod}\; 2;
		\label{eq:4.2.27S}
	\end{equation}
	\item denote
	\begin{equation}
		\mathfrak{O}_{r}= \mathcal{A}_{i,j,r}\frac{\partial^{2}}{\partial 
		Z_{i}\partial Z_{j}} + \mathcal{B}_{i,r}\frac{\partial}{\partial 
		Z_{i}} +\mathcal{C}_{r},
		\label{eq:4.2.28S}
	\end{equation}
	then 
	\begin{equation}
		\mathfrak{L}_{t}=\mathfrak{L}_{0}+\sum_{r=1}^{m} 
		t^{r}\mathfrak{O}_{r}+\mathcal{O}(t^{m+1}).
		\label{eq:4.2.29S}
	\end{equation}
	The reminder term $\mathcal{O}(t^{m+1})$ is a differential 
	operator up to order $2$, and there exists $m'\in \N$ such that 
	for any $k\in\N$, $t<1$, the derivatives of order $\leqslant k$ of the 
	coefficients of $\mathcal{O}(t^{m+1})$ are dominated by 
	$C_{m,k}t^{m+1}(1+|Z|)^{m'}$. Note that since $\mathfrak{L}_{t}$, 
	$\mathfrak{L}_{0}$ are self-adjoint with respect to $\langle\cdot,\cdot\rangle_{\cLL, 
	0}$, so are $\mathfrak{O}_{r}$ and the remainder term $\mathcal{O}(t^{m+1})$ in \eqref{eq:4.2.29S}.
\end{itemize}

\begin{theorem}\label{thm:24-4.3.1}
Fix $\rho_0\in\{2,\ldots, \rho_\Sigma\}$. Let $W: [0,1]\ni s\mapsto W(s)\in \Sigma$ be a smooth path such 
	that $W(s)\in \Sigma_{\rho_{{0}}}$ for all $s\in [0,1]$.
	For $r\in\N$, there exists a smooth function 
	$\mathfrak{F}_{z,r}(Z,Z')$ on $\R^{2}\times\R^{2}$ which is also 
	smooth in $z\in W([0,1])$ such that 
	for any $k,m,m'\in\N$, $q>0$, there exists $C>0$ such that if 
	$p\geqslant 1$, $Z,Z'\in T_{z}\Sigma$, $|Z|,|Z'|\leqslant 
	q/p^{1/\rho_{{0}}}$, 
	\begin{equation}
		\begin{split}
			&\sup_{|\beta|+|\beta'|\leqslant 
			m}\Big\|\frac{\partial^{|\beta|+|\beta'|}}{\partial 
			Z^{\beta}\partial 
			Z^{\prime,\beta'}}\Big(\frac{1}{p^{2/\rho_{{0}}}} 
			B_{p}(\exp_{z}(Z),\exp_{z}(Z'))\\
			&-\sum_{r=0}^{k}\mathfrak{F}_{z,r}(p^{1/\rho_{{0}}}Z,p^{1/\rho_{{0}}}Z')\kappa^{-1/2}(Z)\kappa^{-1/2}(Z')p^{-r/\rho_{{0}}}\Big)\Big\|_{\mathscr{C}^{m'}(W)}\leqslant C p^{-\frac{k-m+1}{\rho_{{0}}}},
		\end{split}
		\label{eq:4.2.30P}
	\end{equation}
	where $\beta,\beta'\in \N^{2}$ are multi-indices, 
	and the norm $\mathscr{C}^{m'}(W([0,1]))$ is taken with respect to the smooth path $s\mapsto 
	W(s)$ since all the objects inside the big bracket of the 
	left-hand side depend smoothly on $z_{0}\in W([0,1])$.
	
	Moreover, we have the following results:
	\begin{enumerate}
		\item for $r=0$,
			\begin{equation}
		\mathfrak{F}_{z,0}(Z,Z')=B^{R^{L}_{0}}_{z}(Z,Z'),
	\end{equation}
 where $R^L_0$ is the model curvature form on $\Sigma_0=T_{z}\Sigma$ given in \eqref{eq:4.2.8P} for the point $z$, and $B^{R^{L}_{0}}_{z}(Z,Z')$ denotes the corresponding model Bergman kernel as in \eqref{eq:3.2.9S};
	\item 	each $\mathfrak{F}_{z,r}(Z,Z')$ defines a linear operator 
	$\mathfrak{F}_{z,r}$ on $\mathcal{L}^{0,0}_{2}(\R^{2}, 
	E_{z})$,  and $\mathfrak{F}_{z,r}$ is computable by a 
	certain algorithm (cf. \cite[Subsection 4.1.7]{MM07}) in terms of $\mathfrak{L}_{0}\;$, $B^{R^{L}_{0}}$, and 
	$\mathfrak{O}_{j}\;$, $j\leqslant r$;
	\item if $r$ is odd, then $\mathfrak{F}_{z,r}(Z,Z')$ is odd 
	function in $(Z,Z')$, in particular, 
	$\mathfrak{F}_{z,r}(0,0)=0$.
	\end{enumerate}
\end{theorem}
\begin{proof}
Note that when we construct the local operators near each point $z$ in the image of the path $W$, that is $W([0,1])\subset \Sigma_{\rho_0}$, we need to choose small number $\varepsilon>0$, as the explanation before \eqref{eq:exponantialmap}, to be such that for $z\in W([0,1])$, the ball $\mathbb{B}^\Sigma(z,4\varepsilon)$ does not intersect with $\Sigma_j$ with $j>\rho_0$. 

Note that for each $z_0\in W([0,1])$, we have $\rho_{z_0}=\rho_0$. The structure of our operator $\mathfrak{L}_{t}$ given in \eqref{eq:4.2.29S} are the same as in 
	\cite[Theorem 4.1.7]{MM07} (except the different bounds on the 
	degrees in $Z$ of $\mathcal{B}_{i,r}$, $\mathcal{C}_{r}$), so 
	that the Sobolev estimates for the resolvent 
	$(\lambda-\mathfrak{L}_{t})^{-1}$ as well as the asymptotic 
	expansions for $\mathfrak{B}_{0,t,z_{0}}$ obtained in \cite[Subsections 
	4.1.4 \& 4.1.5]{MM07} still hold true. In particular, the operators
$\mathfrak{F}_{z_{0},r}$, $r\in\N$, are defined in the same way with 
smooth Schwartz kernels $\mathfrak{F}_{z_{0},r}(Z,Z')$ respectively, 
and $\mathfrak{F}_{z_{0},0}=B^{R^{L}_{0}}$.
Then \eqref{eq:4.2.30P} with $m'=0$ follows from \cite[Theorem 4.1.18]{MM07}, 
\eqref{eq:4.2.19P} and \eqref{eq:4.2.25P} with 
$t=p^{-\sfrac{1}{\rho_{{0}}}}$.

For higher $m'\geqslant 1$, we can see it as follows: if the path $W$ is a 
constant point $z_{0}$, then it is clear that \eqref{eq:4.2.30P} 
holds with $m'\geqslant 1$; if $W$ is not a constant path, with the 
assumption that $W([0,1]) \subset \Sigma_{\rho_{z_{0}}}$, the spectral gaps 
of the modified operators $\widetilde{\square}_{p}$  with $z_{0}\in W([0,1])$ are given by the same power of $p$, so that we can always use the 
same rescaling factor $t=p^{-\sfrac{1}{\rho_{z_{0}}}}$ to construct 
our operators $\mathfrak{L}_{t}$ as a smooth family parametrized by 
$z_{0}\in W([0,1])$. Then we can proceed as in \cite[Proofs of Theorems 
4.1.16 \& 4.1.24]{MM07}
by considering the derivatives of $(\lambda-\mathfrak{L}_{t})^{-k}$ with respect to $s\in 
[0,1]$ via $z_0=W(s)$. Note that the smooth dependence of $B^{R^{L}_{0}}$ on 
$z_{0}\in W([0,1])$ is already proved in Lemma \ref{prop:3.2.1}. In this way,
we conclude \eqref{eq:4.2.30P} with general $m'\in\N$.

Finally, we prove the parity of $\mathfrak{F}_{z_{0},r}$. Consider 
the symmetry $S_{-1}: \R^{2}\ni Z\mapsto -Z\in\R^{2}$. Since the homogeneous polynomial $R^{L}_{0}(\omega,\bar{\omega})$ is even, that is, it is invariant by $S_{-1}$, we get 
that $\mathfrak{F}_{z_{0},0}=B^{R^{L}_{0}}$ is invariant under the $S_{-1}$-conjugation. By the structure of $\mathfrak{O}_{r}$ 
given in \eqref{eq:4.2.26S} - \eqref{eq:4.2.28S}, we get that 
\begin{equation}
	S_{-1}\mathfrak{O}_{r}S_{-1}=(-1)^{r}\mathfrak{O}_{r}.
\end{equation}
Then using the iterative formula for $\mathfrak{F}_{z_{0},r}$ in 
\cite[(4.1.89), (4.1.91)]{MM07}, by induction from $r=0$, we get
\begin{equation}
	S_{-1}\mathfrak{F}_{z_{0},r}S_{-1}=(-1)^{r}\mathfrak{F}_{z_{0},r}.
\end{equation}
In this way, we complete our proof of the theorem.
\end{proof}

In fact, using the heat kernel approach to $B_{p}$ as in 
\cite[Section 4.2]{MM07}, we can improve the expansion 
\eqref{eq:4.2.30P} so that we get an analogue of \cite[Theorem 
4.2.1]{MM07} as follows.

\begin{theorem}\label{thm:4.3.3feb24}
	Fix $\rho_0\in\{2,\ldots, \rho_\Sigma\}$ and let $W: [0,1]\ni s\mapsto W(s)\in \Sigma$ be a smooth path such 
	that $W(s)\in \Sigma_{\rho_{{0}}}$ for all $s\in [0,1]$. There 
	exists $C''>0$ such that for any $k,m,m'\in\N$, $q>0$, there 
	exists $C>0$ such that if 
	$p\geqslant 1$, $Z,Z'\in T_{z}\Sigma$, $z\in W([0,1])$, $|Z|,|Z'|\leqslant 2\varepsilon$, 
	\begin{equation}
		\begin{split}
			&\sup_{|\beta|+|\beta'|\leqslant 
			m}\left\|\frac{\partial^{|\beta|+|\beta'|}}{\partial 
			Z^{\beta}\partial 
			Z^{\prime,\beta'}}\Big(\frac{1}{p^{2/\rho_{{0}}}} 
			B_{p}(\exp_{z}(Z),\exp_{z}(Z')) \right. \\
			&\qquad \qquad \left. -\sum_{r=0}^{k}\mathfrak{F}_{z,r}(p^{1/\rho_{{0}}}Z,p^{1/\rho_{{0}}}Z')\kappa^{-1/2}(Z)\kappa^{-1/2}(Z')p^{-r/\rho_{{0}}}\Big) \right\|_{\mathscr{C}^{m'}(W)}\\
			&\leqslant C 
			p^{-\frac{k-m+1}{\rho_{{0}}}}\left(1+p^{1/\rho_{{0}}}|Z|+p^{1/\rho_{{0}}}|Z'|\right)^{M_{k+1,m,m'}}\exp\left\{-C'' p^{1/\rho_{{0}}}\left|Z-Z'\right|\right\}+\mathcal{O}(p^{-\infty}),
		\end{split}
		\label{eq:4.2.14S}
	\end{equation}
	where
	\begin{equation}
		M_{k+1,m,m'}=2(k+m'+\rho_{{0}}+1)+m.
		\label{eq:4.2.15P}
	\end{equation}
\end{theorem}
\begin{proof}
	This is just a consequence of the results of \cite[Section 
	4.2]{MM07} together with the spectral gap \eqref{eq:4.2.4S}: 
	applying \eqref{eq:4.1.16P} and \eqref{eq:3.2.17paris} to 
	$\mathfrak{L}_{t}$, then we can use the 
	heat kernel estimates to get suitable bounds on 
	$\mathfrak{B}_{0,t,z_{0}}(Z,Z')$. Note 
	that since the vanishing order of $R^{L}_{0}$ at $Z=0$ is 
	$\rho_{{0}}-2$, so that the power of $(1+|Z|+|Z'|)$ in 
	\cite[Theorem 4.2.5]{MM07} is replaced by 
	$2(r+\rho_{{0}}+m')+m$, which gives \eqref{eq:4.2.15P}. At 
	last, we apply \cite[(4.2.32)]{MM07} with 
	$t=p^{-\sfrac{1}{\rho_{{0}}}}$ to conclude this theorem.
\end{proof}

\begin{remark}\label{rmk:4.3.3}
For the case $z \in \Sigma_2$ (i.e. $\ui R^L_{z}>0$) in \eqref{eq:4.2.14S}, the results in \cite[Theorem 4.1.21]{MM07} still hold. In particular, we have a formula 
\begin{equation}
    \mathfrak{F}_{z,r}(Z,Z')=\mathcal{F}_{z,r}(Z,Z')B^{R^{L}_{0}}_{z}(Z,Z'),
    \label{eq:4.3.16-july}
\end{equation}
 where $\mathcal{F}_{z,r}(Z,Z')$ is a polynomial in $Z,Z'$ with degree $\leqslant 3r$, and $B^{R^{L}_{0}}_{z}(Z,Z')$ has the property
\begin{equation}
   | B^{R^{L}_{0}}_{z}(Z,Z') |=\frac{\boldsymbol{c}(z)}{2\pi} \exp\left\{-\frac{\boldsymbol{c}(z)}{4}\left|Z-Z'\right|^2\right\}
   \label{eq:4.3.17-24}
\end{equation}
with $\boldsymbol{c}(z)=\frac{\ui R^L_{z}}{\omega_{\Sigma}(z)}$.
\end{remark}

\begin{remark}
Note that by our assumption on the small number $\varepsilon>0$ taken in the beginning of the proof of Theorem \ref{thm:24-4.3.1}, we have
$$\bigcup_{z\in W([0,1])}\mathbb{B}^\Sigma(z,2\varepsilon)\subset \Sigma_{\leqslant \rho_{0}}.$$
This means that all the points involved in the expansion \eqref{eq:4.2.14S} can only have the vanishing order $\leqslant \rho_{0}$ for $R^L$. 

   When fix a nonzero $Z=Z'$ in \eqref{eq:4.2.14S}, the term $(1+2p^{1/\rho_{{0}}}|Z|)^{M_{k+1,m,m'}}$ is large enough to cover the difference between $\mathcal{O}(p^{2/\rho_0})$ and $\mathcal{O}(p^{2/\rho_Z})$ with possibly $\rho_Z< \rho_0$, so that the result \eqref{eq:4.2.14S} is not useful to obtain the accurate asymptotic expansion of $B_{p}(\exp_{z}(Z),\exp_{z}(Z))$ when $\rho_0>2$.
\end{remark}
%

\subsection{Proofs of 
Theorem \ref{thm:bkexpansion}, Corollary \ref{lem:bkupper}, and Proposition \ref{cor:1.2.3} }\label{ss:4.4june}

Now we prove Theorem \ref{thm:bkexpansion} as a consequence of Theorem \ref{thm:24-4.3.1}.
\begin{proof}[Proof of Theorem \ref{thm:bkexpansion}]
	We take $Z=Z'=0$, $m=0$ in \eqref{eq:4.2.30P}, note that 
	$\mathfrak{F}_{z_{0},2r+1}(0,0)=0$, $r\in\N$, $z_{0}\in W([0,1])$, then we get 
	\eqref{eq:bkexpansion} by setting
	\begin{equation}
		b_{r}(z_{0})=\mathfrak{F}_{z_{0},2r}(0,0),\; z_{0}\in W([0,1]).
	\end{equation}
	For the second part, on $\D^\ast(a_j, 1/4)$, the estimates \eqref{eq:24-3.1.6} and \eqref{eq:bk-on-puncture} hold, from them we conclude \eqref{eq:24-1.2.5}. This way, we complete our proof.
\end{proof}

\begin{proof}[Proof of Corollary \ref{lem:bkupper}]
After fixing $t$ and $\gamma$ as in the corollary, we consider suifficiently large $p\gg 1$ and set
\begin{equation}
\begin{split}
       K_{1,p} &:= \bigcup_{j=1}^N \mathbb{D}(a_j,1/6)\setminus  \mathbb{D}(a_j,te^{-p^\gamma});\\
       K_2 &:= \overline{\Sigma}\setminus \left(\bigcup_j \D(a_j, 1/6) \right).
\end{split}
\end{equation}
Then $\Sigma_{p,t,\gamma}=K_{1,p} \cup K_2$. 

By \eqref{eq:24-1.2.5}, we conclude that the following identity hold uniformly for $x\in K_{1,p}$ as $p \to +\infty$
\begin{equation}
    B_p(x)=\frac{1}{2\pi}(1+o(1))p.
    \label{eq:4.4.3june}
\end{equation}

Now we deal with the points in $K_2$ which is a compact subset of $\Sigma$ independent of $p$. By Theorem \ref{thm:HM-thm}-(ii), taking any sequence $\{\varepsilon_j>0\}_{j\in\N}$ with $\lim_{j \to +\infty} \varepsilon_j=0$, we have an increasing sequence of integers $\{p_j\}_j$ with $p_j \to +\infty$ such that for any $p\geqslant p_j$
\begin{equation}
    \sup_{x\in K_2} B_p(x)\leqslant (C_0+\varepsilon_j) p.
\end{equation}
Then we conclude, as $p \to +\infty$,
\begin{equation}
    \sup_{x\in K_2} B_p(x)\leqslant C_0(1+o(1))p.
\end{equation}

Combining the above result with \eqref{eq:4.4.3june}, we prove this corollary.
\end{proof}

\begin{proof}[Proof of Proposition \ref{cor:1.2.3}]
	Fix $0<r \leqslant e^{-1}$. For $z_j \in V_j\subset \Sigma$ near a puncture, \eqref{eq:bk-on-puncture}, together with \eqref{eq:3.1.5} and \eqref{eq:3.1.7}(see also \cite[Corollary 
3.6]{AMM21}) implies that
\begin{equation}
\label{eq:4.4.1}
		\sup_{|z_j| \leqslant r} B_p(z_j) = \left( \frac{p}{2\pi} \right)^{\sfrac{3}{\,2}} + \mathcal{O}(p) \qquad \text{ as } p \to \infty \, .
	\end{equation}

Away from the punctures, on the compact subset $K:=\overline{\Sigma}\setminus \cup_j \D(a_j, r)$ of $\Sigma$, we apply \eqref{eq:2.3.2june} (from \cite[Corollary 1.4]{MR3194375}) or Corollary \ref{lem:bkupper} to it, then there exists $C>0$ such that
\begin{equation}
    \sup_{x\in K} B_p(x)\leqslant Cp.
    \label{eq:4.4.7june}
\end{equation}

 Combining \eqref{eq:4.4.1} with \eqref{eq:4.4.7june}, we get \eqref{eq:1.2.11june}.
\end{proof}

We can describe the derivatives of the Bergman kernel in a coordinate-free fashion by considering the associated jet-bundles (see \hyperref[ssec:jet]{Appendix}). A pointwise asymptotic expansion also exists for derivatives of the Bergman kernel functions. 

\begin{theorem}\label{thm:4.4.1june}
	For all $\ell \in \N_0\,$, the $\ell$-th jet of the on-diagonal Bergman kernel has a pointwise asymptotic expansion
	\begin{equation}
		j^\ell B_p(x) / j^{\ell-1} B_p(x) = p^{(2+\ell)/\rho_x} \left[ \sum_{j=0}^k c^\ell_j(x) p^{-j/\rho_x} \right] + \mathcal{O}(p^{-(k-\ell-1)/\rho_x})
  \label{eq:4.4.8june}
	\end{equation}
	for all $k \in \N$ with the coefficients $c^\ell_j(x)\in\C$. 
 
The leading term is given by
	\begin{equation}
		c^\ell_0(x) = j^\ell B^{R^L_0}_x(0)/j^{\ell-1} B^{R^L_0}_x(0)
	\end{equation}
	in terms of the $\ell$-th jet of the model Bergman kernel on the tangent space at $x \in \Sigma$ with respect to the geodesic coordinates $Z=(Z_1,Z_2)$ (see also Theorem \ref{thm:24-4.3.1}).
 In particular, if $\ell$ is odd, then $c^\ell_0(x)=0$.
\end{theorem}

\begin{proof}
	This is a consequence of Theorem \ref{thm:24-4.3.1} via taking the Taylor expansion for the Bergman kernel function $B_p(\exp_x(Z)):=B_p(\exp_x(Z),\exp_x(Z))$ in variable $Z$ at $Z=0$. For the leading term, we have
 \begin{equation}
 \begin{split}
         & j^\ell_{Z=0} \left[B^{R^L_0}_x(p^{1/\rho_x}Z)\kappa^{-1}(Z)\right]/j^{\ell-1}_{Z=0}\left[B^{R^L_0}_x(p^{1/\rho_x}Z)\kappa^{-1}(Z)\right]\\
         &=p^{\ell/\rho_x} j^\ell B^{R^L_0}_x(0)/j^{\ell-1} B^{R^L_0}_x(0)+\mathcal{O}_x(p^{(\ell-1)/\rho_x}).
 \end{split}
 \end{equation}
 In this way, we conclude \eqref{eq:4.4.8june} and the formula for $c^\ell_0(x)$. If $\ell$ is odd, using the fact that $B^{R^L_0}_x(Z)$ is an even function (by Lemma \ref{prop:3.2.1}) in $Z$, we get $c^\ell_0(x)=0$. 
\end{proof}

Theorem \ref{thm:4.4.1june} extends \cite[Theorem 3.1]{MS23} for compact Riemann surfaces.

\subsection{Normalized Bergman kernel: proof of Theorem \ref{thm:1.3.1}}
Different from \cite[Theorem 1.8]{Drewitz_2023}, the line bundle $(L,h)$ here is semipositive and hence no longer uniformly positive in $\Sigma$, this is the reason we only make the statement for a subset $U\subset \Sigma_2\,$, see also \cite[Theorem 1.20]{Drewitz:2024aa} for an analogous result of normalized Berezin-Toeplitz kernels.


\begin{proof}[Proof of Theorem \ref{thm:1.3.1}]
    By Theorem \ref{thm:4.3.3feb24}, we see that, for the points where $\ui R^L$ is strictly positive in $U$, the near-diagonal expansions of $B_p(x,y)$ behave the same as in \cite[Theorems 4.2.1 and 6.1.1]{MM07}. Using analogous arguments as in \cite[Subsection 2.3]{Drewitz_2023} and \cite[Subsection 2.4]{Drewitz:2024aa} together with the off-diagonal estimate \eqref{eq:3.3.8P}, we can obtain the estimates in Theorem \ref{thm:1.3.1} - \ref{1.4.1-1} and \ref{1.4.1-2}. Note that instead of $b>\sqrt{16k/\varepsilon_{0}}$ in \cite[Theorem 1.8]{Drewitz_2023}, we improve the condition to $b\geqslant \sqrt{12k/\varepsilon_{0}}$, and here we also state a sharper estimate in Theorem \ref{thm:1.3.1} - \ref{1.4.1-3} for the remainder term $R_p$ than \cite[Theorem 1.8]{Drewitz_2023}. Therefore, we reproduce the proof in detail as follows.

First of all, since $\overline{U}\subset \Sigma_2$, by Theorem \ref{thm:bkexpansion}, there exists a constant $c>0$ such that for all point $x\in U$ and for $p\gg 1$,
\begin{equation}
    B_p(x)=B_p(x,x)\geqslant cp.
    \label{eq:4.5.1-july-24}
\end{equation}

Now we start with a proof of \ref{thm:1.3.1} - \ref{1.4.1-1}. Note that $U$ is relatively compact, so Proposition \ref{Prop: off-diag-Berg} is applicable. Fix $k\geqslant 1$ and let $\varepsilon>0$ be the sufficiently small quantity stated in Proposition \ref{Prop: off-diag-Berg}. Then for
$x,y\in U$ with $\operatorname{dist}(x,y) \geqslant \varepsilon$, we have
\begin{equation}
	|B_{p}(x,y)|\leqslant C_{k,\varepsilon,K}\, p^{-k+1}.
	\label{eq:3.1.1}
\end{equation}

Recall that $\varepsilon_0:=\inf_{x\in U} \bb{c}(x)>0$. Now we fix $b\geqslant \sqrt{12k/\varepsilon_{0}}$, and a large enough $p_{0}\in\mathbb{N}$ such that
\begin{equation}
	b\,\sqrt{\frac{\log{p_{0}}}{p_{0}}}\leqslant \frac{\varepsilon}{2}.
	\label{eq:3.1.2}
\end{equation}

For $p>p_{0}$, if $x,y\in U$ is such that 
$b\sqrt{\log{p}/{p}}\leqslant \operatorname{dist}(x,y)< \varepsilon$, since we work on $\overline{U}\subset \Sigma_2$, we take advantage of the expansion in \eqref{eq:4.2.14S} with the first $2k+1$ terms and with $\rho_0=2$, $m=m'=0$, $x_{0}=x$,  $y=\exp_{x}(Z)$, and $Z\in T_{x}\Sigma$, in order to obtain
\begin{equation}
	\begin{split}
		&\left|\frac{1}{p}B_{p}(x,y)-\sum_{r=0}^{2k}\mathfrak{F}_{x,r}(0,\sqrt{p}Z)\kappa^{-1/2}(Z)p^{-r/2}\right|\\
		&\leqslant 
		Cp^{-k-1/2}(1+\sqrt{p}|Z| )^{4k+6}\exp\left\{-C'\sqrt{p}|Z| \right\}+\mathcal{O}(p^{-k-1}).
	\end{split}
	\label{eq:3.1.3}
\end{equation}

There exists a constant $C_{k}>0$ such that for any $r>0$,
\begin{equation}
	(1+r)^{4k+6}\exp(-C'r)\leqslant C_{k}.
\end{equation}

Note that $|Z| =\operatorname{dist}(x,y)$. By Remark \ref{rmk:4.3.3}, we have the formula \eqref{eq:4.3.16-july} for $\mathfrak{F}_{x,r}$ with the polynomial factor $\mathcal{F}_{x,r}(Z,Z')$, and that the degree of $\mathcal{F}_{x,r}(Z,Z')$ is not greater than $3r$, and the fact that $\varepsilon >|Z| \geqslant b\sqrt{\log{p}/{p}}$, we get for $r=0,\ldots, 2k$,
\begin{equation}
	|\mathfrak{F}_{x,r}(0,\sqrt{p}Z)p^{-r/2}|\leqslant C 
	p^{r}\exp \left\{-\frac{\bb{c}(x)}{4} 
	b^{2}\log{p} \right\},
\end{equation}
where the constant $C=C_U>0$ does not depend on $x\in U$.

Since we take $b \geqslant \sqrt{12k/\varepsilon_{0}}$, then for $r=0,\ldots, 
2k$, we get
\begin{equation}
	\left|p^{r}\exp \left\{-\frac{\bb{c}(x)}{4} 
	b^{2}\log{p} \right\} \right|\leqslant 
	p^{-k}.
	\label{eq:3.1.6}
\end{equation}
Finally, combining \eqref{eq:4.5.1-july-24}--\eqref{eq:3.1.6}, we get the 
desired estimate in Theorem \ref{thm:1.3.1} - \ref{1.4.1-1}.

    Let us prove Theorem \ref{thm:1.3.1} - \ref{1.4.1-2}. Fix $b \geqslant \sqrt{12k/\varepsilon_{0}}$, and we only consider $p\gg 1$.
    Recall that the constant $C_0$ is defined in \eqref{eq:1.2.9-24}, then set
    \begin{equation}
        M_b=\lceil \pi b^2 C_0 \rceil \in \bN.
    \end{equation}
    
    Then for $x\in U\subset \Sigma_2$ and $Z\in T_z\Sigma$ with $|Z|\leqslant b\sqrt{\log {p}/{p}}$, set $y=\exp_x(Z)\in U$, then $\dist(x,y)=|Z|$. Then 
    \begin{equation}
        \exp\left\{\frac{\boldsymbol{c}(x)p}{4}\operatorname{dist}(x,y)^{2} \right\} \leqslant p^{M_b/2}.
        \label{eq:4.5.2-24}
    \end{equation}
    Take the expansion \eqref{eq:4.2.14S} with $\rho_0=2$ and $k=M_b$, $m=m'=0$, we get
  \begin{equation}
        \left|\frac{1}{p} B_{p}(x,y)-\sum_{r=0} ^{M_b}\mathfrak{F}_{z,r}(0,\sqrt{p}Z)\kappa^{-1/2}(Z)p^{-r/2}\right|\leqslant C p^{-\frac{M_b+1}{2}}+\mathcal{O}(p^{-\infty}).
\label{eq:4.5.3-24}	
 \end{equation}
By Remark \ref{rmk:4.3.3}, we get for $r\geqslant 1$,
\begin{equation}
   \exp \left\{\frac{\boldsymbol{c}(x)p}{4}\operatorname{dist}(x,y)^{2} \right\} \left| \mathfrak{F}_{z,r}(0,\sqrt{p}Z)\kappa^{-1/2}(Z)p^{-r/2} \right| \leqslant C_r \left| \log p \right|^{3r/2} p^{-1/2}.
   \label{eq:4.5.4-24}
\end{equation}

Combining \eqref{eq:4.5.2-24} - \eqref{eq:4.5.4-24}, we get
\begin{equation}
	\begin{split}
		\frac{  \exp\left\{\frac{\boldsymbol{c}(x)p}{4}\operatorname{dist}(x,y)^{2}\right\} B_{p}(x,y)}{\sqrt{B_{p}(x)}\sqrt{B_{p}(y)}}&=\frac{\frac{\bb{c}(x)}{2\pi}\kappa^{-1/2}(Z)+\mathcal{O}(p^{-1/2+\varepsilon})}{\sqrt{\frac{\bb{c}(x)}{2\pi}+\mathcal{O}(p^{-1})}\sqrt{\frac{\bb{c}(y)}{2\pi}+\mathcal{O}(p^{-1})}}\\
		&=1+\mathcal{O}(|Z| +p^{-1/2+\varepsilon})\\
&= 1+\mathcal{O}(p^{-1/2+\varepsilon}) \text{ as } p \to +\infty.
	\end{split}
	\label{eq:3.1.8bis}
\end{equation}
The term $\mathcal{O}(p^{-1/2+\varepsilon})$ in the last line of \eqref{eq:3.1.8bis} represents the 
function $R_p$, so Theorem \ref{thm:1.3.1} - \ref{1.4.1-2} and \ref{1.4.1-3} follow.
\end{proof}

Analogously to \cite[Proposition 2.8]{SZ08} and \cite[Lemma 2.13]{Drewitz:2024aa}, we have the following results, and we refer to \cite[Proof of Lemma 2.13]{Drewitz:2024aa} for a proof.
\begin{lemma}\label{lm:2.14feb24}
With the same assumptions in Theorem \ref{thm:1.3.1}, the term $R_p(x,y)$ 
satisfies the following estimate: there exists $C_1=C_1(\varepsilon, U)>0$ 
such that for all sufficiently large $p$, $x,y \in U$ with $\dist(x,y)\leqslant b\sqrt{\log p / p}\,$,
\begin{equation}
|R_p(x,y)|\leqslant C_1 p^{1/2+\varepsilon} \dist(x,y)^2.
\label{eq:2.75feb24}
\end{equation}
For given $k, \ell\in \N$, there exists a sufficiently large $b>0$ 
such that there exists a constant $C_2>0$ such that
for all $x,y\in U$, $\mathrm{dist}(x,y)\geqslant b\sqrt{\log{p}/p}\,$, we have for $p\gg 1$
\begin{equation}
    \left|\nabla^\ell_{x,y}N_{p}(x,y) \right| \leqslant C_2 p^{-k}.
\label{eq:2.76feb24}
\end{equation}
\end{lemma}

\section{Equidistribution and smooth statistics of random zeros}
 Marinescu and Savale  \cite[Theorem 1.4 and Section 6]{MS23} proved a equidistribution result for the zeros of Gaussian random holomorphic sections of the semipositive line bundles over a compact Riemann surface. In this section, we apply our results of Section \ref{ss:bkestimates} to prove a refined equidistribution result for the random zeros of $s_p\in H^0_{(2)}(\Sigma, L^p\otimes E)$. Furthermore, we will follow the work of \cite{SZZ08, SZ08, MR2742043} and \cite{Drewitz_2023, DrLM:2023aa, Drewitz:2024aa} to study the large deviations and smooth statistics of these random zeros.

\subsection{On {$\mathcal{L}^1$}-norm of logarithm of Bergman kernel function}

An important ingredient to study the semi-classical limit of zeros of $\bb{S}_p$ (see Definition \ref{defn:Gaussian}) is to study the function $\log B_p(x)$ as $p \to +\infty$.

For $t\in\;]0,1[\;$, $\gamma\in \;]0,\frac{1}{2}[\;$, as in \eqref{eq:1.2.8-intro-24}, we set
\begin{equation}
    \Sigma_{p,t,\gamma}=\Sigma \setminus \bigcup_{j=1}^{N} 
\mathbb{D}^{\ast}(a_{j},te^{-p^\gamma}).
\end{equation}
We have the following result for the $\mathcal{L}^1$-norm of $\log B_p$ on $\Sigma_{p,t,\gamma}$.
\begin{theorem}\label{thm:5.1.1-24}
    Let $\Sigma$ be a punctured Riemann surface, and let $L$ be a holomorphic line bundle as above such that $L$ carries a singular Hermitian metric $h_{L}$ satisfying conditions \ref{itm:1} and \ref{itm:2}. Let $E$ be a holomorphic line bundle on 
	$\Sigma$ equipped with a smooth Hermitian metric $h^{E}$ such 
	that $(E,h^{E})$ on each chart $V_{j}$ is exactly a trivial Hermitian line bundle. 
 Then for the Bergman kernel functions $B_p(x)$ associated to $H^0_{(2)}(\Sigma, L^p\otimes E)$, there exists a constant $C=C(t,\gamma)>0$ such that for all $p\gg 1$
    \begin{equation}
    \label{eq:log-bk}
        \int_{\Sigma_{p,t,\gamma}} | \log B_p(z)|\omega_\Sigma(x)\leqslant  C \log p.
    \end{equation}
\end{theorem}
\begin{proof}
For a compact Riemann surface with a semipositive line bundle, this theorem follows easily from the uniform two-sided bounds on $B_p$ in \cite[Lemma 3.3]{MS23}, and the analogous arguments, combining with \eqref{eq:24-1.2.5}, shall prove this theorem. But in the sequel, we will sketch a different approach which is independent of the uniform estimates as in \cite[Subsection 3.1]{MS23}.

By Proposition \ref{cor:1.2.3}, there exists a constant $C>0$ such that
\begin{equation}
    \sup_{x\in \Sigma} \log B_p \leqslant \frac{3}{2}\log p +C.
 \end{equation}
Thus, in order to prove \eqref{eq:log-bk}, it remains to bound the negative part of $\log B_p$.

At first, we claim that there exists a smooth Hermitian metric $\tilde{h}$ on $L\to \Sigma$ such that for a small $\varepsilon>0$ and on $\Sigma$, we have
\begin{equation}
    h\leqslant \tilde{h}\,, \quad \ui \widetilde{R}^L\geqslant \varepsilon \omega_\Sigma .
    \label{eq:5.1.3-24-june}
\end{equation}

In fact, since $L$ is positive in $\overline{\Sigma}$, we can always take a smooth Hermitian metric $\hat{h}$ on $L$ such that $\ui \widehat{R}^L>0$ on $\overline{\Sigma}$ (see \cite{MR615130}). For each $z\in \Sigma$, take $e_L(z)$ a nonzero element of $L_z$, then set
\begin{equation}
    \widehat{F}(z):=-\log\frac{|e_L(z)|_{\hat{h}}}{|e_L(z)|_h}.
\end{equation}
Then $\widehat{F}$ is a smooth real function on $\Sigma$ and tends to $+\infty$ at punctures. Then on $\Sigma$, 
\begin{equation}
    \ui \widehat{R}^L=2\ui\partial\db \widehat{F} +\ui R^L >0.
\end{equation}

 Now we modify $\widehat{F}$ to a new function $\widetilde{F}$ such that $\widetilde{F}$ is a smooth function on $\overline{\Sigma}$ with the properties:
\begin{enumerate}
    \item $\max_{z\in\overline{\Sigma}} |\widetilde{F}|\leqslant M_0$, where $M_0\gg 1$ is some constant.
    \item $\ui\partial\db \widetilde{F}\equiv 0$ on each local chart $\{0<|z_j|<r_0\}\subset V_j$, where $0<r_0<e^{-1}$ is given, and $V_j$ is the local chart in the assumption \ref{itm:2}.
    \item $\widetilde{F}=\widehat{F}$ on the subset  $\Sigma\setminus \cup_j\{0<|z_j|<2r_0\}$.
\end{enumerate}
Hence there exists $\delta\geqslant 1$ such that 
\begin{equation}
 \sup_{z\in \Sigma}\left|   \frac{\ui\partial\db \widetilde{F}(z)}{\omega_\Sigma(z)}\right|\leqslant \delta.
\end{equation}

Now we set a new smooth metric on $L\to\Sigma$,
\begin{equation}
    \tilde{h}(\cdot,\cdot)_z:=e^{(-\widetilde{F}(z)+M_0)/2 \delta} h(\cdot,\cdot)_z.
    \label{eq:5.1.6-24-June}
\end{equation}
It is clear that $h\leqslant \tilde{h}$, and we have
\begin{equation}
    \widetilde{R}^L = \frac{1}{2\delta}\ui\partial\db \widetilde{F} +R^L,
\end{equation}
which implies that the metric $\tilde{h}$ satisfies the second condition in \eqref{eq:5.1.3-24-june}.

Moreover, choosing properly $\varepsilon>0$, and fix a large $p_0\in \bN$, we have for $p\geqslant p_0$ and globally on $\Sigma$, 
\begin{equation}
    (p-p_0)\ui R^L +p_0 \ui \widetilde{R}^L \geqslant p_0 \varepsilon \omega_\Sigma.
    \label{eq:5.1.3-24}
\end{equation}

 Let $x \in \Sigma$ and $U_0 \subset \Sigma$ be a small coordinate neighborhood of $x$ on which there exist holomorphic frames $e_L$ of $L \to U_0$ and $e_E$ of $E\to U_0$. Let $\psi, \tilde{\psi}, \psi_E \in \cC^\infty(U_0)$ be the subharmonic weights of $h$, $\tilde{h}$ and $h^E$, respectively, on $U_0$ relative to $e_L$, $e_E$, that is, $|e_L|_{h}^2 = e^{-2\psi}$ and etc. A suitable scalar multiplication of the section $e_L$ allows us to assume that $\psi \leqslant 0$. The condition that $h \leqslant \tilde{h}$ implies $\tilde{\psi} \leqslant \psi$. 
 
 Consider a $p_0$ (that will be chosen momentarily) and write $L^p = L^{p-p_0} \otimes L^{p_0}$. Now for $p > p_0$ on $L^p\otimes E$, recall that $h_p := h^{\otimes p}\otimes h^E$, and we set a new metric
    \begin{equation}
    \label{eq:log-bk2}
        H_p := h^{\otimes(p-p_0)} \otimes \tilde{h}^{\otimes p_0}\otimes h^E \, .
    \end{equation}
    Then by \eqref{eq:5.1.3-24} ($c_1(E,h^E)$ on $\Sigma$ can be properly bounded), for $p> p_0$,
    \begin{equation}
        c_1(L^p\otimes E, H_p) \geqslant p_0 \varepsilon \omega_\Sigma \, ,
        \label{eq:5.1.5-24june}
    \end{equation}
    where $\varepsilon>0$ is chosen sufficiently small.     The local weight of the metric $H_p$ on $U_0$ with respect to the frame $e_L^p\otimes e_E$ is $\Psi_p := (p-p_0)\psi + p_0 \tilde{\psi}+\psi_E$.

    Now as in the proof of \cite[Theorem 4.3]{DMM16}, we need to prove that there exist constants $C_1 > 0$, $p_0\gg 1$ such that for $p > 2 p_0$ and all $z\in U_0$, there is a section $s_{z,p} \in H^0_{(2)}(\Sigma, L^p)$, such that $s_{z,p}(z) \neq 0$ and
    \begin{equation}
    \label{eq:log-bk3}
        \int_\Sigma |s_{z,p}|_{H_p}^2 \omega_\Sigma \leqslant C_1 |s_{z,p}(z)|_{H_p}^2 \, .
    \end{equation}

The technical part is to prove the existence of $s_{z,p}$. Since \eqref{eq:5.1.5-24june} holds globally on $\Sigma$ and $(\Sigma,\omega_\Sigma)$ is complete, we can proceed as in  \cite[Proof of Theorem 5.1]{CM15} and \cite[(4.23) - (4.31)]{DMM16}. More precisely, one can construct the local holomorphic sections near $x$ as in \eqref{eq:log-bk3} by the Ohsawa--Takegoshi extension theorem \cite{OT87}, then applying the $L^2$-estimates for $\db$-operator on complete K\"{a}hler manifold (see \cite[Theorem 4.1 - (ii)]{DMM16} or \cite[Th\'{e}or\`{e}me 5.1]{MR690650}) to modify these local holomorphic sections to finally obtain global ones as wanted for \eqref{eq:log-bk3}. We may and will choose $s_{z,p}$ such that
    \begin{equation}
        \int_\Sigma |s_{z,p}|_{H_p}^2 \omega_\Sigma =1\; ,\quad  |s_{z,p}(z)|_{H_p}^2 \geqslant \frac{1}{C_1}\, .
     \label{eq:log-bk3-24}
    \end{equation}

 Since $h \leqslant \tilde{h}$ on $\Sigma$, the first property of \eqref{eq:log-bk3-24} and the definition of $H_p$ imply that
    \begin{equation}
    \int_\Sigma |s_{z,p}|_{h_p}^2 \omega_\Sigma  \leqslant 1 \, .
    \label{eq:5.1.15-24-june}
    \end{equation}
     Then the second property of \eqref{eq:log-bk3-24} implies that
    \begin{equation}
        |s_{z,p}(z)|_{h_p}^2 \geqslant \frac{1}{C_1} e^{2p_0 (\tilde{\psi}(z) - \psi(z))} \, .
    \end{equation}

Note that the quantity $e^{2p_0 (\tilde{\psi}(z) - \psi(z))}$, defined on $U_0$, actually is a global function on $\Sigma$, by the definition of $\tilde{h}$ in \eqref{eq:5.1.6-24-June},
\begin{equation}
    e^{2p_0 (\tilde{\psi}(z) - \psi(z))} = h_z^{\otimes p_0}/\tilde{h}_z^{\otimes p_0}=e^{p_0(\widetilde{F}(z)-M_0)/2 \delta}.
\end{equation}

    Recall the variational characterization of the Bergman kernel,
    \begin{equation}
        B_p(z) = \max \lb |s_p(z)|_{h_p}^2: s \in H^0_{(2)}(\Sigma, L^p\otimes E), \|s_p\|_{\cLL} = 1 \rb \, .
        \label{eq:5.1.18-24-june}
    \end{equation}
    Note that each time we work on a small local chart of a point $x\in \Sigma$, then we can use finitely many such local charts to cover the set $\Sigma\setminus \cup_j V_j$. As a consequence, we can choose uniformly the constant $C_1 \gg 0$ for all points $z\in \Sigma\setminus \cup_j V_j$, from \eqref{eq:5.1.15-24-june} - \eqref{eq:5.1.18-24-june}, we get
    \begin{equation}
        \log B_p(z) \geqslant \log |s_{z,p}(z)|_{h_p}^2 \geqslant \frac{p_0}{2 \delta}(\widetilde{F}(z)-M_0) - \log C_1 =: H(z) \, ,
        \label{eq:5.1.18-24}
    \end{equation}
    where $H \leqslant 0$. For the point $z\in \Sigma_{p,t,\gamma}\cap V_j$, we need use \eqref{eq:24-3.1.6} and \eqref{eq:bk-on-puncture} to get a lower bound for $\log B_p(z)$. So that \eqref{eq:5.1.18-24} holds uniformly for all $z\in \Sigma_{p,t,\gamma}$ for $p\gg 1$.

     Since $\widetilde{F}$ is smooth on $\overline{\Sigma}$ and $\int_\Sigma \omega_\Sigma <\infty$, then $H \in \mathcal{L}^1(\Sigma, \omega_\Sigma)$, so that we get the inequality \eqref{eq:log-bk}.
\end{proof}

 \begin{remark}
     As we saw from the above, Theorem \ref{thm:5.1.1-24} is closely related to the situations solved in \cite[Theorem 5.1]{CM15} or in \cite[Theorems 4.3 and 4.5]{DMM16}. If we regard $L$ as a holomorphic line bundle on $\overline{\Sigma}$ with singular metric $h$, the results in \cite[Theorem 5.1]{CM15} or in \cite[Theorem 4.3]{DMM16} can apply if we use a smooth K\"{a}hler metric on $\overline{\Sigma}$. However, here $\omega_\Sigma$ on $\overline{\Sigma}$ becomes singular. If we work on the noncompact model $\Sigma$ with smooth K\"{a}hler metric $\omega_\Sigma$, then \cite[Theorem 4.5]{DMM16} applies only on the open subset away from the vanishing points $\Sigma_\ast=\{z\in\Sigma\;: \; R^L_z=0 \}$ of $R^L$. Therefore, we cannot apply \cite[Theorem 5.1]{CM15} or \cite[Theorems 4.3 and 4.5]{DMM16} directly to obtain our Theorem \ref{thm:5.1.1-24}, but the basic strategy of the proof remains the same.
 \end{remark}


\subsection{On Tian's approximation theorem}\label{ss5.2Tian}
Tian's approximation theorem and its analogues are the key step to obtain the equidistribution result of random zeros for $\bb{S}_p$. Now, let us work out a version of Tian's approximation theorem in our setting. For each $p\gg 1$, consider the Kadaira map,
\begin{equation}
    \Phi_p: \Sigma 
    \longdashrightarrow[->, >=latex]{
    \hphantom{12345}
    }
    \P(H^0_{(2)}(\Sigma, L^p\otimes E)^\ast).
\end{equation}
We will use $\omega_{\FS}$ to denote the Fubini-Study metric on $\P(H^0_{(2)}(\Sigma, L^p\otimes E)^\ast)$ (see \cite[Subsection 5.1.1]{MM07}). If $U$ is a relatively compact open subset of $\Sigma$, then for sufficiently large $p$, $\Phi_p|_U$ is well-defined, and the pull-back $\Phi_p^\ast \omega_{\FS}|_U$ is a smooth form on $U$. In general,  $\Phi_p^\ast \omega_{\FS}$ defines a measure on $\Sigma$ (which might be singular), that is called the induced Fubini-Study current (or form) on $\Sigma$. It is well-known that 
\begin{equation}
    \Phi_p^\ast \omega_{\FS}= p c_1(L,h)+c_1(E,h^E)+\frac{\ui}{2\pi}\partial\overline{\partial}\log B_p(x).
    \label{eq:5.2.2-24}
\end{equation}

For any open subet $U\subset\Sigma$, recall that the norm $\|\cdot\|_{U,-2}$ for the measures or distributions on $\Sigma$ was defined in \eqref{eq:1.3.4june-24}.
\begin{definition}[Convergence speed]\label{defn:speed}
    Let $\{c_p\}_p$ be a sequence of positive numbers converging to $0$ (as $p\to +\infty$), and let $\{T_p\}_p$ and $T$ be measures on $\Sigma$ with full measures bounded by a fixed constant. We say that the sequence $\{T_p\}_p$ converges on $U$ to $T$ with speed $\mathcal{O}(c_p)$ if there exists a constant $C>0$ such that $\|T_p-T\|_{U,-2}\leqslant Cc_p$ for all sufficiently large $p$. 
\end{definition}

\begin{theorem}[Tian's approximation theorem]\label{thm:5.2.2-24}
	Let $\Sigma$ be a punctured Riemann surface, and let $L$ 
	be a holomorphic line bundle as above such that $L$ carries a singular 
	Hermitian metric $h_{L}$ satisfying conditions \ref{itm:1} and 
	\ref{itm:2}. Let $E$ be a holomorphic line bundle on 
	$\Sigma$ equipped with a smooth Hermitian metric $h^{E}$ such 
	that $(E,h^{E})$ on each chart $V_{j}$ is exactly the trivial Hermitian line bundle.
   We have the convergences of the induced Fubini-Study forms as follows. 
    \begin{enumerate}[label=(\roman*)]
        \item\label{5.2.1-1}  For any relatively compact open subset $U\subset \Sigma\;$, we have the convergence $$\frac{1}{p} \Phi_p^\ast \omega_{\FS} \longrightarrow c_1(L,h^L)$$
        in the norm $\|\cdot\|_{U,-2}$ as $p \to \infty$, with speed $\mathcal{O}(\log p/{p})$ on $U$. In particular, we have the weak convergence of measures on $\Sigma\,$,
        $$\frac{1}{p} \Phi_p^\ast \omega_{\FS} \longrightarrow c_1(L,h^L).$$
        \item\label{5.2.1-2} For any relatively compact open subset $U\subset \Sigma_2\;$, for any $\ell\in\bN\,$, there exists $C_{\ell,U}>0$ such that for $p\gg 1$,
         \begin{equation}
            \left\| \frac{1}{p} \Phi_p^\ast \omega_{\FS}-c_1(L,h^L) \right\|_{\mathscr{C}^\ell(U)}\leqslant \frac{C_{\ell, U}}{p}.
        \end{equation}
        \item\label{5.2.1-3} Fix $x\in\Sigma\;$, there exists $C_x>0$ such that for all $p\gg 1$, we have
        \begin{equation}
             \left| \frac{1}{p} (\Phi_p^\ast \omega_{\FS})(x) - c_1(L,h)(x) \right| \leqslant \frac{C_x}{\sqrt{p}}.
             \label{eq:5.2.5-24}
        \end{equation}
    \end{enumerate}
\end{theorem}
\begin{proof}
 By \eqref{eq:5.2.2-24}, we have
    $$ \frac{1}{p} \Phi_p^\ast \omega_{\FS} - c_1(L,h^L)=\frac{1}{p}c_1(E,h^E)+\frac{\ui}{2\pi p}\partial\overline{\partial}\log B_p(x). $$
    Note that any compact set in $\Sigma$ will lie in $\Sigma_{p,t,\gamma}$ for all $p\gg 1$,
    then \ref{5.2.1-1} follows directly from Theorem \ref{thm:5.1.1-24} and the definition of $\|\cdot\|_{U,-2}$.

    When the open subset $U$ is relatively compact in $\Sigma_2$, then the asymptotic expansion $B_p(x)$ on $U$ behaves the same as in \cite[Theorems 4.1.1 and 6.1.1]{MM07}, so that \ref{5.2.1-2} follows from the same arguments for \cite[Theorem 5.1.4 and Corollary 6.1.2]{MM07}. 

    Now we consider \ref{5.2.1-3}. If $x\in \Sigma_2$, then \eqref{eq:5.2.5-24} follows from \ref{5.2.1-2}. If $x\in \Sigma\setminus\Sigma_2$, then by Theorems \ref{thm:bkexpansion} and \ref{thm:4.4.1june}, we conclude that 
         \begin{equation}
             \left| \frac{1}{p} (\Phi_p^\ast \omega_{\FS})(x) - c_1(L,h)(x) \right| \leqslant \frac{C_x}{p^{1-\,\sfrac{2}{\rho_x}}},
             \label{eq:5.2.7-24}
        \end{equation}
        then by $\rho_x\geqslant 4$, we get \eqref{eq:5.2.5-24}. In this way, we complete the proof.
\end{proof}

The original Tian's approximation theorem, started with Tian \cite{Tia90} and further developed by Ruan \cite{MR1638878}, Catlin \cite{MR1699887}, and Zelditch \cite{Zel98}, is for the case of positive line bundles on compact K\"{a}hler manifolds. Then Ma and Marinescu \cite{MM07} extended it for the uniformly positive line bundles on complete Hermitian manifolds. For big or semipositive line bundles equipped with possibly singular Hermitian metrics, the $(1,1)$-current versions of Tian's approximation theorem have been widely studied, such as by Coman and Marinescu \cite{CM13,CM15}, Dinh, Ma, and Marinescu \cite{DMM16}.

\subsection{Equidistribution of random zeros and convergence speed}\label{ss:equidis}
In this subsection, we give a proof of Theorem \ref{thm:1.3.2-24}. We only consider $p\gg 1$. The standard Gaussian holomorphic section $\bb{S}_p$ is defined in Definition \ref{defn:Gaussian}. By \cite[Subsection 5.3]{MM07} (see also \cite[Theorem 1.1]{DrLM:2023aa}), we know that $\E[[\Div(\bb{S}_p)]]$ exists as a positive distribution (hence a measure) on $\Sigma$, and we have the identity
\begin{equation}
        \E[[\Div(\bb{S}_p)]]=\Phi_p^\ast\omega_{\FS}= p c_1(L,h)+c_1(E,h^E)+\frac{\ui}{2\pi}\partial\overline{\partial}\log B_p(x).
    \label{eq:5.3.1-24}
\end{equation}

Let $V$ be a Hermitian vector space of complex dimension $d+1$. On projective space $\P(V^\ast)$, let $\sigma_\FS$ denote the normalized Fubnini-Study volume form on $\P(V^\ast)$ so that it defines a uniform probability measure on $\P(V^\ast)$, that is,
\begin{equation}
    \sigma_\FS:=\omega_{\FS}^d.
\end{equation}
Meanwhile, for a non-zero $\xi\in V^\ast$, let $H_\xi=\ker \xi$ be the hyperplane in $V$ so that it defines a positive $(1,1)$-current $[H_\xi]$ on $\P(V)$. Similar to \eqref{eq:1.3.4june-24}, we can define the norm $\|\cdot\|_{U,-2}$ for $(1,1)$-currents.
\begin{theorem}[{\cite[Theorem 4]{DMS12}}]\label{thm:5.3.1-24}
Let $(X,\omega)$ be a Hermitian complex manifold of dimension $n$ and let $U$ be a relatively compact open subset of $X$. Let $V$ be a Hermitian vector space of complex dimension $d+1$. There exists a constant $C>0$ independent of $d$ such that for every $\gamma>0$ and every holomorphic map $\Phi: X\longrightarrow \P(V)$ of generic rank $n$, we can find a subset $E\subset \P(V^\ast)$ satisfying the following properties:
\begin{enumerate}
    \item $\sigma_\FS(E)\leqslant C d^2 e^{-\gamma/C}$.
    \item If $[\xi]$ is outside $E$, the current $\Phi^\ast([H_\xi])$ is well-defined and we have
    \begin{equation}
        \left\| \Phi^\ast([H_\xi])-\Phi^\ast\omega_{\FS} \right\|_{U,-2} \leqslant \gamma.
    \end{equation}
\end{enumerate}
\end{theorem}

Now we can give the proof of Theorem \ref{thm:1.3.2-24}.
\begin{proof}[Proof of Theorem \ref{thm:1.3.2-24}]
    At first, Theorem \ref{thm:1.3.2-24} - \ref{1.3.2-1} follows from Theorem \ref{thm:5.2.2-24} - \ref{5.2.1-1} and \eqref{eq:5.3.1-24}.

    Let us focus on the proof of Theorem \ref{thm:1.3.2-24} - \ref{1.3.2-2}. Consider the probability space $(\P(H^{0}_{(2)}(\Sigma, L^p\otimes E)),\sigma_{\FS})$, to each $[s_p]\in \P(H^{0}_{(2)}(\Sigma, L^p\otimes E))$, we associated with the measure defined by its zero divisor $\Div(s_p)$; this way, we constructed a random variable $\bb{\mu}_p$ valued in the measures on $\Sigma$. Then $\bb{\mu}_p$ has the same probability distribution as $[\Div(\bb{S}_p)]$. So, now we proceed with the proof for the sequence $\{\bb{\mu}_p\}_p$ using the arguments as in \cite[Proof of Theorem 2]{DMS12}.

Let $U'$ be a relatively compact open subset in $\Sigma$ such that $\overline{U}\subset U'$. For each $p\gg 1$, take $V=H^{0}_{(2)}(\Sigma, L^p\otimes E)^\ast$ in Theorem \ref{thm:5.3.1-24} and map $\Phi$ is given by the Kodaira map $\Phi_p$, when we restrict the map to $U'$, so that \ref{thm:5.3.1-24} applies.  Note that for $[s_p]\in \P(H^{0}_{(2)}(\Sigma, L^p\otimes E))$, the positive $(1,1)$-current (hence measure) $\Phi_p^\ast([H_{s_p}])$ on $U'$ is exactly the measure $[\Div(s_p)]|_{U'}$. 

Since the constant $C$ in Theorem \ref{thm:5.3.1-24} is independent of the choices of $d$ or $\gamma$. We take the sequence $\gamma_p=4C \log{p}$. We conclude that for all $p\gg 1$,
\begin{equation}
    \sigma_{\FS}\left(\left\|\frac{1}{p}\bb{\mu}_p-\frac{1}{p}\Phi_p^\ast\omega_{\FS}\right\|_{U,-2}>\frac{4C \log{p}}{p}\right)\leqslant \frac{C'}{p^2},
    \label{eq:5.3.4-24}
\end{equation}
with certain constant $C'>0$. Then by the equivalence between $[\Div(\bb{S}_p)]$ and $\bb{\mu}_p$ and Theorem \ref{thm:5.2.2-24} - \ref{5.2.1-1}, we get for $p\gg 1$,
\begin{equation}
    \P_p\left(\left\|\frac{1}{p}[\Div(\bb{S}_p)]-c_1(L,h)\right\|_{U,-2}>\frac{\widetilde{C}\log{p}}{p}\right)\leqslant \frac{C'}{p^2},
    \label{eq:5.3.5-24}
\end{equation}
Since $\sum_p \frac{C'}{p^2}<\infty$, we conclude exactly \ref{thm:1.3.2-24} - \ref{1.3.2-2}.   
\end{proof}

\begin{remark}
    The probability inequality \eqref{eq:5.3.5-24} has a similar nature as our large deviation estimates \eqref{eq:1.4.4-24} (whose proof is given in the next subsection). In fact, from \eqref{eq:1.4.4-24}, one can also deduce the equidistribution result for $\bb{S}_p$ on $U$ but without the convergence speed $\mathcal{O}(\log{p} / p)$. 
    If we take the sequence $\lambda_p\cong \delta p$ in \eqref{eq:5.3.4-24} and \eqref{eq:5.3.5-24}, then we get
    \begin{equation}
    \P_p\left(\left\|\frac{1}{p}[\Div(\bb{S}_p)]-c_1(L,h)\right\|_{U,-2}>\delta\right)\leqslant C'p^2 e^{-c\delta p},
    \label{eq:5.3.6-24}
\end{equation}
    For a given $\delta$, the above inequality is less sharp than \eqref{eq:1.4.4-24}.
\end{remark}

\subsection{Large deviation estimates and hole probability}\label{ss:LDE}
In this subsection, we will prove Theorem \ref{thm:5.3.1} and Proposition \ref{prop:1.2.3}, which consists of the arguments in \cite[Subsection 3.3 - 3.6]{Drewitz_2023} with small modifications. We always assume the geometric conditions in Subsection \ref{ss:1.1amiens}.

For an open subset $U\subset \Sigma$, $s_p\in H^0_{(2)}(\Sigma, L^p\otimes E)$, set
\begin{equation}
    \mathcal{M}^U_p(s_p):=\sup_{x\in U} |s_p(x)|_{h_p}.
\end{equation}
The following proposition is an extension of \cite[Theorem 1.4 and Proposition 1.9]{Drewitz_2023} for semipositive line bundles, as an application of Proposition \ref{cor:1.2.3} and Theorem \ref{thm:1.3.1}.

\begin{proposition}\label{prop:1.3.3}
	Let $U$ be a relatively compact open subset in $\Sigma$. For any $\delta>0$, there exists $C_{U,\delta}>0$ such that for all $ p\gg 1$,
	\begin{equation}
		\P_{p} \left( \left\{ s_{p}\;:\; \left|\log{\mathcal{M}^{U}_{p}(s_{p})} \right| \geqslant 
		\delta p \right\} \right)\leqslant e^{-C_{U,\delta}p^{2}}\,.
  		\label{eq:5.4.2-24}
	\end{equation}
 As a consequence, there exists $C'_{U,\delta}>0$ such that for all $ p\gg 1$,
	\begin{equation}
\P_{p}\left(\left\{s_{p}\;:\; 
		\int_{U}\big|\log{|s_{p}|_{_{h_{p}}}}\big| \,\omega_{\Sigma}\geqslant 
		\delta p\right\} \right)\leqslant e^{-C'_{U,\delta}p^{2}}.
		\label{eq:1.3.5}
	\end{equation}
\end{proposition}
\begin{proof}
    At first, the proof of \eqref{eq:1.3.5} follows from the same arguments as in \cite[Subsection 3.4]{Drewitz_2023} and \eqref{eq:5.4.2-24}. So we now focus on proving \eqref{eq:5.4.2-24}.

    As explained in \cite[Subsection 3.3]{Drewitz_2023}, the proof of \eqref{eq:5.4.2-24} consists of two parts:
    \begin{enumerate}
        \item \label{s-1-1} Using the uniform upper bound on $B_p(x)$ from Proposition \ref{cor:1.2.3} and proceeding as in \cite[Subsection 3.1]{Drewitz_2023} (in particular, \cite[Corollary 3.6]{Drewitz_2023}), then we get
        $$\P_{p} \left( \left\{ s_{p}\;:\; \mathcal{M}^{U}_{p}(s_{p})\geqslant 
		e^{\delta p} \right\} \right) \leqslant e^{-C_{U,\delta}p^{2}}\,.$$
        \item \label{s-1-2} Since $\Sigma_2$ is an open dense subset of $\Sigma$, then for any (non-empty) open subset $U$, we can always find a small open ball in $\bB \subset U\cap \Sigma_2$ such that the expansion in Theorem \ref{thm:1.3.1} for $N_p(x,y)$ holds for $x,y\in \bB$. Then we consider a sequence of lattices $\Gamma_p$ in $\bB$ with mesh $\sim \frac{1}{\sqrt{p}}$ and proceed as in \cite[Subsection 3.3]{Drewitz_2023}, we conclude
        $$\P_{p} \left( \left\{ s_{p}\;:\; \mathcal{M}^{U}_{p}(s_{p})\leqslant 
		e^{-\delta p} \right\} \right)\leqslant e^{-C_{U,\delta}p^{2}}\,.$$
    \end{enumerate}
    In this way, we get \eqref{eq:5.4.2-24}. The proposition is proved.
\end{proof}
\begin{remark}
    Since Proposition \ref{cor:1.2.3} gives the global uniform upper bound for $B_p(x)$, if $U$ is an open subset but not relatively compact in $\Sigma$,  \eqref{eq:5.4.2-24} still holds.
\end{remark}

Now we are ready to prove Theorem \ref{thm:5.3.1}.
\begin{proof}[Proof of Theorem \ref{thm:5.3.1}]

Let us start with Theorem \ref{thm:5.3.1} - \ref{1.4.2-1}.
Fix $\varphi\in \mathscr{C}^\infty_{\mathrm{c}}(\Sigma)$ with $\supp \varphi\subset U$, by Poincar\'{e}-Lelong formula \eqref{eq:1.3.3-24june}, we have
\begin{equation}
\begin{split}
\left\langle\frac{1}{p}[\Div(\bb{S}_{p})],\varphi\right\rangle
-\int_{\Sigma}\varphi c_{1}(L,h)=\frac{\sqrt{-1}}{p\pi}\int_{\Sigma}
\log{|\bb{S}_{p}|_{h_{p}}}\,\partial\overline{\partial}\varphi +\frac{1}{p}\langle c_1(E,h^E),\varphi\rangle.
		\end{split}
		\label{eq:3.4.14Jan}
	\end{equation}
Since $\varphi$ has a compact support in $U$, so has $\partial\overline{\partial}\varphi$. 
Then
	\begin{equation}
	\begin{split}
	\left|\frac{\sqrt{-1}}{p\pi}\int_{\Sigma}\log{|\bb{S}_{p}|_{h_{p}}}\,
		\partial\overline{\partial}\varphi\right|\leqslant 
		\frac{\|\varphi\|_{\cC^2(U)}}{p\pi}\int_{U}
		\left|\log{|\bb{S}_{p}(x)|_{h_{p}}}\right|\,\omega_{\Sigma}(x).
		\end{split}
		\label{eq:3.4.17Jan}
	\end{equation}

We fix a sufficiently small $\varepsilon>0$ such that
$$\delta-2\varepsilon > 0.$$ 
Since the term $\frac{1}{p} c_1(E,h^E)$ converges to $0$ as 
$p \to \infty$, there exists an integer $p_0\in\N$
(depending on $(E,h^E)$) such that for all $p\geqslant p_0$,
\begin{equation}
\left|\frac{1}{p}\langle c_1(E,h^E),\varphi\rangle\right|
\leqslant \frac{\varepsilon \|\varphi\|_{\cC^2(U)}}{\pi}\,\cdot
\end{equation} 
Applying \eqref{eq:1.3.5} to the right-hand side of 
\eqref{eq:3.4.17Jan} with $\delta-2\varepsilon$, we get, for $p\gg 1$,
\begin{equation}
		\P\left(\frac{1}{p}
		\int_{U}\Big|\log{\big|\bb{S}_{p}(x)\big|_{h_{p}}}\Big| 
		\,\omega_{\Sigma}(x)>
		\delta-2\varepsilon \right)\leqslant e^{-Cp^{2}}.
		\label{eq:4.67Jan}
	\end{equation}
For $p\geqslant p_0$, except the event from \eqref{eq:4.67Jan} of probability 
$\leqslant e^{-Cp^{2}}$, we have that, for all $\varphi\in \cC^\infty_{\mathrm{c}}(U)$,
\begin{equation}
\begin{split}
&\left| \left\langle \frac{1}{p}[\Div(\bb{S}_{p})] 
-  c_{1}(L,h),\varphi \right\rangle \right|\\
&\leqslant   
\frac{\|\varphi\|_{\cC^2(U)}}{p\pi}\int_{U}
\left|\log{|\bb{S}_{p}(x)|_{h_{p}}}\right|\,\omega_{\Sigma}(x)+ 
\left|\frac{1}{p} \left\langle c_1(E,h^E),\varphi \right\rangle \right|\\
&\leqslant \frac{1}{\pi} \left( \|\varphi\|_{\cC^2(U)}(\delta-2\varepsilon)+
\varepsilon \|\varphi\|_{\cC^2(U)} \right) \leqslant 
\|\varphi\|_{\cC^2(U)}\frac{\delta-\varepsilon}{\pi},
\end{split}
\end{equation}
Equivalently, except the event in \eqref{eq:4.67Jan} of probability 
$\leqslant e^{-Cp^{2}}$, we have 
\begin{equation}
\left\|\frac{1}{p}[\Div(\bb{S}_{p})] - c_{1}(L,h_L)\right\|_{U,-2}\leqslant \frac{\delta-\varepsilon}{\pi}.
\end{equation}
Hence \eqref{eq:1.4.4-24} follows. 

Now we consider Theorem \ref{thm:5.3.1} - \ref{1.4.2-2}. If $U$ is still relatively compact in $\Sigma$, then \eqref{eq:1.4.5DLM} follows from \eqref{eq:1.4.4-24} and the arguments as in \cite[Subsection 3.6]{Drewitz_2023}. However, here we allow $U$ to contain the punctures. Since the line bundle $L$ is positive on $\overline{\Sigma}$, the arguments \cite[Subsection 3.5]{Drewitz_2023} (to control the vanishing order at punctured points) together with Proposition \ref{prop:1.3.3} show that \cite[Theorem 1.10]{Drewitz_2023} still holds in our case. As a consequence, the arguments as in \cite[Subsection 3.6]{Drewitz_2023} still apply and we get \eqref{eq:1.4.5DLM} in full generality. Finally, using Borel-Cantelli type arguments to \eqref{eq:1.4.5DLM}, we get \eqref{eq:1.4.8-24}. 
\end{proof}

\begin{proof}[Proof of Proposition \ref{prop:1.2.3}]
    The upper bound \eqref{eq:1.4.9-24june} follows directly from \eqref{eq:1.4.5DLM} with $\delta=\mathrm{Area}^L(U)$. The lower bound \eqref{eq:1.2.4} follows from the same arguments as in \cite[Subsection 4.2.4]{SZZ08} (see also \cite[Subsection 3.7]{Drewitz_2023}).
\end{proof}

\subsection{Smooth statistics: leading term of number variances}\label{ss:nv}

Following Shiffman and Zelditch \cite[\S 3]{SZ08}, we now introduce the variance current of $[\mathrm{Div}(\bb{S}_{p})]$. Let $\pi_1, \pi_2: \Sigma \times \Sigma \longrightarrow \Sigma$ denote the projections to the first and second factors. Then if $S$ and $T$ are two distributions on $\Sigma$, then we define a distribution on $\Sigma \times \Sigma$ as follows
\begin{equation}
S\boxtimes T:= \pi_1^* S\wedge \pi_2^* T.
\end{equation} 
In particular, $[\mathrm{Div}(\bb{S}_{p})]\boxtimes [\mathrm{Div}(\bb{S}_{p})]$ defines a random distribution on $\Sigma \times \Sigma$. In the same time, we introduce the following notation: for a current $T$ on $\Sigma \times \Sigma$, we write
\begin{equation}
\partial T=\partial_1 T+\partial_2 T,
\end{equation}
where $\partial_1$, $\partial_2$ denote the corresponding $\partial$-operators on the first and second factors of $\Sigma \times \Sigma$. Similarly, we also write $\overline{\partial} T=\overline{\partial}_1 T+\overline{\partial}_2 T$.

\begin{definition}
The variance current of $[\mathrm{Div}(\bb{S}_{p})]$, denoted as $\mathbf{Var}[\bb{S}_{p}]$, is a distribution on $\Sigma \times \Sigma$ defined by
\begin{equation}
\mathbf{Var}[\bb{S}_{p}]:=\E\big[[\mathrm{Div}(\bb{S}_{p})]\boxtimes [\mathrm{Div}(\bb{S}_{p})]\big]-\E\left[[\mathrm{Div}(\bb{S}_{p})]\right]\boxtimes \E\left[[\mathrm{Div}(\bb{S}_{p})]\right]
\end{equation}
\end{definition}

Now we consider only the real test functions. For $\varphi\in \cC^\infty_{\mathrm{c}}(\Sigma,\R)$, we have
\begin{equation}
\mathrm{Var} \left[ \left\langle [\mathrm{Div}(\bb{S}_{p})],\varphi \right\rangle \right]= \left\langle \mathbf{Var}[\bb{S}_{p}], \varphi\boxtimes \varphi \right\rangle. 
\end{equation}

For $t\in [0,1]$, we set the function
\begin{equation}
\widetilde{G}(t):=-\frac{1}{4\pi^2} \int_0^{t^2} \frac{\log(1-s)}{s} \, \d s = \frac{1}{4\pi^2} \sum_{j=1}^\infty \frac{t^{2j}}{j^2}.
\label{eq:6.6geb}
\end{equation}
This is an analytic function with radius of convergence $1$. Moreover, for $t\sim 0$, we have $\widetilde{G}(t)=\mathcal{O}(t^2)$.

Recall that $N_p(z,w)$ is the normalized Bergman kernel defined in \eqref{eq:1.3.1}.
\begin{definition}[{cf. \cite[Theorem 3.1]{SZ08}}]
For $(z,w)\in \Sigma \times \Sigma$, define
\begin{equation}
Q_{p}(z,w):=\widetilde{G}(N_{p}(z,w))=-\frac{1}{4\pi^2} \int_0^{N_{p}(z,w)^2} \frac{\log(1-s)}{s} \, \d s.
\end{equation}
\end{definition}

Following the calculations in \cite[\S 3.1]{SZ08} and using Theorem \ref{thm:1.3.1} and Lemma \ref{lm:2.14feb24}, 
we have the following results for $Q_{p}(z,w)$ on the open set $\Sigma_2\times \Sigma_2$.
\begin{proposition}[{cf. \ \cite[Lemmas 3.4, 3.5 and 3.7]{SZ08}}]
\label{prop:6.3feb}
Let $U$ be a relatively compact open subset of $X$ such that $\overline{U}\subset \Sigma_2$.
\begin{enumerate}[label=(\roman*)]
\item\label{5.5.4-1} Then there exists an integer $p_0\in \N$ such that for all $p\geqslant p_0$, $B_{p}(z)$ never vanishes on $\overline{U}$. Moreover, for all $p\geqslant p_0$, the function $Q_{p}(z,w)$ is smooth in the region $U\times U\setminus \Delta_U$ ($\Delta_U$ denotes the diagonal) and it is $\cC^1$ on $U\times U$.
\item\label{5.5.4-2} Fix $b\gg 0$ and $\varepsilon>0$, then for all sufficiently large $p$ and for $x\in U$, $Z\in T_x \Sigma$ with $|Z|\leqslant b\sqrt{\log{p}}$, we have
\begin{equation}
Q_{p}(x,\exp_x(Z/\sqrt{p}))=\widetilde{G} \left( \exp \left\{-\bb{c}(x)|Z|^2/4 \right\} \right)+\mathcal{O}(p^{-1/2+\varepsilon}),
\end{equation}
where $\bb{c}(x)$ is defined in \eqref{eq:1.2.6june}.
\item\label{5.5.4-3} For given $k, \ell\in \N$, there exist a sufficiently large $b>0$ such that there exist a constant $C>0$ such that for all $z,w\in U$, $\mathrm{dist}(z,w)\geqslant b\sqrt{\log{p}/p}$, we have
\begin{equation}
|\nabla^\ell_{z,w}Q_{p}(z,w)|\leqslant C p^{-k}.
\label{eq:6feb24}
\end{equation}
\end{enumerate}
\end{proposition}

The same proof of \cite[Theorem 3.1]{SZ08} (see also \cite[\S 3.1]{MR2742043}) together with Proposition \ref{prop:6.3feb} - \ref{5.5.4-1} shows the following result.
\begin{theorem}[cf. {\cite[Theorem 3.1]{SZ08}}]\label{thm:6.4number}
We assume the same conditions on $\Sigma$, $L$ and $E$ as in Theorem \ref{thm:1.2.1}. Let $U$ be a relatively compact open subset of $\Sigma$. Then for sufficiently large $p$, we have the identity of distribution on $U\times U$,
\begin{equation}\label{eq:6.11feb24}
\mathbf{Var}[\bb{S}_{p}]|_{U\times U}=-\partial_1\overline{\partial}_1\partial_2\overline{\partial}_2 Q_{p}(z,w)|_{U\times U}=(\sqrt{-1}\partial\overline{\partial})_z  (\sqrt{-1}\partial\overline{\partial})_w Q_{p}(z,w)|_{U\times U}.
\end{equation}
\end{theorem}

Recall that the operator $\fL(\varphi)$ and the test function space $\mathcal{T}^3(L,h)$ are defined in Definition \ref{defn:1.5.1-24}. Now we give the proof of Theorem \ref{thm:1.5.2-24ss}.
\begin{proof}[Proof of Theorem \ref{thm:1.5.2-24ss}]
Fix $\varphi\in  \mathcal{T}^3(L,h) $ with $\partial\overline{\partial}\varphi\not\equiv 0$, and let $U$ be a relatively compact open subset of $\Sigma$ such that $\supp \varphi \subset U$. Note that $U$ may contain the vanishing points of $R^L$. 

Since $\fL(\varphi)$ vanishes identically near $\Sigma_\ast$, then there exists a sufficiently small $\delta>0$, such that 
\begin{equation}
    \fL(\varphi)|_{V(R^L,\delta)}\equiv 0,
    \label{eq:5.5.11-24}
\end{equation}
where $V(R^L,\delta):=\{z\in\Sigma\;:\; \mathrm{dist}(z,\Sigma_\ast)\leqslant \delta\}$ is the closed tubular neighbourhood of $\Sigma_\ast$ in $\Sigma$. 
We write
\begin{equation}
    U=U_1(\delta) \cup U_2(\delta),
    \label{eq:5.5.12-24}
\end{equation}
where $U_1(\delta):=U\cap V(R^L,\delta)$, and $U_2(\delta) = U\cap (\Sigma\setminus V(R^L,\delta))$ is a relatively compact open subset of $\Sigma_2$.

Then by \eqref{eq:6.11feb24}, \eqref{eq:5.5.11-24} and \eqref{eq:5.5.12-24}, we have
 \begin{equation}\label{eq:6.15feb}
 \begin{split}
 \mathrm{Var} \big[ \left\langle \left[ \mathrm{Div}(\bb{S}_{p}) \right], \varphi \right\rangle \big] &= - \int_{U\times U} (\partial\overline{\partial}\varphi(z))\wedge
 (\partial\overline{\partial}\varphi(w))\widetilde{G}(N_{p}(z,w))\\
 &=- \int_{U_2(\delta)\times U_2(\delta)} (\partial\overline{\partial}\varphi(z))\wedge
 (\partial\overline{\partial}\varphi(w))\widetilde{G}(N_{p}(z,w))
  \end{split}
 \end{equation}

Therefore, the calculation reduces for the subset $U_2(\delta)$. By construction of $U_2(\delta)$, Proposition \ref{prop:6.3feb} - \ref{5.5.4-2} and \ref{5.5.4-3} hold uniformly for $z,w \in U_2(\delta)$. Then we can proceed as in \cite[\S 3.1]{MR2742043} (see also \cite[Proof of Theorem 6.4]{Drewitz:2024aa}), we conclude \eqref{eq:6.12feb24}.
\end{proof}

\begin{remark}
Note that following the work of Shiffman \cite{MR4293941}, one can obtain the full expansion of the variance $\mathrm{Var}\big[\langle [\mathrm{Div}(\bb{S}_{p})],\varphi\rangle\big]$ and calculate the subleading term. 
\end{remark}

For better understanding on the vanishing points of $R^L$ and the space $\mathcal{T}^3(L,h)\,$, let us introduce an intuitive but nontrivial lemma; we refer to the short article \cite{MR2351134} for a proof.
\begin{lemma}\label{lm:5.5.1-24june}
    Let $\alpha$ be a smooth $(1,1)$-form on $\Sigma$ such that it only vanishes on a compact subset of $\Sigma$ and with finite vanishing orders. Set
    $V(\alpha):= \left\{z\in\Sigma\;:\; \alpha(z)=0 \right\}$, and for $\delta>0$, set
    $$V(\alpha,\delta) = \left\{z\in\Sigma\;:\; \mathrm{dist}(z,V(\alpha))\leqslant \delta \right\}\subset \Sigma.$$
    Then there exist constants $\delta_0\in\,]0,1[\, , C_0>0$ such that for any $0<\delta<\delta_0$, we have
    \begin{equation}
        \int_{V(\alpha,\delta)}\omega_\Sigma\leqslant C_0 \delta.
    \end{equation}
\end{lemma}

 As a consequence of the above lemma, there are always test functions $\varphi$ in $\mathcal{T}^3(L,h)$ such that the vanishing points of $\fL(\varphi)$ near $\Sigma_\ast$ have arbitrarily small size. For example, consider the set $U_1(\delta)$ given in \eqref{eq:5.5.12-24}, by Lemma \ref{lm:5.5.1-24june}, there exists a constant $C_U>0$ independent of $\delta$ such that
\begin{equation}
    \int_ {U_1(\delta)}\omega_\Sigma \leqslant C_U \delta.
\end{equation}

If $\psi$ is an arbitrary real test function on $\Sigma$ with support in $U$, then we can modify the values of $\psi$ on $U_1(\delta)$ to construct a real test function $\widetilde{\psi}_\delta$ such that: it coincides with $\psi$ outside $U_1(\delta)$ and is locally constant on $U_1(\delta/2)$; it satisfies
$$ \left\| \psi-\widetilde{\psi}_\delta \right\|_{\cC^0(\Sigma)}\leqslant \|\psi\|_{\cC^0(\Sigma)}.$$
This way, we get $\widetilde{\psi}_\delta\in  \mathcal{T}^3(L,h)$, and
\begin{equation}
  \P_{\infty}\left( \limsup_{p \to +\infty}\frac{1}{p} \left|Y_p(\psi)-Y_p(\widetilde{\psi}_\delta) \right|\leqslant C_U \delta \|\psi\|_{\cC^0(\Sigma)} \right)=1.
\end{equation}
Since $\delta$ is arbitrarily small, we can view $\frac{1}{p}Y_p(\widetilde{\psi}_\delta)$ as a $\delta$-approximation of $\frac{1}{p}Y_p(\psi)$.

\subsection{Smooth statistics: central limit theorem for random zeros}\label{ss:CLT-24}
Let us recall the main result of \cite[\S 2.1]{STr}. Let $(T,\mu)$ be a measure space with a finite positive measure $\mu$ (with $\mu(T)>0$). We also fix a sequence of measurable functions $A_k: T \longrightarrow \C$, $k\in\N$ such that on $T$,
\begin{equation}
\sum_k |A_k(t)|^2\equiv 1.
\label{eq:6.25feb24}
\end{equation}
We consider a complex-valued Gaussian process on $T$ defined as
\begin{equation}
W(t):=\sum_k \eta_k A_k(t),
\label{eq:6.26feb24}
\end{equation}
where $\{\eta_k\}$ is a sequence of i.i.d.\ standard complex Gaussian variables. Then for each $t\in T$, $W(t)\sim \mathcal{N}_{\C}(0,1)$. The covariance function for $W$ is $\rho_W: T\times T\longrightarrow \C$ given by
\begin{equation}
\rho_W(s,t):=\E \left[ W(s)\overline{W(t)} \right] = \sum_k A_k(s)\overline{A_k(t)}.
\end{equation}

Let $\{W_p\}_{p\in \N}$ be a sequence of independent Gaussian processes on $T$ described as above, and let $\rho_p(s,t)$ ($p\in\N$) denote the corresponding covariance functions. We also fix a non-trivial real function $F\in\mL^2(\R_+, e^{-r^2/2} r \, \d r)$, and a bounded measurable function $\psi: T\rightarrow \R$, set
\begin{equation}
Z_p:=\int_T F \left( \left| W_p(t) \right| \right) \psi(t)\d \mu(t).
\label{eq:6.28feb24}
\end{equation}

Sodin and Tsirelson proved the following result.

\begin{theorem}[{\cite[Theorem 2.2]{STr}}]\label{thm:STr2004}
With the above construction suppose that
\begin{equation}
\tag{i}
\liminf_{p\rightarrow +\infty} \frac{\int_{T}\int_{T}
\left| \rho_p(s,t) \right|^{2\alpha}\psi(s)\psi(t) \, \d \mu(s) \, \d \mu(t)}{\sup_{s\in T}\int_T \left|\rho_p(s,t) \right| \,\d\mu(t)}>0,
\end{equation}
for $\alpha=1$ if $f$ is monotonically increasing, or for all 
$\alpha\in \N$ otherwise;
\begin{equation}
\tag{ii}
\lim_{p\rightarrow+\infty} \sup_{s\in T}
	\int_T \left| \rho_p(s,t) \right| \, \d \mu(t)=0.
\end{equation}	
Then the distributions of the random variables
\begin{equation}
\frac{Z_p-\E[Z_p]}{\sqrt{\mathrm{Var}[Z_p]}}
\end{equation}
converge weakly to  the (real) standard Gaussian distribution $\mathcal{N}_\R(0,1)$ as $p \to +\infty$.
\end{theorem}

Now we are ready to present the proof of Theorem \ref{thm:3.5.1ss}.
\begin{proof}[Proof of Theorem \ref{thm:3.5.1ss}]

Let us use the same notation as in the proof of Theorem \ref{thm:1.5.2-24ss}. Fix $\varphi\in  \mathcal{T}^3(L,h) $ with $\partial\overline{\partial}\varphi\not\equiv 0$, and fix a sufficiently small $\delta>0$ as desired.

By \eqref{eq:1.3.3-24june}, \eqref{eq:1.5.1-24} and \eqref{eq:5.5.11-24} - \eqref{eq:5.5.12-24}, we have
\begin{equation}
    Y_p(\varphi)=\int_{U_2(\delta)} \frac{1}{\pi}\log \left|\bb{S}_p(x) \right|_{h_p} \left( \sqrt{-1}\partial\overline{\partial}\varphi \right)(x) +\left\langle pc_1(L,h)+c_1(E,h^E),\varphi \right\rangle.
    \label{eq:5.6.6-24june}
\end{equation}

Let $\mathfrak{f}:\overline{U_2(\delta)} \longrightarrow L$, $\mathfrak{e}:\overline{U_2(\delta)} \longrightarrow E$ be the continuous sections 
such that $|\mathfrak{f}(z)|_{h}\equiv 1$, $|\mathfrak{e} (z)|_{h^E}\equiv 1$ on $\overline{U_2(\delta)}$. 
For each $p$, fix an orthonormal basis 
$\{S^p_j\}_{j=1}^{d_p}$ of $H^0_{(2)}(\Sigma, L^p\otimes E)$.
Then on $\overline{U_2(\delta)}$, we write
\begin{equation}
S^p_j(z)=a^p_j(z) \, \mathfrak{f}^{\otimes p}(z)\otimes \mathfrak{e}(z).
\end{equation}
Then we can set 
$A^p_j(z)={ a^p_j(z)}/{\sqrt{B_{p}(z)}}$,
which forms a sequence of measurable functions on 
$U_2(\delta)$ satisfying \eqref{eq:6.25feb24}. 
Then we have the identity on $U_2(\delta)$
\begin{equation}
\frac{\bb{S}_{p}(z)}{\sqrt{B_{p}(z)}}
=W_p(z) \, \mathfrak{f}^{\otimes p}(z)\otimes \mathfrak{e}(z),
\label{eq:6.33feb24}
\end{equation}
where $W_p$ is the Gaussian process on $U_2(\delta)$ constructed
as in \eqref{eq:6.26feb24}. The covariance function 
$\rho_p(z,w)$ for $W_p$ is given by
\begin{equation}
\left| \rho_p(z,w) \right| = N_{p}(z,w).
\end{equation}

We take $F(r)=\log{r}$, $(T,\mu)=(U_2(\delta),c_1(L,h)|_{U_2(\delta)})$, 
$\psi(z)=\frac{1}{\pi} \fL(\varphi)(z)$ which satisfies the 
conditions in Theorem \ref{thm:STr2004}. Then let $Z_p(\varphi)$ be the random variable defined
as in \eqref{eq:6.28feb24} for $W_p$ on $U_{2}(\delta)$.

Then \eqref{eq:5.6.6-24june} and \eqref{eq:6.33feb24} imply that
\begin{equation}
Y_{p}(\varphi)=Z_p(\varphi) + C_p,
\end{equation}
where $C_p$ is a deterministic constant. Thus the asymptotic 
normality of $Y_{p}(\varphi)$ is equivalent to that of 
$Z_p(\varphi)$. 

Therefore, the last step is to check the conditions (i) and (ii) 
in Theorem \ref{thm:STr2004} for $N_{p}(z,w)$ with $z,w\in U_2(\delta)$ and for $(T,\mu)=(U_2(\delta),c_1(L,h)|_{U_2(\delta)})$. 
Since $U_2(\delta)$ is a relatively compact open subset of $\Sigma_2$, Theorem \ref{thm:1.3.1} applies and we proceed
as in the last part of \cite[\S 4 Proof of Theorem 1.2]{MR2742043} to complete the proof.
\end{proof}

\appendix
\numberwithin{equation}{section}
\section{Jet-bundles and induced norms on them}\label{ssec:jet}

In this appendix, we introduce the necessary notation and notions for the jet bundles on $\Sigma$. 
Let $(F,h^{F})$ be a real (or complex) vector bundle on $\Sigma$ with 
a Euclidean (or Hermitian) inner product $h^{F}$.

For $x\in\Sigma$, let $\cG_{x}(F)$ denote the germs of local sections of 
$F$ at $x$. For $\ell\in\N$, $s\in \cG_{x}(F)$, the $\ell$-th jet of $s$ at 
$x$, denoted by $j^{\ell}_{x}s$, is the equivalence class of $s$ in 
$\cG_{x}(F)$ under the equivalence relation: two 
germs are equivalent if on some open coordinate chart containing $x$ where the bundle $F$ is 
trivialized, they have the same Taylor expansions at $x$ up to order 
$\ell$. Let $J^{\ell}(F)_{x}$ denote the vector space of all $\ell$-th jets 
$j^{\ell}_{x}s$, $s\in \cG_{x}(F)$. Then $J^{\ell}(F)_{x}$ is finite dimensional, and actually the fibration $\coprod_{x\in\Sigma} 
J^{\ell}(F)_{x}\to \Sigma$ defines in a natural way a smooth 
vector bundle on $\Sigma$, which is denoted by $J^{\ell}(F)$ and called 
the $\ell$-th jet bundle of $F$ on $\Sigma$. Note that $J^{0}(F)$ is 
just $F$ itself.

For an integer $\ell>0$, let $\pi_{\ell-1}^\ell : J^\ell(F) \longrightarrow 
J^{\ell-1}(F)$ 
denote the obvious projection of vector bundles. Observe that there 
exists a short exact sequence of vector bundles over $\Sigma$ (cf. \cite[pp.121]{KMS93})
\begin{align}
	\begin{xy}
		\xymatrix{
		0 \to S^{\ell}T^\ast\Sigma\otimes F\ar[r]^(0.6){\mathrm{incl}}& J^\ell(F) 
		\ar[r]^(0.4){\pi_{\ell-1}^\ell} & J^{\ell-1}(F) \to 0 \, ,
		}
	\end{xy}
\end{align}
where $S^{\ell}T^\ast\Sigma$ is the $\ell$-th symmetric tensor power of 
$T^{\ast}\Sigma$. The map $\mathrm{incl}$ is defined as follows: 
for $x\in \Sigma$, we fix a local chart $U$ around $x$ where $F$ is trivialized 
as $F_{x}$; then one 
element $\xi$ in $(S^{\ell}T^\ast\Sigma\otimes F)_{x}$ can be constructed 
as $df_{1}\odot df_{2}\odot\cdots \odot df_{\ell}\otimes v$, where $\odot$ 
denotes the symmetric tensor product, $v\in 
F_{x}$ and $f_{1}, \ldots, f_{\ell}$ are smooth functions on $U$ which 
vanish at $x$. Then we define $\mathrm{incl}(\xi):= 
j^{\ell}_{x}(f_{1}f_{2}\cdots f_{\ell}\otimes v)$. As a consequence, we 
have the identification of the vector bundles over $\Sigma$ as follows,
\begin{equation}
	S^{\ell}T^\ast\Sigma\otimes F \cong \faktor{J^{\ell}(F)}{J^{\ell-1}(F)}.
	\label{eq:3.2.2}
\end{equation}

We equip the vector bundle $S^{\ell}T^\ast\Sigma\otimes F$ with the 
metric induced by $g^{T\Sigma}$ and $h^{F}$. For $s\in \cG_{x}(F)$, 
let $j^\ell_x s / j^{\ell-1}_x s\in (S^{\ell}T^\ast\Sigma\otimes F)_{x}$ be 
the unique element determined by isomorphism \eqref{eq:3.2.2}, 
and let $|j^\ell_x s / j^{\ell-1}_x s |$ denote the corresponding norm. 
For $x\in\Sigma$, let $(Z_{1},Z_{2})\in\R^{2}\cong T_{x}\Sigma$ 
denote the normal (geodesic) coordinate centred at $x$. Then for any 
germ $s\in \cG_{x}(F)$, we have
\begin{equation}
	\abs{j^\ell_x s/ j^{\ell-1}_x s}^{2} := \sum_{{\scriptscriptstyle 
	\underset{\abs{\alpha} = \ell}{\alpha \in \N^2}}} 
	\frac{1}{\alpha!}\abs{ \frac{\partial^{|\alpha|} 
	s}{\partial Z^{\alpha}} (0) }_{h^{F}_{x}}^{2}.
\end{equation}
This way, we can define a norm on $J^{\ell}(F)$ as follows, for 
$s\in\cG_{x}(F)$,
\begin{equation}
	\abs{j^\ell_x s }^{2} := \sum_{k=0}^{\ell} \|j^k_x s/ j^{k-1}_x s\|^{2},
\end{equation}
where $|j^0_x s/ j^{-1}_x s|:=|s(x)|_{h^{F}}.$

\printbibliography[title=References] 

\end{document}